\documentclass[11pt,a4paper]{amsart}
\usepackage[T1]{fontenc}
\usepackage{lmodern,amsmath,amssymb,amsthm,mathtools}
\usepackage[margin=27mm,headheight=14pt]{geometry}
\usepackage{microtype,booktabs,longtable,array,enumitem,fancyhdr,float,needspace}
\usepackage{tikz,tikz-cd}
\usetikzlibrary{arrows.meta,calc}
\usepackage[colorlinks=true,linkcolor=blue!45!black,citecolor=blue!45!black,urlcolor=blue!45!black]{hyperref}
\hypersetup{pdftitle={Hausdorff Dimension of Escaping Sets of Transcendental Entire Functions in the Speiser Class},pdfsubject={Hausdorff dimension of escaping sets in the Speiser class}}
\newtheorem{theorem}{Theorem}[section]
\newtheorem{proposition}[theorem]{Proposition}
\newtheorem{lemma}[theorem]{Lemma}

\theoremstyle{definition}

\numberwithin{equation}{section}
\DeclareMathOperator{\diam}{diam}
\DeclareMathOperator{\dist}{dist}
\DeclareMathOperator{\area}{area}
\DeclareMathOperator{\Rea}{Re}
\DeclareMathOperator{\Ima}{Im}

\newcommand{\HHm}{\mathbb H_l}
\newcommand{\dimH}{\dim_{\mathrm H}}
\newcommand{\calS}{\mathcal S}
\newcommand{\calB}{\mathcal B}
\newcommand{\C}{\mathbb C}
\newcommand{\R}{\mathbb R}
\newcommand{\Z}{\mathbb Z}

\newcommand{\HH}{\mathbb H_r}

\setlist[enumerate]{label=(\roman*),leftmargin=*,itemsep=3pt}
\title{Hausdorff Dimension of Escaping Sets of\\Transcendental Entire Functions in the Speiser Class}

\author{YUMING FU}
\address{School of Mathematical Sciences, Shenzhen University, Shenzhen, 518060, P. R. China}
\email{yumingfuxy@szu.edu.cn}

\author{WEIWEI SU}
\address{School of Mathematical Sciences, Nanjing University, Nanjing, 210093, P. R. China}
\email{suweiwei2024@163.com}

\subjclass[2020]{Primary: 37F10; Secondary: 37F20, 37F50}
\keywords{Transcendental entire functions; escaping sets; Hausdorff dimension; Speiser class; quasiconformal folding.}
\date{\today}

\begin{document}
\maketitle
\begin{abstract}
For every \(d\in[1,2]\), we construct a transcendental entire function \(f\) whose finite singular value is contained in \(\{-1,0,1\}\) and whose escaping set has Hausdorff dimension \(d\). For \(d<2\), the construction combines spiral corridors with quasiconformal folding. A narrowing Archimedean spiral gives dimension one, while a logarithmic spiral with adjustable width gives intermediate dimensions. The endpoint \(d=2\) is realized by the exponential map. This gives an affirmative answer to a question of Bishop posed in 2018. .
\end{abstract}
\maketitle


\section{Introduction}\label{sec:introduction}

Let $f:\C\to\C$ be a transcendental entire function, and let $f^n$ denote its $n$th iterate. Its \emph{escaping set} is
\begin{equation}\label{eq:escaping-set}
 I(f)=\{z\in\C:|f^n(z)|\to\infty\text{ as }n\to\infty\}.
\end{equation}
A \emph{critical value} of $f$ is a value $f(c)$ with $f'(c)=0$. A \emph{finite asymptotic value} of $f$ is a point $a\in\C$ for which there is a curve $\gamma:[0,\infty)\to\C$ such that $|\gamma(t)|\to\infty$ and $f(\gamma(t))\to a$ as $t\to\infty$. The set of singular values of $f$ is
\[
 S(f)=\overline{\{\text{critical values and finite asymptotic values of }f\}},
\]
where the closure is taken in $\C$.

The Speiser class $\calS$ consists of the transcendental entire functions for which $S(f)$ is finite, and the Eremenko--Lyubich class $\calB$ consists of those for which $S(f)$ is bounded. Thus $\calS\subset\calB$. Write $J(f)$ for the Julia set, the set of points having no neighborhood on which the iterates of $f$ form a normal family. For $f\in\calB$, one has $I(f)\subset J(f)$, and the dynamics near infinity can be studied in logarithmic coordinates (see \cite[Section~2 and Theorem~1]{EremenkoLyubich}). Writing $\dim_H$ for Hausdorff dimension, we ask which values $\dim_H I(f)$ can take for $f\in\calS$.

Every transcendental entire function satisfies $1\le\dimH I(f)\le2$ (see \cite[Theorem~9.1]{BergweilerRempe}). The upper endpoint is already attained by exponential maps: McMullen's construction yields \(\dim_H I(\lambda e^z)=2\) for \(\lambda\ne0\) (see \cite[proof of Theorem~1.2]{McMullen}). The growth condition imposes a more general restriction: if $f\in\calB$ has finite order, that is,
\[
 \rho(f)=\limsup_{r\to\infty}\frac{\log\log M(r,f)}{\log r}<\infty
 \quad\text{with}\quad M(r,f)=\max_{|z|=r}|f(z)|,
\]
then \(\dim_{\mathrm H}I(f)=2\). This conclusion follows from the independent proofs of Barański \cite{Baranski2008} and Schubert \cite{Schubert2007}. It also follows directly from the growth theorem of Bergweiler, Karpińska and Stallard \cite[Theorem~1.2]{BergweilerKarpinskaStallard} with \(q=1\). Consequently, every example in $\calB$ with escaping dimension below two has infinite order.

Without the finite-order restriction, every $d\in[1,2]$ is realized by $I(f)$ for some $f\in\calB$ (see \cite[Corollary~1.2]{RempeStallard}). In particular, Rempe and Stallard constructed a function in $\calB$ with escaping dimension one. This differs from the Julia-set problem: Stallard proved that $\dimH J(f)>1$ for every $f\in\calB$ (see \cite[Theorem~2]{Stallard1996}). Since $I(f)\subset J(f)$ in this class, the dimension of the escaping set can be strictly smaller than that of the Julia set.

Requiring the set of singular values to be finite imposes additional constraints on the construction. Albrecht and Bishop constructed functions in $\calS$ with three singular values and Julia-set dimension arbitrarily close to one \cite[Theorem~1.1]{AlbrechtBishop}. Their result also gives escaping dimensions arbitrarily close to one. The question whether every \(d\in[1,2]\) can be realized as \(\dim_H I(f)\) for an entire function \(f\in\mathcal S\) is recorded in \cite[Question 9.3]{BergweilerRempe}, where it is attributed to Albrecht and Bishop. The endpoint case \(d=1\) was explicitly posed earlier by Bishop in his 2018 research proposal, Question 16.

Let $f_{\mathrm{mer}}$ be a transcendental meromorphic function of finite order $\rho$ whose finite singular values form a bounded set. Suppose that $\infty$ is not an asymptotic value of $f_{\mathrm{mer}}$ and that all but finitely many poles of $f_{\mathrm{mer}}$ have multiplicity at most $M$, where $M$ is a positive integer. Bergweiler and Kotus \cite[Theorems~1.1 and~1.2]{BergweilerKotus} proved that
\[
 \dimH I(f_{\mathrm{mer}})\le \frac{2M\rho}{2+M\rho}.
\]
For every $0<\rho<\infty$ and every positive integer $M$, they constructed examples attaining equality. Under the same hypotheses on the singular values and poles, Mayer and Urba\'nski \cite[Theorem~1.3]{MayerUrbanski} expressed $\dimH I(f_{\mathrm{mer}})$ as the critical exponent of a series determined by the pole locations, their multiplicities, and the leading Laurent coefficients at the poles. Their formula also holds for functions of infinite order.

For meromorphic functions with finitely many singular values, Aspenberg and Cui used quasiconformal surgery to realize every escaping dimension in $[0,2]$ \cite[Theorem~1.1]{AspenbergCui2021}. They later showed that four singular values suffice (see \cite[Theorem~1.1]{AspenbergCui2023}), and Bergweiler and Cui showed that three suffice (see \cite[Corollary~1.3]{BergweilerCui}). These counts are taken on the Riemann sphere. More generally, for every $0\le d_I\le d_J\le2$ with $d_J>0$, Bergweiler and Cui constructed a meromorphic function $f_{\mathrm{mer}}$ with exactly three singular values such that $\dimH I(f_{\mathrm{mer}})=d_I$ and $\dimH J(f_{\mathrm{mer}})=d_J$ (see \cite[Theorem~1.4]{BergweilerCui}). 

The following theorem shows that every $d\in[1,2]$ occurs as the escaping dimension of a transcendental entire function with at most three finite singular values.

\begin{theorem}\label{thm:construction}
For every $d\in[1,2]$, there exists a transcendental entire function
$f_d$ such that
\[
 S(f_d)\subset\{-1,0,1\}\quad\text{and}\quad \dimH I(f_d)=d.
\]
For $d=2$, one may take $f_2(z)=e^z$.
\end{theorem}

The case $d=1$ gives a function in $\calS$ with $\dimH I(f)=1<\dimH J(f)$, by \cite[Theorem~2]{Stallard1996}. For $1<d<2$, put $q=(d-1)^{-1}$. The construction also satisfies
\[
 \log\log M(r,f_d)\le C(\log r)^q,
\]
for \(r\) sufficiently large. Thus these examples attain the growth lower bound $1+1/q$ from \cite[Theorem~1.2]{BergweilerKarpinskaStallard}. Together with the growth estimate above, these examples show that the lower bound $1+1/q$ is sharp in the Speiser class, answering \cite[Question~9.5]{BergweilerRempe} in the affirmative.

The proof uses quasiconformal folding (see \cite{BishopFolding}), in the $L/R$ form of \cite[Theorem~3.4]{AlbrechtBishop} and with the corrected boundary identifications of \cite[Sections~1 and~4]{BishopCorrections}. Related constructions from plane trees include Cui's realization, under a topological uniformness condition, of an infinite tree as $f^{-1}([-1,1])$ up to a homeomorphism of the plane, where $f$ is entire with only two critical values $\pm1$ and no finite asymptotic values \cite[Theorem~1.1]{CuiTrees}.

We begin with a simply connected spiral corridor $\Omega$. Its outside domain $L$ has a conformal coordinate in a left half-plane. After inserting slits, the inside domain $R$ has a conformal coordinate in a right half-plane. The exponential models on these domains are modified near their common graph to obtain a quasiregular map $g$. A quasiconformal homeomorphism $\varphi$ solving the Beltrami equation of $g$ then gives the entire function 
\begin{equation}\label{eq:folding-factorization} 
f=g\circ\varphi^{-1}.
\end{equation}

Logarithmic spirals also occur in the quasiconformal surgery of Aspenberg and Cui (see \cite[Section~3]{AspenbergCui2023}), where elliptic functions are used to construct meromorphic functions with prescribed escaping dimension. Here the spiral geometry controls the widths of logarithmic tracts of entire functions and the estimates for their inverse branches. The two spiral geometries lead to different covering estimates. For $d=1$, the corridor follows the Archimedean spiral $te^{it}$ and has width scale $\eta(t)=e^{-t}$. Its rapid narrowing gives finite $s$-cost covers for every $s>1$. For $1<d<2$, a logarithmic spiral has width scale
\[
 \eta(r)\asymp\frac{r}{(\log r)^p},\qquad
 p=q-1=(d-1)^{-1}-1.
\]
This width law also occurs in the examples in $\calB$ described in \cite[Proposition~4.1]{RempeStallard}. Here it is incorporated into a folding construction with a fixed finite singular set. The widening corridor is covered at its local width scale; one inverse image of this cover has finite $s$-cost for every $s>d$. The endpoint $d=2$ is supplied by the exponential map.

The common analytic step is a simultaneous matching of boundary scales. The outer coordinate prescribes intervals of length \(2\pi\), whereas the same physical edges may have very different conformal lengths when measured from inside the corridor. We first normalize the entrances between successive slits, and then construct one graded subdivision that works on both sides of each slit. Separate estimates, away from and near the slit tips, control this subdivision and the geometry of nearby edges. This order of the argument establishes the entrance estimates before using the boundary-length comparisons they support.

To pass from the geometric model to an entire function, we control both conformal boundary lengths and the folding correction. A metric harmonic-measure estimate places the support of the Beltrami coefficient of the change of conformal coordinates in a fixed vertical strip. Distortion estimates and quasisymmetry then give uniform conformal side-length bounds after straightening in the narrowing model. A uniform logarithmic-area estimate treats individual edges, neighboring slits, and intermediate physical scales. Together with \cite[Lemma~15.1]{AlbrechtBishop}, it makes \(\varphi\) bilipschitz on the region relevant to large-value preimages. The bounded interpolation in the folding construction keeps all modified values uniformly bounded. It therefore identifies the full preimage of the exterior of a sufficiently large disk and gives all inverse branches there.

For the upper bounds, Theorem~\ref{prop:covering-principle} separates the finite initial covering cost from subsequent contraction under inverse branches. In the logarithmic model, grouping boundary intervals at the scale of a complete slit and its following trunk edge gives the required inverse-derivative sums. The resulting estimate for the first covering cost is summable whenever \(s>d\), after allowing an arbitrarily small exponent loss. Each function is fixed before the covering exponent \(s\) is chosen. The lower bounds follow from the universal estimate \(\dim_H I(f)\ge1\) and, for intermediate dimensions, the growth theorem \cite[Theorem~1.2]{BergweilerKarpinskaStallard}.

The paper is organized as follows. Section~\ref{sec:background} collects the definitions, distortion estimates and covering principle used in the proofs. Section~\ref{sec:models} constructs the two spiral models and their boundary meshes. Section~\ref{sec:slits} establishes the common slit estimates and simultaneous subdivision. Section~\ref{sec:graph-geometry} proves the geometric, harmonic-measure and logarithmic-area estimates, and Section~\ref{sec:folded-tract} constructs the folded entire functions and their large-value inverse branches. Subsection~\ref{sec:dimone} treats the narrowing model. Subsection~\ref{sec:intermediate} proves the boundary-block estimates for the logarithmic model and uses them to obtain the upper bound and verifies the growth condition for the lower bound. Subsection~\ref{sec:scope} treats the endpoint \(d=2\).

\section{Preliminaries}\label{sec:background}

\subsection{Metric notation, Hausdorff dimension and logarithmic tracts}
\label{subsec:definitions}

All distances and diameters are Euclidean unless stated otherwise. Write $D(z,r)=\{w:|w-z|<r\}$ and let $dA$ denote planar Lebesgue area measure. For two positive numbers, $u_n\asymp v_n$ means $cv_n\le u_n\le Cv_n$ with constants independent of $n$, whereas $u_n\sim v_n$ means $u_n/v_n\to1$. We write $u=O(v)$ if $|u|\le C|v|$, and $u=o(v)$ if $u/v\to0$.

For a set $E\subset\C$ and real numbers $s,\varepsilon>0$, define
\[
 \mathcal H^s_\varepsilon(E)=
 \inf\left\{\sum_j(\diam E_j)^s:
 E\subset\bigcup_j E_j,\quad \diam E_j\le\varepsilon\right\},
\]
where covers are countable. We define the Hausdorff measure and its associated dimension by
\[
 \mathcal H^s(E)=\lim_{\varepsilon\to 0^+}
 \mathcal H^s_\varepsilon(E)\quad\text{and}\quad
 \dimH E=\inf\{s>0:\mathcal H^s(E)=0\};
\]
see \cite[Sections~2.1--2.2]{Falconer}. For a cover \(\{E_j\}\), its \(s\)-cost is defined as \(\sum_j(\operatorname{diam}E_j)^s\) (see \cite[Section~1]{Kolossvary}), and its mesh is defined as \(\sup_j\operatorname{diam}E_j\) (see \cite[Definition~2.5]{Colvin}).

A map $F:E\to\C$ is $L$-bilipschitz if
\[
 L^{-1}|z-w|\le |F(z)-F(w)|\le L|z-w|\text{ for }z,w\in E.
\]
Lipschitz maps do not increase Hausdorff dimension, and bilipschitz maps preserve it (see \cite[Corollary~2.4]{Falconer}).

\begin{lemma}\label{lem:analytic-preimages}
If $F$ is a nonconstant entire function and $E\subset\C$, then $\dimH F^{-1}(E)\le\dimH E$.
\end{lemma}
\begin{proof}
Since $F'$ is not identically zero, its zeros are isolated, so $F$ has at most countably many critical points. Near each noncritical point, $F$ is biholomorphic onto a disk. On a smaller concentric disk, the inverse branch has bounded derivative and is therefore Lipschitz. These neighborhoods admit a countable subcover.

Thus the noncritical part of $F^{-1}(E)$ is a countable union of Lipschitz images of subsets of $E$, while the set of critical points has Hausdorff dimension zero. The Lipschitz inequality and countable stability of Hausdorff dimension give $\dimH F^{-1}(E)\le\dimH E$ as in \cite[Section~2.2]{Falconer}.
\end{proof}

For $R>0$, let $A_R=\{w:|w|>R\}$. A logarithmic tract over $A_R$ is an unbounded simply connected component $D$ of $f^{-1}(A_R)$ for which \(f|_D=\exp\circ\mathcal T,\) where \(\mathcal T:D\to\{\Rea\zeta>\log R\}\) is conformal. Thus $f:D\to A_R$ is a universal covering (see \cite[Section~2]{EremenkoLyubich}). If $\mathcal P=\mathcal T^{-1}$,
the inverse branches with image in $D$, on a simply connected subdomain of $A_R$, are
\[
 h_k(w)=\mathcal P(\log w+2\pi i k)\text{ for } k\in\Z,
\]
with a fixed holomorphic branch of $\log$. Changing this branch only relabels $k$. 

\subsection{Conformal coordinates and local estimates}
\label{subsec:conformal-estimates}

Throughout the construction, set
\[
S=\{x+iy:|y|<\frac\pi2\},\quad
 \HH=\{\Rea\zeta>0\},\quad\text{and}\quad \HHm=\{\Rea\zeta<0\}.
\]
Then $\exp:S\to\HH$ is conformal. For a simply connected proper domain $D$, write $\omega_D(z,E)$ for harmonic measure, with boundary sets interpreted in the prime-end boundary when necessary. It is obtained from Poisson measure by a conformal map and is conformally invariant. The two approaches to the interior of a slit give distinct boundary sides, even when they belong to the same domain (see \cite[Chapters~I and~III]{GarnettMarshall}). In particular,
\begin{equation}\label{eq:poisson-half-plane}
 \omega_{\HH}(x+iy,i[a,b])=
 \frac1\pi\int_a^b\frac{x}{x^2+(t-y)^2}\,dt,\text{ for }x>0.
\end{equation}
Thus $i[y-x,y+x]$ has harmonic measure $1/2$ at $x+iy$.

Let $\rho_D$ denote hyperbolic distance with curvature $-4$. The densities in the unit disk and in $\HH$ are $1/(1-|z|^2)$ and $1/(2\Rea z)$, respectively. Inclusion decreases hyperbolic distance. Applying the map \(\exp:S\to\HH\) to the points \(u_0+iv_0\) and \(u+iv\) and using the hyperbolic distance formula in \(\HH\) gives
\begin{equation}\label{eq:strip-hyperbolic-distance}
 \cosh\bigl(2\rho_S(u_0+iv_0,u+iv)\bigr)=
 \frac{\cosh(u-u_0)-\sin v\sin v_0}{\cos v\cos v_0}.
\end{equation}
Fix \(z_0=u_0+iv_0\in S\), and let \(z=u+iv\) approach either horizontal boundary of \(S\) while \(u\) remains in a fixed bounded interval. Then $\delta=\pi/2-|v|\to0$ and \(\cos v=\sin\delta\sim\delta\). The numerator in~\eqref{eq:strip-hyperbolic-distance} is bounded above because \(u\) is bounded, and it is bounded below by \(1-|\sin v_0|>0\). It follows that
\[
\rho_S=\tfrac12\log(1/\delta)+O(1).
\]
The same estimate holds near a fixed analytic boundary point in a conformal chart with nonzero boundary derivative.

For a conformal map $\Psi:\HH\to D$, Koebe's quarter theorem gives
\begin{equation}\label{eq:koebe-distance}
 \frac{X}{2}|\Psi'(X+iy)|
 \le\dist(\Psi(X+iy),\partial D)
 \le2X|\Psi'(X+iy)|\quad \text{for}\quad X>0,
\end{equation}
see \cite[Chapter~I, Section~1 and Theorems~4.3 and~4.5]{GarnettMarshall}. The Beurling projection theorem also gives, for a boundary side \(\mathfrak s\) of a geometric edge \(e\) in the simply connected domains considered here,
\begin{equation}\label{eq:beurling-side}
 \omega_D(z,\mathfrak s)
 \le\frac4\pi\arctan\sqrt{\frac{\diam e}{\dist(z,e)}}
 \le\frac4\pi\sqrt{\frac{\diam e}{\dist(z,e)}};
\end{equation}
see \cite[proof of Lemma~2.1]{BishopFolding}.

\subsection{Quasiconformal estimates and dimension lower bounds}
\label{subsec:qc-lower-bounds}

A planar quasiconformal homeomorphism is an orientation-preserving homeomorphism $F\in W^{1,2}_{\mathrm{loc}}$ satisfying
\(F_{\bar z}=\mu_F F_z\) almost everywhere, and \(\|\mu_F\|_\infty<1.\) Its maximal dilatation is $K=(1+\|\mu_F\|_\infty)/(1-\|\mu_F\|_\infty)$. A nonconstant quasiregular map satisfies the same analytic distortion condition without being required to be injective.

Every global $K$-quasiconformal map $F:\C\to\C$ is quasisymmetric \cite[Theorem~11.14]{Heinonen}. More precisely, there is an increasing homeomorphism $\vartheta_K:[0,\infty)\to[0,\infty)$, depending only on $K$, such that
\[
 \frac{|F(z)-F(a)|}{|F(z)-F(b)|}
 \le \vartheta_K\!\left(\frac{|z-a|}{|z-b|}\right)
\]
for all distinct $z,a,b\in\C$.

\begin{lemma}[Normalized logarithmic-area distortion, {\cite[Lemma~15.1]{AlbrechtBishop}}]
\label{lem:normalized-distortion}
Let $A\subset\C$ be measurable, let $B\subset\C\setminus A$, and let $\varphi:\C\to\C$ be $K$-quasiconformal, with $\mu_\varphi=0$ almost everywhere on $\C\setminus A$ and $\varphi(z)/z\to1$ as $|z|\to\infty$. Suppose that
\begin{equation*}
 \sup_{w\in B}\int_A\frac{dA(z)}{|z-w|^2}\le C<\infty.
\end{equation*}
Then there is $L=L(K,C)\ge1$ such that
\[
 L^{-1}|w_1-w_2|\le|\varphi(w_1)-\varphi(w_2)|
 \le L|w_1-w_2|,
\]
for any \(w_1,w_2\in B\).
\end{lemma}

The Teichm\"uller--Wittich--Belinskii theorem states that a global quasiconformal map $\varphi$ satisfying
\[
 \int_{|z|>1}\frac{|\mu_\varphi(z)|}{|z|^2}\,dA(z)<\infty
\]
has $\varphi(z)/z\to c\ne0$ \cite[Chapter~V, Theorem~6.1]{LehtoVirtanen}. Division by $c$ gives the normalization in Lemma~\ref{lem:normalized-distortion}. In particular, finite area of the Beltrami support $E=\operatorname{supp}\mu_\varphi$ guarantees this integrability. Indeed, since $|\mu_\varphi|\le k<1$ and $|z|^{-2}\le1$ for $|z|>1$,
\[
 \int_{|z|>1}\frac{|\mu_\varphi(z)|}{|z|^2}\,dA(z)
 \le k\,\area(E)<\infty.
\]

Every transcendental entire function satisfies $\dimH I(f)\ge1$ (see \cite[Theorem~9.1]{BergweilerRempe}). Let $D$ be a logarithmic tract of $f$ and let $q\ge1$. For every $\varepsilon>0$, if
\begin{equation}\label{eq:tract-growth}
 |f(z)|\le\exp\!\left(\exp\!\left((\log|z|)^{q+\varepsilon}\right)\right)
\end{equation}
for \(z\in D\), and  \(|z|\) sufficiently large,  then \cite[Theorem~1.2]{BergweilerKarpinskaStallard} gives
\begin{equation}\label{eq:tract-lower-bound}
 \begin{gathered}
 \dimH I(f,D)\ge1+\frac1q\quad\text{and}\\
 I(f,D)=\{z\in D:f^n(z)\in D\ (n\ge0),\ |f^n(z)|\to\infty\}.
 \end{gathered}
\end{equation}
In particular, the same lower bound holds for $\dimH I(f)$.

\subsection{A covering principle for escaping orbits}
\label{subsec:covering-principle}

Let $f\in\calB$, and choose $R_0>0$ such that $S(f)\subset D(0,R_0)$. For $s>0$ and $|w|>R_0$, define
\begin{equation*}
 P_s(w)=\sum_{f(z)=w}|f'(z)|^{-s},
\end{equation*}
and every point $z\in f^{-1}(w)$ satisfies $f'(z)\ne0$. On a simply connected domain disjoint from $S(f)$, the sum equals $\sum_h|h'(w)|^s$, where $h$ ranges over all inverse branches of $f$. Such branches exist by the covering property outside $S(f)$
(see \cite[Section~2]{EremenkoLyubich}). For $R>R_0$, put
\begin{equation}\label{eq:exterior-orbits}
 X_R=\{z\in\C:|f^n(z)|>R\text{ for every }n\ge0\}.
\end{equation}
Then $f(X_R)\subset X_R$, and every escaping orbit eventually enters $X_R$. This set may also contain nonescaping orbits.

\begin{theorem}[Finite cost followed by contraction]
\label{prop:covering-principle}
Fix $s>0$. Suppose that $P_s(w)\to0$ uniformly as $|w|\to\infty$. Suppose also that, for every sufficiently large $R$, the set $X_R$ admits a countable cover of mesh at most $1/16$ and finite $s$-cost. Then $\dimH I(f)\le s$.
\end{theorem}

\begin{proof}
Set $K_0=48$. Choose $R>R_0+1$ sufficiently large that
\[
 \sup_{|w|>R}P_s(w)\le(2K_0)^{-s}
\]
and $X_R$ admits the cover assumed in the theorem. Discard covering sets that miss $X_R$. For each remaining set of positive diameter, choose a point in its intersection with
$X_R$ and replace the set by the disk centered there with radius twice its diameter. Since the original mesh is at most $1/16$,
these disks have radii at most $1/8$, and their total $s$-cost is at most $4^s$ times the original cost. Replace the singleton
sets $\{z_j\}$ by disks $D(z_j,2^{-j-2})$. These disks also have radii less than $1/4$ and finite total $s$-cost.

Let $D=D(w,r)$ be one such disk. The disk $D(w,2r)$ avoids $S(f)$, so every inverse branch $h$ on $D$ extends univalently to $D(w,2r)$. Koebe distortion gives
\[
 \diam h(D)\le12|h'(w)|\,\diam D.
\]
For each branch with $h(D)\cap X_R\ne\varnothing$, choose a point of this intersection and cover $h(D)$ by the disk $D_h$ centered there with radius $2\diam h(D)$. We obtain
\[
 \diam D_h\le K_0|h'(w)|\,\diam D\quad\text{and} \quad
 \sum_h(\diam D_h)^s\le K_0^sP_s(w)(\diam D)^s,
\]
where the sum is taken over inverse branches $h$ satisfying $h(D)\cap X_R\ne\varnothing$.

Applying this construction to every disk \(D_j\) gives another cover of $X_R$, since $f(X_R)\subset X_R$. We obtain
\[
\sum_j\sum_h(\diam D_{j,h})^s
\le2^{-s}\sum_j(\diam D_j)^s.
\]
Moreover, since $|h'(w_j)|^s\le P_s(w_j)\le(2K_0)^{-s}$, we have
\[
 \diam D_{j,h}
 \le K_0|h'(w_j)|\,\diam D_j
 \le\frac12\diam D_j,
\]
where $w_j$ is the center of $D_j$. All new centers lie in $X_R$, and the construction can be repeated.

If the initial cost and mesh are $C_s<\infty$ and $\varepsilon_0$, the $n$th cover has cost at most $2^{-ns}C_s$ and mesh at most
$2^{-n}\varepsilon_0$. Hence $\mathcal H^s(X_R)=0$. Finally,
\[
 I(f)\subset\bigcup_{m\ge0}f^{-m}(X_R).
\]
Lemma~\ref{lem:analytic-preimages} and countable stability of Hausdorff dimension give $\dimH I(f)\le s$.
\end{proof}

\section{The two spiral models and their boundary meshes}\label{sec:models}
We construct two spiral models independently, corresponding to $d=1$ and $1<d<2$, respectively. The model for $1<d<2$ is not obtained by varying the parameters of the model for $d=1$. In both cases, we use the notation $\Omega$, $L$, and $\psi$, together with the same notation for the mesh sequences, to denote the corresponding objects.

\subsection{The narrowing spiral}\label{subsec:narrow-model}
For a large starting parameter $T_0>0$, set
\begin{equation}\label{eq:narrow}
 \gamma(t)=te^{it},\quad \eta(t)=e^{-t}\quad\text{and}\quad
 \Omega=\{(t+u)e^{it}:t>T_0,\ |u|<\eta(t)\}.
\end{equation}
Join the boundary curves $\Gamma_\pm(t)=(t\pm\eta(t))e^{it}$ by a smooth finite cap. For sufficiently large $T_0$, distinct turns of the spiral corridor do not overlap. Its width is exponentially small, whereas successive turns are separated by approximately $2\pi$. Both $\Omega$ and $L=\C\setminus\overline\Omega$ are Jordan domains on the sphere with infinity on their boundaries.

Choose a conformal $\tau_L:L\to\HHm$ taking a fixed cap point to $0$ and infinity to infinity. Set $\sigma_L=\log(-\tau_L)$, using the logarithm on the right half-plane. Choose $\sigma_R:\Omega\to S$ taking the cap point to the left end and infinity to the right end, and write $\psi=\sigma_R^{-1}$.
Let $A_t(z)=e^{-it}(z-\Gamma_+(t))$. The image $A_t(L)$ remains a spiral domain; the selected gap is studied only in fixed local windows.

\begin{lemma}[Local coordinate estimates]\label{lem:narrow-scale}
On either analytic tail, the following estimates hold as $t\to\infty$:
\begin{equation}\label{eq:narrow-der}
 \textup{(1)}\quad |\sigma_L'(\Gamma_\pm(t))|\longrightarrow\frac12.
 \qquad
 \textup{(2)}\quad \eta(t)|\sigma_R'(\Gamma_\pm(t))|\longrightarrow\frac\pi2.
\end{equation}
\begin{equation}\label{eq:narrow-coord}
 \textup{(3)}\quad
 \Rea\sigma_L(\Gamma_\pm(t))=\frac{t^2}{4}+o(t^2)\quad\text{and}\quad
 \Rea\sigma_R(\Gamma_\pm(t))=\left(\frac\pi2+o(1)\right)te^t.
\end{equation}
Boundary derivatives refer to local holomorphic continuations across the analytic tails.
\end{lemma}
\begin{proof}
For a holomorphic function $\sigma=u+iv$ at $z=x+iy$, the Cauchy--Riemann equations give
\[
 u_x=v_y,\quad u_y=-v_x,\quad
 \sigma'=v_y+i v_x\quad\text{and}\quad|\sigma'|=|\nabla v|.
\]
Thus derivative estimates reduce to estimates for the imaginary part, which is harmonic. Write $v_L=\Ima\sigma_L$ and $v_R=\Ima\sigma_R$. Boundary orientation gives
\[
 v_L(\Gamma_\pm(t))=\pm\frac\pi2,\quad
 v_R(\Gamma_\pm(t))=\mp\frac\pi2\quad\text{and}\quad
 |v_L|,|v_R|\le\frac\pi2.
\]
Also,
\[
 \Gamma_\pm'(t)=e^{it}\bigl(1\mp\eta(t)+i(t\pm\eta(t))\bigr)
 \quad\text{and}\quad |\Gamma_\pm'(t)|\sim t\text{ as }t\to\infty.
\]

\emph{(1) The outside coordinate.}
Let
\[
 V_t(w)=v_L(A_t^{-1}(w)).
\]
We compare $V_t$ with
\[
 V_0(w)=\frac\pi2-\frac{\Rea w}{2}\quad\text{in}\quad
 G=\{0<\Rea w<2\pi\}.
\]
\emph{Claim.} $V_t\to V_0$ locally uniformly in $G$ as $t\to\infty$.

Fix $R>2$. Direct expansion of $\Gamma_\pm$ gives
\[
 A_t(\Gamma_+(t+s/t))=is+O_R(t^{-1})\quad\text{and}\quad
 A_t(\Gamma_-(t+2\pi+s/(t+2\pi)))=2\pi+is+O_R(t^{-1}),
\]
uniformly for real $|s|\le R+1$ as $t\to\infty$.
Thus, in $|y|\le R$, the selected walls are graphs $x=a_t(y)$ and $x=b_t(y)$ with
\[
 \varepsilon_t(R):=\max_{|y|\le R}\max\{|a_t(y)|,|b_t(y)-2\pi|\}
 \to 0\text{ as }t\to\infty.
\]
These are the outer boundary of one turn and the inner boundary of the next.
Since the spiral turns are radially ordered and disjoint, the region
\[
 Q_t(R)=\{x+iy:|y|<R,\ a_t(y)<x<b_t(y)\}
\]
lies in $A_t(L)$ for large $t$; the finite cap is then outside this window. 
Figure~\ref{fig:narrow-spiral} shows the selected gap and the two local charts used for reflection across its walls. In the figure, the wall charts are
\[
 g_{t,0}(\zeta)
 = A_t\bigl(\Gamma_+(t+\zeta/t)\bigr)\quad\text{and}\quad
g_{t,1}(\zeta)
 = A_t\bigl(\Gamma_-(t+2\pi+\zeta/(t+2\pi))\bigr).
\]
Here $D^+$ and $D^-$ denote the upper and lower halves, respectively, of a small disk centered at $0$ in the $\zeta$-plane.

\begin{figure}[htbp]
  \centering
  \includegraphics[width=0.9\linewidth]{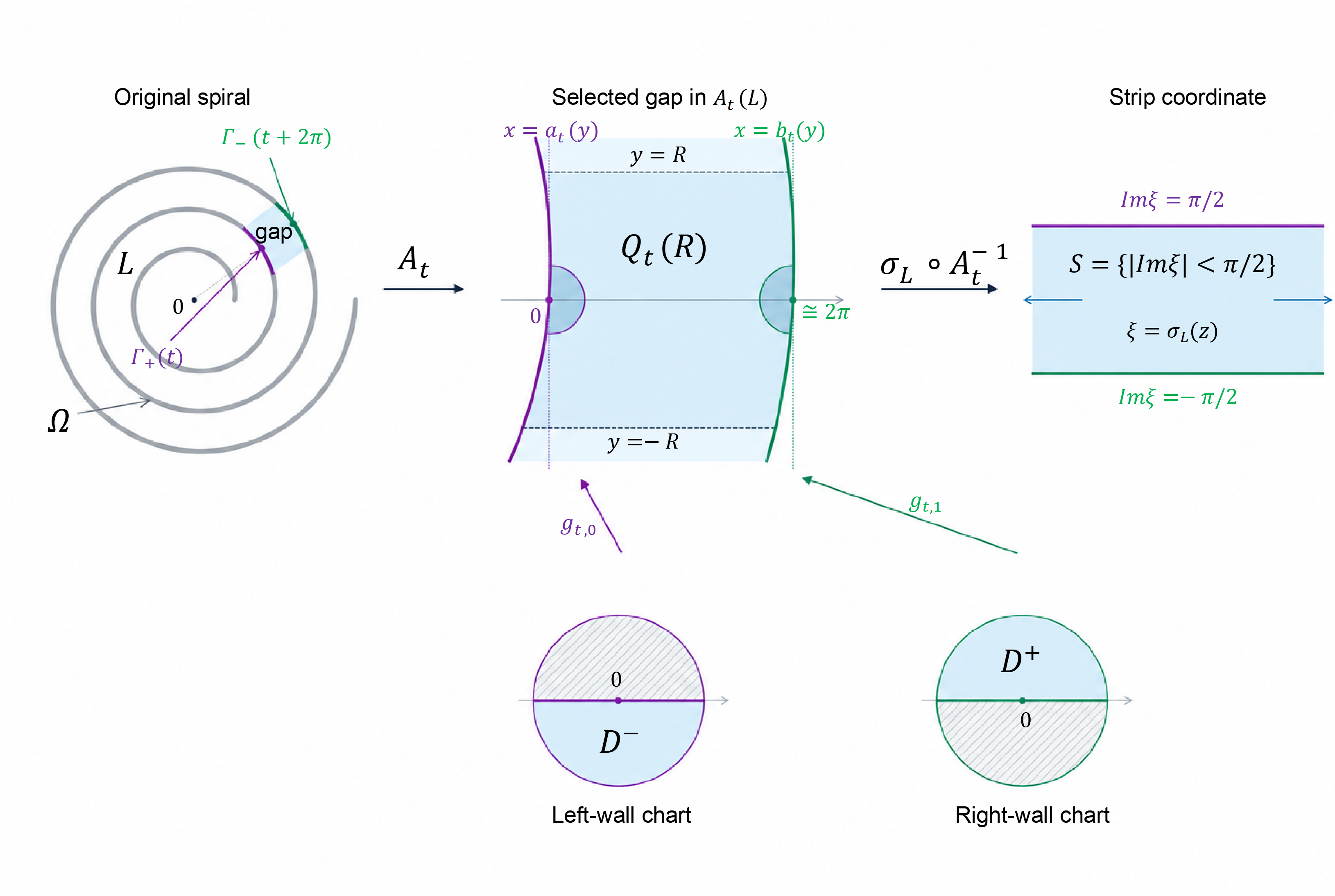}
  \caption[The outside coordinate of the narrowing spiral]{The outside coordinate of the narrowing spiral. The gray spiral represents $\Omega$, and the blue region on the left marks the selected gap in $L$.
    In the middle panel, the dashed lines $y=\pm R$ bound the window $Q_t(R)$ in height. The right panel shows the strip $S=\sigma_L(L)$. In the lower panels, $g_{t,0}$ maps $D^-$ and $g_{t,1}$ maps $D^+$ into the gap. Hatching indicates the halves used for Schwarz reflection.}
  \label{fig:narrow-spiral}
\end{figure}

Set $E_t=V_t-V_0$. On the left and right walls, respectively,
\[
 E_t(a_t(y)+iy)=\frac{a_t(y)}2\quad\text{and}\quad
 E_t(b_t(y)+iy)=\frac{b_t(y)-2\pi}{2}.
\]
Thus $|E_t|\le\varepsilon_t(R)/2$ there. On $y=\pm R$, the bounds for $V_t$ and $V_0$ give $|E_t|\le\pi+\varepsilon_t(R)/2$. When $\varepsilon_t(R)<1$, the harmonic function
\[
 B_R(x,y)=2\cos\bigl((x-\pi)/8\bigr)
             \frac{\cosh(y/8)}{\cosh(R/8)}
\]
is positive on $Q_t(R)$ and at least $1$ on its upper and lower sides. Applying the maximum principle to
$\pm E_t-\varepsilon_t(R)/2-\pi B_R$ yields
\[
 |V_t(x+iy)-V_0(x+iy)|\le\frac{\varepsilon_t(R)}2+\pi B_R(x,y).
\]
For every fixed $0\le Y<R$, this gives
\begin{equation}\label{eq:narrow-comparison}
 \sup_{\substack{w\in\overline{Q_t(R)}\\|\Ima w|\le Y}}
 |V_t(w)-V_0(w)|\le\frac{\varepsilon_t(R)}2+C_Ye^{-R/8}.
\end{equation}
Letting $t\to\infty$ and then $R\to\infty$, we obtain $V_t\to V_0$ locally uniformly in $G$. Interior estimates for harmonic functions and the Cauchy--Riemann equations then give
\[
 e^{it}\sigma_L'(A_t^{-1}(w))
 =\partial_y V_t(w)+i\partial_x V_t(w)
 \longrightarrow-\frac i2
\]
locally uniformly in $G$ as $t\to\infty$; in particular, $|\sigma_L'(A_t^{-1}(w))|\to1/2$ there as $t\to\infty$.

For the boundary derivatives, consider
\[
 \sigma_L\bigl(\Gamma_+(t+\zeta/t)\bigr)-\frac{\pi i}{2}\quad\text{and}\quad
 \sigma_L\bigl(\Gamma_-(t+2\pi+\zeta/(t+2\pi))\bigr)+\frac{\pi i}{2}
\]
for complex $\zeta$ near $0$. These functions extend across the real diameter by Schwarz reflection. The preceding expansions and~\eqref{eq:narrow-comparison} give uniform convergence of their imaginary parts to $\Ima\zeta/2$ on a fixed smaller disk as $t\to\infty$. Interior harmonic estimates, the Cauchy--Riemann equations and the chain rule then give
\begin{equation}\label{eq:narrow-outside-tangent}
 \frac{\sigma_L'(\Gamma_\pm(t))\Gamma_\pm'(t)}{t}
 =e^{it}\sigma_L'(\Gamma_\pm(t))
  \frac{e^{-it}\Gamma_\pm'(t)}{t}
 \longrightarrow\frac12\text{ as }t\to\infty.
\end{equation}
This proves (1).

\emph{(2) The inside coordinate.}
Normalize the thin band by
\[
 B_t(z)=\frac{e^{-it}(z-\gamma(t))}{\eta(t)}\quad\text{and}\quad
 W_t(w)=v_R(B_t^{-1}(w)).
\]
Its local walls tend to $\Rea w=\pm1$ as $t\to\infty$. The left and right boundary values are $\pi/2$ and $-\pi/2$. An argument similar to that in (1) gives
\[
 W_t(w)\longrightarrow-\frac\pi2\Rea w\quad\text{and}\quad
 \eta(t)e^{it}\sigma_R'(B_t^{-1}(w))\longrightarrow-\frac{\pi i}{2},
\]
locally uniformly for $-1<\Rea w<1$ as $t\to\infty$. As in (1), Schwarz reflection gives
\begin{equation}\label{eq:narrow-inside-tangent}
 \frac{\eta(t)}t\sigma_R'(\Gamma_\pm(t))\Gamma_\pm'(t)
 \longrightarrow\frac\pi2\text{ as }t\to\infty.
\end{equation}
Taking moduli proves (2).

\emph{(3) Growth along the tails.}
The signed tangential limits~\eqref{eq:narrow-outside-tangent} and~\eqref{eq:narrow-inside-tangent} imply, as $t\to\infty$,
\[
 \frac{d}{dt}\Rea\sigma_L(\Gamma_\pm(t))=\left(\frac12+o(1)\right)t\quad\text{and}\quad
 \frac{d}{dt}\Rea\sigma_R(\Gamma_\pm(t))
 =\left(\frac\pi2+o(1)\right)te^t.
\]
Integration proves~\eqref{eq:narrow-coord}, since
$\int^t s\,ds=t^2/2+O(1)$ and $\int^t se^s\,ds=(t-1)e^t+O(1)$; the integrated errors are $o(t^2)$ and $o(te^t)$, respectively, as $t\to\infty$.
\end{proof}

\begin{lemma}[Univalent extension and Lipschitz bounds]\label{lem:narrow-lip}
There are constants $r,c,C>0$ such that, for all sufficiently large $t$, $\psi=\sigma_R^{-1}$ has a univalent holomorphic continuation to $D(\sigma_R(\Gamma_\pm(t)),r)$. For any $w_1,w_2$ in this disk,
\begin{equation}\label{eq:narrow-wall-chart}
 c\eta(t)|w_1-w_2|\le |\psi(w_1)-\psi(w_2)|
 \le C\eta(t)|w_1-w_2|.
\end{equation}
\end{lemma}
\begin{proof}
Choose a fixed $0<\rho<1/4$ so that, for all sufficiently large $t$, the parametrizations $\zeta\mapsto\Gamma_\pm(t+\eta(t)\zeta/t)$ map the upper and lower half-disks in $D(0,2\rho)$, respectively, into $\Omega$. On those half-disks, put
\[
 F_{t,\pm}(\zeta)=
 \sigma_R\bigl(\Gamma_\pm(t+\eta(t)\zeta/t)\bigr)
 -\sigma_R(\Gamma_\pm(t)).
\]
The explicit boundary parametrizations satisfy
\[
 B_t\bigl(\Gamma_\pm(t+\eta(t)\zeta/t)\bigr)
 \longrightarrow\pm1+i\zeta
\]
uniformly on fixed disks as $t\to\infty$. The same comparison as in Lemma~\ref{lem:narrow-scale} therefore gives $\Ima F_{t,\pm}(\zeta)\to(\pi/2)\Ima\zeta$ uniformly on smaller half-disks, including the real diameter, as $t\to\infty$. Schwarz reflection extends these functions to $D(0,2\rho)$. Interior harmonic estimates, the Cauchy--Riemann equations and $F_{t,\pm}(0)=0$ give
\[
 F_{t,\pm}(\zeta)\longrightarrow\frac\pi2\zeta\quad\text{and}\quad
 F_{t,\pm}'(\zeta)\longrightarrow\frac\pi2
 \text{ as }t\to\infty,
\]
where both convergences are uniform on $D(0,\rho)$ and its boundary. In particular, for all sufficiently large $t$,
\[
 \left|F_{t,\pm}'(\zeta)-\frac\pi2\right|<\frac\pi8
\quad\text{for}\quad |\zeta|\le\rho.
\]

First, the derivative estimate implies that $F_{t,\pm}$ is univalent on $D(0,\rho)$ and satisfies
\[
 |F_{t,\pm}(\zeta)|\ge\frac{3\pi}{8}|\zeta|\quad\text{and}\quad |F_{t,\pm}'(0)|>\frac{3\pi}{8}.
\]

\emph{Next, we construct the holomorphic continuation of $\psi$ by taking an inverse.}
Koebe's one-quarter covering theorem gives
\[
 D\left(0,\frac{3\pi\rho}{32}\right)
 \subset D\left(0,\frac{\rho|F_{t,\pm}'(0)|}{4}\right)
 \subset F_{t,\pm}(D(0,\rho)).
\]
Set \(r=\pi\rho/16.\) Then $F_{t,\pm}^{-1}$ is holomorphic on a neighborhood of the closed disk $|\xi|\le r$.
To estimate its values, fix $|\xi|\le r$ and put $\zeta=F_{t,\pm}^{-1}(\xi)\in D(0,\rho)$.
Substituting this point into the lower bound for $|F_{t,\pm}(\zeta)|$ gives
\[
 |\xi|=|F_{t,\pm}(\zeta)|
 \ge\frac{3\pi}{8}|\zeta|
 =\frac{3\pi}{8}|F_{t,\pm}^{-1}(\xi)|.
\]
Dividing by $3\pi/8$ and using $r=\pi\rho/16$, we obtain
\[
 |F_{t,\pm}^{-1}(\xi)|
 \le\frac{8}{3\pi}|\xi|
 \le\frac{8}{3\pi}\frac{\pi\rho}{16}
 =\frac\rho6
 \quad\text{for}\quad|\xi|\le r.
\]
Solving the defining relation for $F_{t,\pm}$ now gives the continuation formula
\[
 \psi\bigl(\sigma_R(\Gamma_\pm(t))+\xi\bigr)
 =\Gamma_\pm\left(t+\frac{\eta(t)}tF_{t,\pm}^{-1}(\xi)\right)
 \quad\text{for}\quad |\xi|<r.
\]
The right side is holomorphic on the whole disk.

It can be verified that this formula agrees with the original $\psi$ on the strip side.

\emph{Finally, we prove the univalence and two-sided Lipschitz bounds for the continued $\psi$.}
As $t\to\infty$, the explicit formula for $\Gamma_\pm$ gives $e^{-it}\Gamma_\pm'(t+\eta(t)\zeta/t)/t\to i$ uniformly on fixed small disks, while $F_{t,\pm}'\to\pi/2$. Since $|F_{t,\pm}^{-1}(\xi)|\le\rho/6$ for $|\xi|\le r$, differentiation of the continuation formula gives
\[
 \frac{e^{-it}}{\eta(t)}
 \psi'\bigl(\sigma_R(\Gamma_\pm(t))+\xi\bigr)
 =\frac{e^{-it}\Gamma_\pm'\bigl(t+\eta(t)F_{t,\pm}^{-1}(\xi)/t\bigr)}
         {tF_{t,\pm}'(F_{t,\pm}^{-1}(\xi))}
 \longrightarrow\frac{2i}{\pi}
\]
uniformly for $|\xi|\le r$ as $t\to\infty$. Hence, for large $t$,
\[
 \left|\frac{e^{-it}}{\eta(t)}\psi'(w)-\frac{2i}{\pi}\right|
 <\frac1\pi
 \quad\text{for}\quad w\in D(\sigma_R(\Gamma_\pm(t)),r).
\]
For any $w_1,w_2$ in this disk, integration along the segment joining them yields
\[
 \begin{aligned}
 &\left|\frac{e^{-it}}{\eta(t)}
       \bigl(\psi(w_1)-\psi(w_2)\bigr)-\frac{2i}{\pi}(w_1-w_2)\right|\\
 =&\left|\int_{w_2}^{w_1}
       \left(\frac{e^{-it}}{\eta(t)}\psi'(w)-\frac{2i}{\pi}\right)dw\right|
 \le\frac1\pi|w_1-w_2|.
 \end{aligned}
\]
Applying the reverse triangle inequality gives the lower bound
\[
 \begin{aligned}
 \frac{|\psi(w_1)-\psi(w_2)|}{\eta(t)}
 &\ge\frac2\pi|w_1-w_2|-\left|\frac{e^{-it}}{\eta(t)}
       \bigl(\psi(w_1)-\psi(w_2)\bigr)-\frac{2i}{\pi}(w_1-w_2)\right|\\
 &\ge\left(\frac2\pi-\frac1\pi\right)|w_1-w_2|
 =\frac1\pi|w_1-w_2|.
 \end{aligned}
\]
If $w_1\ne w_2$, this gives
\[
 |\psi(w_1)-\psi(w_2)|\ge\frac{\eta(t)}\pi|w_1-w_2|>0.
\]
Thus the holomorphic continuation of $\psi$ is injective, and hence univalent, on the entire disk $D(\sigma_R(\Gamma_\pm(t)),r)$, not just on its intersection with $S$.

The triangle inequality gives the upper bound
\[
 \begin{aligned}
 \frac{|\psi(w_1)-\psi(w_2)|}{\eta(t)}
 &\le\frac2\pi|w_1-w_2|+\left|\frac{e^{-it}}{\eta(t)}
       \bigl(\psi(w_1)-\psi(w_2)\bigr)-\frac{2i}{\pi}(w_1-w_2)\right|\\
 &\le\left(\frac2\pi+\frac1\pi\right)|w_1-w_2|
 =\frac3\pi|w_1-w_2|.
 \end{aligned}
\]
Multiplying by $\eta(t)$ proves the asserted Lipschitz estimates with $c=1/\pi$ and $C=3/\pi$. All constants are independent of $t$ and of the choice of tail, once $t$ is sufficiently large.
\end{proof}

Fix either tail $\Gamma_\pm$ and choose its boundary vertices $\Gamma_\pm(t_n)$ by
\begin{equation}\label{eq:narrow-vertices}
 \tau_L(\Gamma_\pm(t_n))=\mp 2\pi ni.
\end{equation}
The arc $\Gamma_\pm([t_n,t_{n+1}])$ is the $n$th \emph{trunk edge}. Set
\begin{equation*}
 x_n=\Rea\sigma_R(\Gamma_\pm(t_n)),\quad
 \delta_n=x_{n+1}-x_n,
 \quad 
 \Delta_n=\int_{t_n}^{t_{n+1}}|\Gamma_\pm'(s)|\,ds
 =\int_{x_n}^{x_{n+1}}\left|\psi'\left(x\mp\frac{\pi i}{2}\right)\right|\,dx.
\end{equation*}
Thus $\Delta_n$ is the physical arclength, $\delta_n$ is the inside strip length, and the outside $\tau_L$-length is exactly $2\pi$. Both tails are indexed separately.

\begin{lemma}[Boundary mesh of the narrowing spiral]\label{lem:narrow-mesh}
On either fixed tail of the narrowing spiral, as $n\to\infty$,
\begin{equation}\label{eq:narrow-mesh}
 \Delta_n\sim\frac2n,\quad
 t_n\sim2\sqrt{\log n},\quad
 \delta_n\sim\frac{\pi e^{t_n}}n\quad\text{and}\quad
 \eta(t_n)\delta_n\asymp\Delta_n.
\end{equation}
More precisely, as $n\to\infty$,
\[
 t_{n+1}-t_n\sim\frac{2}{nt_n}\quad\text{and}\quad
 \frac{\eta(t_n)\delta_n}{\Delta_n}\longrightarrow\frac\pi2.
\]
For every fixed exponent $\beta>0$, as $n\to\infty$,
\begin{equation}\label{eq:narrow-small}
 x_n^\beta\delta_n\longrightarrow0\quad\text{and}\quad
 \log(1/\delta_n)=O\bigl((\log x_n)^2\bigr).
\end{equation}
In fact, $\log(1/\delta_n)\sim\tfrac14(\log x_n)^2$ as $n\to\infty$.
\end{lemma}
\begin{proof}
Since $\sigma_L=\log(-\tau_L)$, we have the exact identities
\[
 \Rea\sigma_L(\Gamma_\pm(t_n))=\log(2\pi n),
 \quad
 \Rea\sigma_L(\Gamma_\pm(t_{n+1}))
 -\Rea\sigma_L(\Gamma_\pm(t_n))
 =\log\left(1+\frac1n\right).
\]
Under the logarithmic change of coordinates, edges of equal length in the $\tau_L$-coordinate have lengths tending to zero in the $\sigma_L$-coordinate as $n\to\infty$.

By~\eqref{eq:narrow-coord}, as $n\to\infty$,
\[
 \log(2\pi n)=\frac{t_n^2}{4}+o(t_n^2).
\]
Since $\log(2\pi n)\sim\log n$ as $n\to\infty$, it follows that
\[
 t_n\sim2\sqrt{\log n}\quad\text{and}\quad\frac{t_{n+1}}{t_n}\longrightarrow1\text{ as }n\to\infty.
\]

We next determine the parameter increment along one edge. By~\eqref{eq:narrow-outside-tangent}, we obtain
\[
 \frac{d}{ds}\Rea\sigma_L(\Gamma_\pm(s))
 =\left(\frac12+o(1)\right)s\text{ as }s\to\infty.
\]
As $n\to\infty$, the error tends to zero uniformly for $s\in[t_n,t_{n+1}]$, because every such $s$ is at least $t_n\to\infty$. Hence, as $n\to\infty$,
\[
 \begin{aligned}
 \log\left(1+\frac1n\right)
 &=\int_{t_n}^{t_{n+1}}
       \frac{d}{ds}\Rea\sigma_L(\Gamma_\pm(s))\,ds=\int_{t_n}^{t_{n+1}}\left(\frac12+o(1)\right)s\,ds\\
 &=\left(\frac14+o(1)\right)(t_{n+1}^2-t_n^2).
 \end{aligned}
\]
Therefore $t_{n+1}^2-t_n^2\sim4/n$ as $n\to\infty$, and
\[
 t_{n+1}-t_n
 =\frac{t_{n+1}^2-t_n^2}{t_{n+1}+t_n}
 \sim\frac{4/n}{2t_n}
 =\frac{2}{nt_n}\longrightarrow0\text{ as }n\to\infty.
\]

We now convert to physical and strip lengths. Since $|\Gamma_\pm'(s)|\sim s$ as $s\to\infty$, the physical arclength satisfies, as $n\to\infty$,
\[
 \Delta_n=\int_{t_n}^{t_{n+1}}|\Gamma_\pm'(s)|\,ds=\bigl(1+o(1)\bigr)\int_{t_n}^{t_{n+1}}s\,ds\sim\frac2n.
\]
For the inside coordinate,~\eqref{eq:narrow-inside-tangent} and $\eta(s)=e^{-s}$ give
\[
 \frac{d}{ds}\Rea\sigma_R(\Gamma_\pm(s))
 =\left(\frac\pi2+o(1)\right)se^s\text{ as }s\to\infty.
\]
Thus, as $n\to\infty$,
\[
 \delta_n
 =\int_{t_n}^{t_{n+1}}
       \frac{d}{ds}\Rea\sigma_R(\Gamma_\pm(s))\,ds
 \sim\frac{\pi e^{t_n}}n.
\]
Since $\eta(t_n)=e^{-t_n}$, we have
\[
 \frac{\eta(t_n)\delta_n}{\Delta_n}
 \longrightarrow\frac\pi2\text{ as }n\to\infty,
\]
which proves~\eqref{eq:narrow-mesh}.

By~\eqref{eq:narrow-coord} and the estimate for $\delta_n$, as $n\to\infty$,
\[
 \log x_n=t_n+\log t_n+\log(\pi/2)+o(1),\quad\log\delta_n=t_n-\log n+\log\pi+o(1).
\]
For a fixed exponent $\beta>0$, this gives, as $n\to\infty$,
\[
 \begin{aligned}
 \log(x_n^\beta\delta_n)
 &=-\log n+(\beta+1)t_n+\beta\log t_n+O(1)\\
 &=-\log n+o(\log n)\longrightarrow-\infty,
 \end{aligned}
\]
because $t_n\sim2\sqrt{\log n}$ as $n\to\infty$. Hence $x_n^\beta\delta_n\to0$ as $n\to\infty$.
Also, as $n\to\infty$,
\[
 \log(1/\delta_n)=\log n-t_n+O(1)
 =\frac{t_n^2}{4}+o(t_n^2)\quad\text{and}\quad
 \log x_n\sim t_n.
\]
Consequently, as $n\to\infty$,
\[
 \log(1/\delta_n)\sim\frac14(\log x_n)^2,
\]
which implies the bound in~\eqref{eq:narrow-small}.
\end{proof}

\subsection{The logarithmic spiral}\label{subsec:log-model}
The narrowing spiral of Section~\ref{subsec:narrow-model} is used for dimension one. We now construct the corridor needed for $1<d<2$. Its width will grow with the physical radius $r$, but remain small compared with $r$.

\begin{figure}[H]
  \centering
  \includegraphics[width=0.8\linewidth]{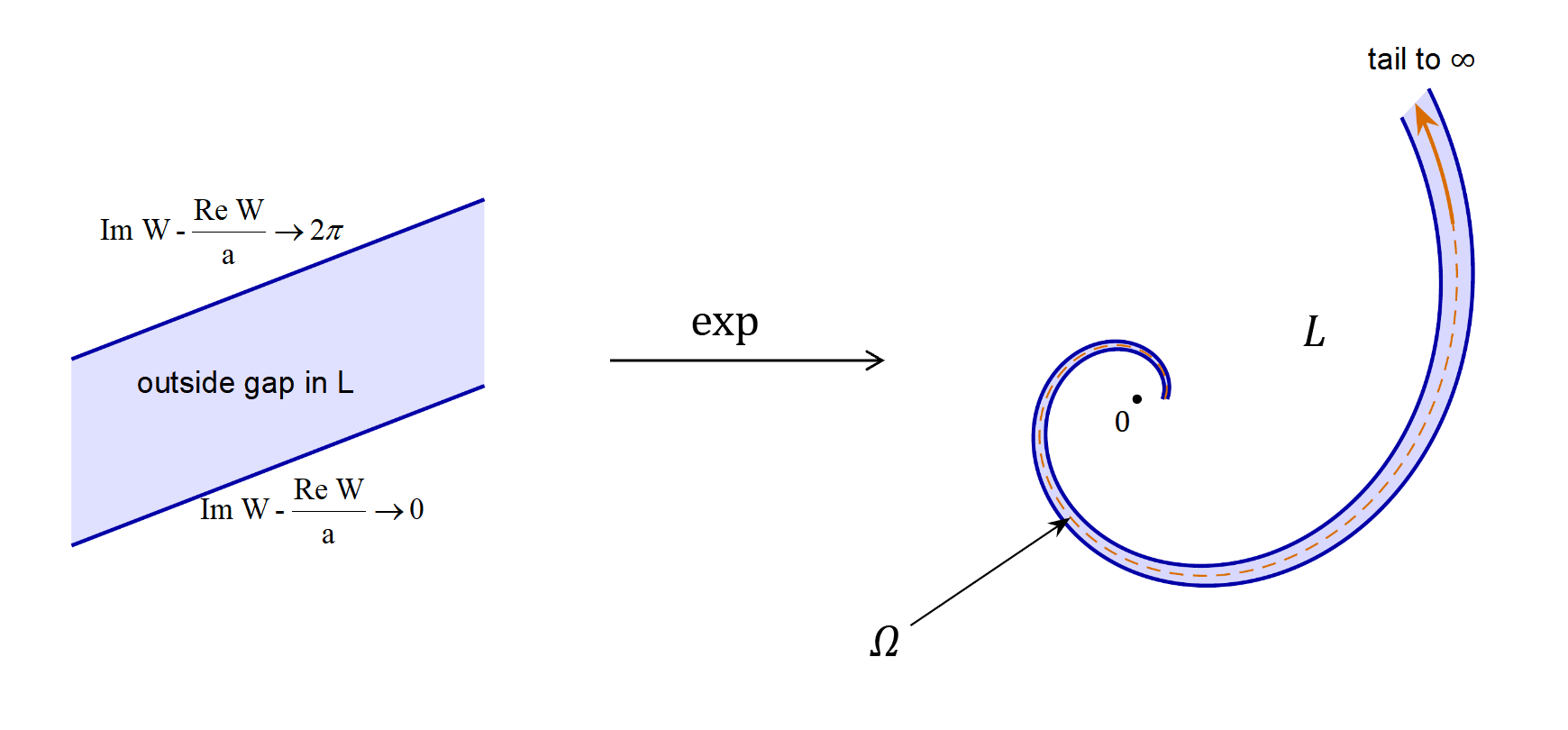}
  \caption[The logarithmic spiral and its outside gap]{%
    The left panel shows the selected outside gap in logarithmic
    coordinates $W=\log z=\log|z|+i\theta$, where $\theta$ is a
    continuous choice of argument along the gap on the spiral tail.
    The indicated limits are taken as $\Rea W\to+\infty$.
    The map $z=e^W$ returns to the physical plane, where the right
    panel shows the spiral corridor $\Omega$ in blue and its
    exterior $L$.}
  \label{fig:log-tail}
\end{figure}

Fix $1<d<2$ and set
\begin{equation}\label{eq:params}
 q=\frac1{d-1}>1\quad\text{and}\quad a=\frac1{\sqrt{2q+3}}.
\end{equation}
Choose $R_0>0$ sufficiently large. With the principal logarithm and power, define the strip parametrization
\begin{equation}\label{eq:log-model}
 \psi(w)=\exp\!\left((a+i)\bigl(R_0+\log(1+e^w)\bigr)^{1/q}\right)\quad\text{and}\quad\Omega=\psi(S).
\end{equation}
The boundary tails are parametrized by $\psi(x\pm\pi i/2)$. As in Section~\ref{subsec:narrow-model}, write $L=\C\setminus\overline\Omega$ and choose $\tau_L:L\to\HHm$ with the cap point sent to $0$ and infinity to infinity; put $\sigma_L=\log(-\tau_L)$. We first estimate the coordinate $\psi$ on the corridor, then the outside coordinate $\tau_L$, and finally the boundary mesh. Figure~\ref{fig:log-tail} shows the spiral geometry.

\begin{lemma}[The inside coordinate]\label{lem:log-scale}
For sufficiently large $R_0$, $\psi$ is conformal on $S$, and $\Omega$ has a smooth cap and two boundary tails. Put
$r(x)=\exp(a(R_0+x)^{1/q})$ and $\eta_{\mathrm{rad}}(r(x)):=\eta(x)=r(x)x^{1/q-1}$,
where $r(x)$ is a reference radius.
Uniformly for $|y|\le\pi/2$ as $x\to\infty$,
\begin{align}
 &\log|\psi(x+iy)|=a(R_0+x)^{1/q}+O(x^{1/q-1}),\label{eq:log-radius}\\
 &\text{and}\quad|\psi'(x+iy)|\asymp\eta(x).\label{eq:log-width}
\end{align}
Consequently,
\begin{equation}\label{eq:radius-width}
 x\asymp\bigl(\log|\psi(x+iy)|\bigr)^q\quad\text{and}\quad\eta(x)\asymp\frac{|\psi(x+iy)|}{\bigl(\log|\psi(x+iy)|\bigr)^{q-1}}.
\end{equation}
For sufficiently large $x$, disks of a fixed sufficiently small radius about $x\pm \frac{\pi i}{2}$ are bilipschitz charts for $\psi$ at scale $\eta(x)$.
\end{lemma}

\begin{proof}
First, we prove that $\psi$ is univalent on $S$. Fix a small $\varepsilon>0$ and consider the convex half-strip
\[
 T=\left\{t:\Rea t>-\varepsilon,\ \bigl|\Ima t\bigr|<\frac\pi2+\varepsilon\right\}.
\]
For $F_0(t)=(a+i)(R_0+t)^{1/q}$, we have
\[
 F_0'(t)=\frac{a+i}q(R_0+t)^{1/q-1}.
\]
Since $\arg(R_0+t)\to0$ uniformly for $t\in T$ as $R_0\to\infty$,
we have $\arg F_0'(t)\to\arg(a+i)\in(0,\pi/2)$ uniformly on $T$ as $R_0\to+\infty$.
Thus, when $R_0$ is sufficiently large, there is a constant $c>0$, independent of $R_0$, such that
\[
 \Rea F_0'(t)\ge c|F_0'(t)|>0\quad\text{and}\quad
 |F_0'(t)|\le M:=\frac{\sqrt{1+a^2}}q\,(R_0-\varepsilon)^{1/q-1}
 \text{ for } t\in T.
\]
As $q>1$, we have $M\to0$ as $R_0\to\infty$. For $t_1,t_2\in T$, convexity gives
\[
 F_0(t_2)-F_0(t_1)=(t_2-t_1)A\quad\text{and}\quad
 A=\int_0^1 F_0'\bigl(t_1+s(t_2-t_1)\bigr)\,ds.
\]
The preceding bounds imply
\[
 \Rea A\ge c\int_0^1|F_0'\bigl(t_1+s(t_2-t_1)\bigr)|\,ds
 \ge c|A|>0\quad\text{for}\quad|A|\le M.
\]
We first prove that $e^{F_0}$ is univalent on $T$.
Suppose that $e^{F_0(t_1)}=e^{F_0(t_2)}$. Then
$(t_2-t_1)A=2\pi i k$ for some $k\in\Z$. If $k\ne0$, it follows that
\[
 \pi+2\varepsilon> |\Ima(t_2-t_1)|
 =2\pi|k|\frac{\Rea A}{|A|^2}\ge\frac{2\pi c|k|}{|A|}
 \ge\frac{2\pi c}{M}.
\]
This is impossible for large $R_0$. Thus $k=0$, and $A\ne0$ gives $t_1=t_2$. Hence $e^{F_0}$ is univalent on $T$. Since $w\mapsto\log(1+e^w)$ maps $S$ injectively into $T$, the composition $\psi$ is univalent on $S$.

For fixed bounded $\xi$ with $|\Ima\xi|<\pi/2+\varepsilon/2$, the identity
\[
 \log(1+e^{x+\xi})
 =x+\xi+\log(1+e^{-x-\xi})=x+\xi+O(e^{-x})
\]
holds as $x\to+\infty$. The enlarged strip has height less than $2\pi$, and its image under $\log(1+e^w)$ stays in $T$; the preceding injectivity argument thus gives a univalent holomorphic extension near either boundary tail.

We next estimate the radius and derivative scales. Hence, as $x\to+\infty$, uniformly for bounded $\xi$,
\[
 (R_0+\log(1+e^{x+\xi}))^{1/q}
 =(R_0+x)^{1/q}+O((R_0+x)^{1/q-1}).
\]
Multiplying by $a+i$ and taking real parts with $\xi=iy$ proves~\eqref{eq:log-radius}. Since $q>1$, it also gives $|\psi(x+iy)|/r(x)\to1$ as $x\to+\infty$, uniformly for $|y|\le\pi/2$.
Differentiating gives
\[
 \psi'(w)=\psi(w)\,\frac{a+i}q
       (R_0+\log(1+e^w))^{1/q-1}\frac{e^w}{1+e^w}.
\]
The expansion above implies, as $x\to+\infty$,
\[
 \frac{\psi(x+\xi)}{e^{(a+i)(R_0+x)^{1/q}}}\to1,
 \quad
 \frac{(R_0+\log(1+e^{x+\xi}))^{1/q-1}}{(R_0+x)^{1/q-1}}\to1\quad\text{and}\quad \frac{e^{x+\xi}}{1+e^{x+\xi}}\to1.
\]
Consequently,
\[
 \frac{q\psi'(x+\xi)}{(a+i)(R_0+x)^{1/q-1}e^{(a+i)(R_0+x)^{1/q}}}\longrightarrow1
 \text{ as }x\to+\infty
\]
uniformly for those bounded $\xi$. Since $(R_0+x)/x\to1$ as $x\to+\infty$, taking moduli with $\xi=iy$ gives
\[
 |\psi'(x+iy)|\sim\frac{\sqrt{1+a^2}}q\,r(x)x^{1/q-1}
       =\frac{\sqrt{1+a^2}}q\eta(x)\text{ as }x\to+\infty,
\]
which proves~\eqref{eq:log-width}. By~\eqref{eq:log-radius}, $\log|\psi(x+iy)|\sim ax^{1/q}$ and $|\psi(x+iy)|\sim r(x)$ as $x\to+\infty$, proving~\eqref{eq:radius-width}.

Fix $0<r_*<\varepsilon/4$. As $x\to+\infty$,
\[
 \frac{\psi'(w)}{\psi'(x)}\longrightarrow1
 \quad\text{uniformly for }w\in D\!\left(x\pm\frac{\pi i}{2},r_*\right).
\]
For sufficiently large $x$ and any two points $w_1,w_2$ in the same disk, integration along the segment joining them gives
\begin{equation}\label{eq:log-wall-chart}
 c\eta(x)|w_1-w_2|\le|\psi(w_1)-\psi(w_2)|
 \le C\eta(x)|w_1-w_2|.
\end{equation}
This proves the boundary chart assertion with constants independent of $x$.

Finally, we prove that $\partial\Omega\cup\{\infty\}$ is a Jordan curve, smooth at every finite point. Writing $\zeta=e^w$, formula~\eqref{eq:log-model} defines a holomorphic function of $\zeta$ on $\HH$ that extends across every finite point of $i\R$ with nonzero derivative. Thus the two boundary tails meet smoothly at $\exp((a+i)R_0^{1/q})$, the image of $\zeta=0$. Moreover, $|\psi(w)|\to\infty$ as $\Rea w\to\infty$, uniformly for $|\Ima w|\le\pi/2$. Since $e^{F_0}$ is injective on $T$, the two boundary curves have no other finite intersection. Hence the boundary is a Jordan curve and is real-analytic at every finite point.
\end{proof}

\begin{lemma}[The outside coordinate]\label{lem:outside-log}
Here $\tau_L'$ denotes complex differentiation with respect to $z$, and $\theta_\pm'$ denotes differentiation with respect to $u$.
Along either tail, as $z=\psi(x\pm\pi i/2)$ and $x\to\infty$,
\begin{equation}\label{eq:outside-log-coordinate}
 \frac{z\tau_L'(z)}{\tau_L(z)}\longrightarrow\frac{1-i/a}{2}\quad\text{and}\quad
 \log|\tau_L(z)|=(q+2)\log|z|+o(\log|z|).
\end{equation}
More precisely, writing the corresponding lifted wall as $W=u+i\theta_\pm(u)$, with $u=\log|z|$, one has, as $u\to+\infty$,
\begin{equation}\label{eq:outside-log-slope}
 \theta_\pm'(u)\longrightarrow\frac1a\quad\text{and}\quad
 \frac{d}{du}\log|\tau_L(e^{u+i\theta_\pm(u)})|
       \longrightarrow q+2.
\end{equation}
\end{lemma}

By Lemma~\ref{lem:log-mesh}, the sum of the $s$-th powers of the edge lengths is finite for $s>1+1/(q+1)$.

\begin{proof}
Write
\[
 \widetilde{\Gamma}_\pm(x)
 =(a+i)\bigl(R_0+\log(1+e^{x\pm\pi i/2})\bigr)^{1/q}\quad\text{and}\quad
 e^{\widetilde{\Gamma}_\pm(x)}=\psi(x\pm\pi i/2).
\]
Expanding and differentiating as in Lemma~\ref{lem:log-scale} gives, as $x\to+\infty$,
\[
 \begin{aligned}
 \widetilde{\Gamma}_\pm'(x)
 &=\frac{(R_0+x)^{1/q-1}}q\bigl(a+i+o(1)\bigr),\\
 \Ima\widetilde{\Gamma}_\pm(x)
 -\frac{\Rea\widetilde{\Gamma}_\pm(x)}a
 &=\pm\frac{\pi(a+a^{-1})}{2q}(R_0+x)^{1/q-1}
   +o(x^{1/q-1}).
 \end{aligned}
\]
Thus $u(x)=\Rea\widetilde{\Gamma}_\pm(x)$ is eventually increasing and, writing $\theta_\pm(u(x))=\Ima\widetilde{\Gamma}_\pm(x)$, the chain rule yields the first limit as $x\to+\infty$ and the second as $u\to+\infty$:
\[
 \theta_\pm'(u(x))
 =\frac{\Ima\widetilde{\Gamma}_\pm'(x)}
        {\Rea\widetilde{\Gamma}_\pm'(x)}\longrightarrow\frac1a\quad\text{and}\quad
 \theta_\pm(u)-\frac ua\longrightarrow0.
\]
The signed expansion also gives $\theta_-(u)<\theta_+(u)$. Hence the selected outside gap is bounded by $\widetilde{\Gamma}_+$ and $\widetilde{\Gamma}_-+2\pi i$, with vertical width tending to $2\pi$ as $u\to+\infty$.

Set $\kappa=(1-i/a)/2$ and normalize by
\[
 \xi=\kappa\bigl(W-(a+i)(R_0+x)^{1/q}\bigr).
\]
Since $\kappa(a+i)=(a+a^{-1})/2>0$ and $\Ima(2\pi i\kappa)=\pi$, the normalized gap converges locally to $0<\Ima\xi<\pi$ as $x\to+\infty$. Let
\[
 V_x(\xi)=\Ima\sigma_L\!\left(
    e^{(a+i)(R_0+x)^{1/q}+\xi/\kappa}\right)\quad\text{and}\quad
 V_0(\xi)=\Ima\xi-\frac\pi2.
\]
As in the proof of Lemma~\ref{lem:narrow-scale}\,(1), $V_x\to V_0$ as $x\to+\infty$, locally and up to the walls.

For the boundary estimate, set
\[
 g_{x,\pm}(\zeta)=\widetilde{\Gamma}_\pm\!\left(
      x+\frac{q\zeta}{(R_0+x)^{1/q-1}}\right)
   -\widetilde{\Gamma}_\pm(x),\ F_{x,\pm}(\zeta)=\sigma_L\!\left(
      e^{\widetilde{\Gamma}_\pm(x)+g_{x,\pm}(\zeta)}\right)
   -\sigma_L\!\left(e^{\widetilde{\Gamma}_\pm(x)}\right).
\]
Here $F_{x,\pm}$ is initially defined on the selected outside side, and the logarithm is continued by $\log(1+e^w)=w+\log(1+e^{-w})$. The expansion in Lemma~\ref{lem:log-scale} and the Schwarz reflection argument in the proof of Lemma~\ref{lem:narrow-scale}\,(1) give, on a fixed smaller disk,
\[
 g_{x,\pm}(\zeta)\longrightarrow(a+i)\zeta\quad\text{and}\quad
 F_{x,\pm}(\zeta)\longrightarrow\kappa(a+i)\zeta
 \quad\text{in }C^1\text{ as }x\to+\infty.
\]
At $z=e^{\widetilde{\Gamma}_\pm(x)}$, differentiation with respect to $\zeta$ gives
\[
 F_{x,\pm}'(0)
 =\frac{z\tau_L'(z)}{\tau_L(z)}
  \frac{q\widetilde{\Gamma}_\pm'(x)}
       {(R_0+x)^{1/q-1}}.
\]
As $x\to+\infty$, the left side tends to $\kappa(a+i)$ and the second factor on the right tends to $a+i\ne0$, proving the boundary logarithmic derivative limit.
Finally, along $z=e^{u+i\theta_\pm(u)}$, as $u\to+\infty$,
\[
 \frac{d}{du}\log|\tau_L(z)|
 =\Rea\left[
   \frac{z\tau_L'(z)}{\tau_L(z)}(1+i\theta_\pm'(u))\right]
 \longrightarrow\Rea\left[\kappa\left(1+\frac ia\right)\right]=\frac{1+a^{-2}}2=q+2.
\]
Integration proves the growth estimate.
\end{proof}

\Needspace{8\baselineskip}
\begin{lemma}[Univalent continuation of the outside coordinate]\label{lem:outside-log-univalent}
There is a constant $c>0$ such that, at every sufficiently far boundary point $z=\psi(x\pm\pi i/2)$, the germ of $\tau_L$ from the selected outside side has a univalent holomorphic continuation to $D(z,c|z|)$.
\end{lemma}
\begin{proof}
Use $g_{x,\pm},F_{x,\pm}$ from the preceding proof, with $z=e^{\widetilde{\Gamma}_\pm(x)}$ and $\kappa=(1-i/a)/2$. The inverse-coordinate argument in the proof of Lemma~\ref{lem:narrow-lip}, applied to their $C^1$ limits, gives inverses $g_{x,\pm}^{-1}$ on a common fixed disk and
\[
 F_{x,\pm}\bigl(g_{x,\pm}^{-1}(h)\bigr)\longrightarrow\kappa h
 \quad\text{in}\quad C^1\text{ as }x\to+\infty.
\]
Thus the holomorphic functions
\[
 T_{x,\pm}(h)
 :=\exp\!\left(F_{x,\pm}\bigl(g_{x,\pm}^{-1}(h)\bigr)\right)
 \longrightarrow e^{\kappa h}\text{ as }x\to+\infty
\]
converge in $C^1$ on a fixed disk; on the selected outside side, $T_{x,\pm}(h)=\tau_L(ze^h)/\tau_L(z)$.
After shrinking a fixed radius $r>0$, for all large $x$ one has
\[
 \sup_{|h|<r}|T_{x,\pm}'(h)-\kappa|<\frac{|\kappa|}{2}.
\]
The segment-integration argument in the proof of Lemma~\ref{lem:narrow-lip} makes $T_{x,\pm}$ univalent on $D(0,r)$. Choose $0<c<\min\{1/2,r/2\}$. For $|w-z|<c|z|$, the branch near $0$ satisfies $|\log(w/z)|<2c<r$. Hence
\[
 w\longmapsto\tau_L(z)\,T_{x,\pm}\bigl(\log(w/z)\bigr)
\]
is the required univalent continuation on $D(z,c|z|)$.
\end{proof}

Choose vertices on each boundary tail by
\begin{equation}\label{eq:log-vertices}
 \tau_L\!\left(\psi\left(x_n\pm\frac{\pi i}{2}\right)\right)=\pm2\pi ni,
\end{equation}
and, for one fixed tail, set
\begin{equation*}
 r_n=\left|\psi\left(x_n\pm\frac{\pi i}{2}\right)\right|,\quad
 \delta_n=x_{n+1}-x_n\quad\text{and}\quad
 \Delta_n=\int_{x_n}^{x_{n+1}}
 \left|\psi'\left(x\pm\frac{\pi i}{2}\right)\right|\,dx.
\end{equation*}
Thus $\Delta_n$ is the physical arclength of an edge and $\delta_n$ its strip-coordinate length, as in Section~\ref{subsec:narrow-model}. The two tails are treated separately.

\Needspace{18\baselineskip}
\begin{lemma}[Boundary mesh of the logarithmic spiral]\label{lem:log-mesh}
In the logarithmic spiral, as $n\to\infty$,
\begin{align}
 &r_n=n^{1/(q+2)+o(1)} \quad\text{and}\quad
 \Delta_n\asymp\frac{r_n}{n}
       =n^{-(q+1)/(q+2)+o(1)},\label{eq:outside-mesh}\\
 &\delta_n\asymp\frac{\Delta_n}{\eta(x_n)}\quad\text{and}\quad
 \log\delta_n=-a(q+2) x_n^{1/q}+o(x_n^{1/q}).
       \label{eq:log-small}
\end{align}
More precisely, as $n\to\infty$,
\[
 \Delta_n\sim\frac{\sqrt{1+a^{-2}}}{q+2}\frac{r_n}{n}
 \quad\text{and}\quad \delta_n\asymp\frac{x_n^{1-1/q}}{n}.
\]
Consequently,
\begin{equation}\label{eq:log-mesh-decay}
\begin{gathered}
\sum_n\Delta_n^s<\infty\text{ for }s>1+1/(q+1),\quad
x_n^\beta\delta_n\longrightarrow0\text{ for }\beta>0,\\
\text{and}\quad|\log\delta_n|=O(x_n^{1/q})=o(x_n).
\end{gathered}
\end{equation}
Here the limits and asymptotic estimates are taken as $n\to\infty$.
The exponent $\beta$ in these decay assertions is an arbitrary fixed number; it is not the function $m(x)$.
\end{lemma}
\begin{proof}
The logarithmic radius and the strip abscissa are different parameters. On the chosen wall, let
\[
 u(x)=\log\left|\psi\left(x\pm\frac{\pi i}{2}\right)\right|\quad\text{and}\quad u_n=\log r_n=u(x_n).
\]
Lemma~\ref{lem:log-scale} and its differentiated defining formula~\eqref{eq:log-model} give, as $x\to+\infty$,
\[
 u(x)\sim ax^{1/q}\quad\text{and}\quad
 u'(x)=\Rea\left(
       \frac{\psi'(x\pm\pi i/2)}{\psi(x\pm\pi i/2)}\right)
       \sim \frac aq x^{1/q-1}.
\]
By \eqref{eq:outside-log-coordinate} and~\eqref{eq:outside-log-slope}, along $z=e^{u+i\theta_\pm(u)}$ as $u\to+\infty$,
\[
 \log|\tau_L(z)|=(q+2)\log|z|+o(\log|z|)\quad\text{and}\quad
 \frac{d}{du}\log|\tau_L(e^{u+i\theta_\pm(u)})|
       =q+2+o(1).
\]

We first estimate $r_n$ as $n\to\infty$. At the $n$th vertex, $|\tau_L|=2\pi n$, and so, as $n\to\infty$,
\[
 \log(2\pi n)=(q+2)\log r_n+o(\log r_n).
\]
It follows that, as $n\to\infty$,
\[
 \log r_n=\left(\frac1{q+2}+o(1)\right)\log n\quad\text{and}\quad r_n=n^{1/(q+2)+o(1)}.
\]
Combining this with $u(x_n)\sim ax_n^{1/q}$ as $n\to\infty$ also gives
\begin{equation}\label{eq:log-mesh-index}
 \log n=a(q+2) x_n^{1/q}+o(x_n^{1/q}).
\end{equation}

We next compute the physical arclength of one edge. Integrating along the wall gives
\[
 \log\left(1+\frac1n\right)=\int_{u_n}^{u_{n+1}}\frac{d}{du}\log\left|\tau_L\!\left(e^{u+i\theta_\pm(u)}\right)\right|\,du\\=\int_{u_n}^{u_{n+1}}(q+2+o(1))\,du.
\]
Here the error tends uniformly to zero on $[u_n,u_{n+1}]$ as $n\to\infty$. Hence
\[
 u_{n+1}-u_n\sim\frac1{(q+2)n}\text{ as }n\to\infty.
\] 
Along the wall $z=e^{u+i\theta_\pm(u)}$,
\[
 \left|\frac{dz}{du}\right|=e^u\sqrt{1+\theta_\pm'(u)^2}.
\]
As $n\to\infty$, $u_{n+1}-u_n\to0$, so $e^u/r_n\to1$ uniformly on this edge. Also, $\theta_\pm'(u)\to1/a$ as $u\to+\infty$. Therefore, as $n\to\infty$,
\[
 \Delta_n=\int_{u_n}^{u_{n+1}}
       e^u\sqrt{1+\theta_\pm'(u)^2}\,du\sim r_n\sqrt{1+a^{-2}}\,(u_{n+1}-u_n)\sim\frac{\sqrt{1+a^{-2}}}{q+2}\frac{r_n}{n}.
\]
This proves~\eqref{eq:outside-mesh} and its more precise arclength estimate.

Before converting this length to the strip mesh, we check that the width scale varies little along one edge. The relation~\eqref{eq:log-mesh-index} for $\log n$ implies
\[
 \frac{n^{-1}}{x_n^{1/q-1}}
 =\exp\bigl(-\log n+(1-1/q)\log x_n\bigr)
 \longrightarrow0\text{ as }n\to\infty.
\]
If $\delta_n=x_{n+1}-x_n\ge1$, then, for large $n$,
\[
 u_{n+1}-u_n=\int_{x_n}^{x_{n+1}}u'(x)\,dx\ge\int_{x_n}^{x_n+1}u'(x)\,dx\ge c x_n^{1/q-1}.
\]
This contradicts $u_{n+1}-u_n\sim1/((q+2)n)=o(x_n^{1/q-1})$ as $n\to\infty$. Thus $\delta_n<1$ eventually. For $x_n\le x\le x_{n+1}<x_n+1$, the explicit width formula gives
\[
 \frac{\eta(x)}{\eta(x_n)}
 =\exp\left(a\bigl((R_0+x)^{1/q}-(R_0+x_n)^{1/q}\bigr)\right)
       \left(\frac{x}{x_n}\right)^{1/q-1}
 \longrightarrow1
\]
uniformly as $n\to\infty$, since $q>1$.
Using $|\psi'(x\pm\pi i/2)|\asymp\eta(x)$, we can now integrate to obtain
\[
 \Delta_n=\int_{x_n}^{x_{n+1}}
       \left|\psi'\left(x\pm\frac{\pi i}{2}\right)\right|\,dx\asymp\int_{x_n}^{x_{n+1}}\eta(x)\,dx\asymp\eta(x_n)\delta_n.
\]
Since $r_n/r(x_n)\to1$ as $n\to\infty$ and $\eta(x_n)=r(x_n)x_n^{1/q-1}$,
\[
 \delta_n\asymp\frac{\Delta_n}{\eta(x_n)}
       \asymp\frac{r_n/n}{r(x_n)x_n^{1/q-1}}
       \asymp\frac{x_n^{1-1/q}}n.
\]
Taking logarithms yields, as $n\to\infty$,
\[
 \log\delta_n=-\log n+(1-1/q)\log x_n+O(1)\\=-a(q+2) x_n^{1/q}+o(x_n^{1/q}),
\]
because $\log x_n=o(x_n^{1/q})$ as $n\to\infty$. This proves~\eqref{eq:log-small}. For every fixed $\beta>0$,
\[
 \log(x_n^\beta\delta_n)
 =\beta\log x_n-a(q+2) x_n^{1/q}+o(x_n^{1/q})
 \longrightarrow-\infty\text{ as }n\to\infty.
\]
Hence $x_n^\beta\delta_n\to0$ as $n\to\infty$, and $|\log\delta_n|=O(x_n^{1/q})=o(x_n)$ in the same limit follows from $q>1$.

To prove summability, write
\[
 \Delta_n^s=n^{-s(q+1)/(q+2)+o(1)}\text{ as }n\to\infty.
\]
If $s>1+1/(q+1)$, choose $\varepsilon>0$ with $s(q+1)/(q+2)-\varepsilon>1$. For all sufficiently large $n$,
\[
 \Delta_n^s\le n^{-s(q+1)/(q+2)+\varepsilon}.
\]
Comparison with this convergent series proves $\sum_n\Delta_n^s<\infty$. Finally, since $d=1+1/q$,
\[
 1+\frac1{q+1}<1+\frac1q=d.
\]
No assertion about convergence at $s=1+1/(q+1)$ is needed.
\end{proof}

\subsection{A smooth mesh shared by the two models}\label{subsec:smooth-mesh}
On either boundary line put
\[
 A(x)=\frac{\pm\tau_L(\psi(x\pm \frac{\pi i}{2}))}{2\pi i},
\]
choosing the sign so that $A$ is real increasing and $A(x_n)=n$ by~\eqref{eq:narrow-vertices} and~\eqref{eq:log-vertices}. Boundary values and analytic germs are obtained by reflection. Define
\begin{equation*}
 m(x)=A^{-1}(A(x)+1)-x.
\end{equation*}
For sufficiently large $x$, this is positive and smooth, and $m(x_n)=\delta_n$.
The function $A$ converts the inside boundary parameter $x$ into the outside partition coordinate. At the $n$th vertex, $A(x_n)=n$. It is not the affine map $A_t$ used in Lemma~\ref{lem:narrow-scale}. The defining identity is
\[
 A(x+m(x))=A(x)+1.
\]
Thus, increasing $x$ by $m(x)$ increases $A(x)$ by $1$. In particular, substituting $x=x_n$ gives $x_n+m(x_n)=x_{n+1}$ and hence $m(x_n)=\delta_n$. This defines the mesh between the discrete vertices as well, without declaring it piecewise constant.

We now prove the regularity estimates for the smooth interpolation of the two models' discrete meshes.

\begin{lemma}[Smooth mesh estimates]\label{lem:smooth-mesh}
For either model, on a sufficiently far tail, $m$ is positive and smooth, $m(x_n)=\delta_n$, and $m(x)\to0$ as $x\to\infty$. Moreover,
\begin{equation}\label{eq:mesh-reg}
 1+m'(x)=\frac{A'(x)}{A'(x+m(x))}\quad\text{and}\quad
 |m'(x)|+|m''(x)|\le Cm(x).
\end{equation}
In particular,
\begin{equation}\label{eq:consec-mesh}
 \delta_{n+1}-\delta_n=O(\delta_n^2)\text{ as }n\to\infty.
\end{equation}
All comparison constants are independent of the vertex index.
\end{lemma}
\begin{proof}
The choice of sign makes $A$ real and increasing on the tail. Its local derivative is nonzero by the conformal boundary charts, so $A'(x)>0$. It follows that $m$ is positive and smooth. The identity
\[
 A(x+m(x))=A(x)+1
\]
gives $m(x_n)=x_{n+1}-x_n=\delta_n$.
The two preceding parts imply $\delta_n\to0$ as $n\to\infty$. If $x_n\le x\le x_{n+1}$, monotonicity gives
\[
 A(x)+1\in[n+1,n+2]\quad\text{and}\quad
 x+m(x)=A^{-1}(A(x)+1)\in[x_{n+1},x_{n+2}].
\]
Thus
\[
 0<m(x)\le x_{n+2}-x_n=\delta_n+\delta_{n+1}\longrightarrow0
 \text{ as }x\to+\infty,\ n\to\infty.
\]

\emph{Claim.} There is a constant $r_0>0$ such that, for every sufficiently large real $x$, $A$ has a univalent holomorphic continuation to $D(x,r_0)$.

\emph{Case 1: the narrowing model.}
The boundary-reflection argument in the proof of Lemma~\ref{lem:narrow-scale}\,(1), together with the local inversion argument in Lemma~\ref{lem:narrow-lip}, gives univalent continuations of $\tau_L=-e^{\sigma_L}$ to physical disks of a fixed positive radius. Here the reflected charts are chosen small enough that exponentiation is injective. By Lemma~\ref{lem:narrow-lip}, the images of fixed-radius inside charts have size $O(\eta(t))\to0$ as $t\to+\infty$. They therefore lie in these outside continuation disks for large $t$, and composition proves the claim.

\emph{Case 2: the logarithmic model.}
Put $z_0=\psi(x\pm\pi i/2)$. By Lemma~\ref{lem:outside-log-univalent}, the selected outside branch of $\tau_L$ extends univalently to $D(z_0,c|z_0|)$, with $c>0$ independent of $x$. Lemma~\ref{lem:log-scale} gives a univalent continuation of $\psi$ across the boundary and, for a sufficiently small fixed $r_0>0$,
\[
 \left|\psi\left(x+\zeta\pm\frac{\pi i}{2}\right)-z_0\right|
 \le C\eta(x)|\zeta|\quad\text{for}\quad |\zeta|<r_0.
\]
Moreover,
\[
 \frac{\eta(x)}{|z_0|}\sim x^{1/q-1}\longrightarrow0
 \text{ as }x\to+\infty,
\]
since $q>1$. Hence, for all sufficiently large $x$,
\[
 \psi\left(D\left(x\pm\frac{\pi i}{2},r_0\right)\right)
 \subset D(z_0,c|z_0|).
\]
Consequently,
\[
 A(x+\zeta)
 =\frac{\pm\tau_L\!\left(\psi\left(x+\zeta\pm\frac{\pi i}{2}\right)\right)}{2\pi i}
 \quad\text{for}\quad|\zeta|<r_0
\]
defines a univalent holomorphic continuation to $D(x,r_0)$. This proves the claim.

By Koebe's estimates, $A''/A'$ is uniformly bounded on $D(x,r_0/2)$. This quotient is holomorphic because $A'$ has no zeros there, so Cauchy's estimate also bounds its derivative at $x$. Thus
\[
 \left|\frac{A''(x)}{A'(x)}\right|\le C\quad\text{and}\quad
 \left|\left(\frac{A''}{A'}\right)'(x)\right|\le C,
\]
with $C$ independent of $x$.

We next estimate $m'$. Differentiating the defining identity by the chain rule gives
\[
 A'(x+m(x))(1+m'(x))=A'(x),
\]
which proves the equality in~\eqref{eq:mesh-reg}. Since both derivatives of $A$ are positive, $1+m'(x)>0$, and
\[
 \log(1+m'(x))=\log A'(x)-\log A'(x+m(x))
 =-\int_x^{x+m(x)}\frac{A''(s)}{A'(s)}\,ds.
\]
The bound for $A''/A'$ implies
\[
 e^{-Cm(x)}\le1+m'(x)\le e^{Cm(x)}.
\]
Because $m(x)\to0$ as $x\to+\infty$, this gives $|m'(x)|\le Cm(x)$ after increasing $C$.

To estimate $m''$, differentiate the logarithmic identity once more and separate the two error terms:
\[
 \begin{aligned}
 \frac{m''(x)}{1+m'(x)}
 &=\frac{A''(x)}{A'(x)}
   -\frac{A''(x+m(x))}{A'(x+m(x))}(1+m'(x))\\
 &=\left[\frac{A''(x)}{A'(x)}
   -\frac{A''(x+m(x))}{A'(x+m(x))}\right]
   -\frac{A''(x+m(x))}{A'(x+m(x))}m'(x).
 \end{aligned}
\]
The bound for $(A''/A')'$ controls the bracket by $Cm(x)$. The other term is at most $C|m'(x)|\le Cm(x)$. Since $1+m'(x)=1+o(1)$ as $x\to+\infty$, this factor is bounded and $|m''(x)|\le Cm(x)$. This completes~\eqref{eq:mesh-reg}.

It remains to compare consecutive mesh widths. From $m(x_n)=\delta_n$,
\[
 \delta_{n+1}-\delta_n=\int_{x_n}^{x_{n+1}}m'(x)\,dx.
\]
As $m>0$ and $|m'|\le Cm$, one has $|(\,\log m)'|\le C$. Hence, for $x_n\le x\le x_{n+1}$,
\[
 m(x)\le m(x_n)e^{C(x-x_n)}\le\delta_ne^{C\delta_n}.
\]
Therefore, as $n\to\infty$,
\[
 |\delta_{n+1}-\delta_n|
 \le C\int_{x_n}^{x_{n+1}}m(x)\,dx
 \le C\delta_ne^{C\delta_n}(x_{n+1}-x_n)
 =Ce^{C\delta_n}\delta_n^2=O(\delta_n^2).
\]
This proves~\eqref{eq:consec-mesh}. The proofs apply separately to each tail, so no alignment of the two boundary meshes has been used.
\end{proof}

\section{The common graded slit construction}\label{sec:slits}
This section is in strip coordinates. Lemmas~\ref{lem:central}--\ref{lem:tip} and Lemma~\ref{lem:column-estimates}\,(2) apply to both models of Section~\ref{sec:models}. Lemma~\ref{lem:column-estimates}\,(1) is stated only for the logarithmic model. The key idea is to balance a large entrance value with exponential attenuation along a narrow gap. Throughout this section, tail estimates hold for sufficiently large $n$ or $x$, as applicable. Unless indicated otherwise, the constants are uniform in $n$ and in the spatial variables over the stated ranges.

\subsection{Slit heights and central normalization}
On each sufficiently far tail define
\begin{equation}\label{eq:height}
 h(x)=\frac{m(x)}\pi\bigl(x+\log m(x)\bigr),\quad M_n=e^{x_n}\delta_n\quad\text{and}\quad
 h_n=h(x_n)=\frac{\delta_n}{\pi}\log M_n.
\end{equation}
Here $h_n$ is the slit height in strip coordinates. The physical slit height is comparable to $\eta(t_n)h_n$ and $\eta(x_n)h_n$ in the narrowing and logarithmic models, respectively. Write \(\eta_n=\eta(t_n)\) or \(\eta(x_n)\), according to the model.
The following estimates follow from~\eqref{eq:narrow-small} or~\eqref{eq:log-small}, together with~\eqref{eq:mesh-reg}.
For both models, $h_n>0$ and, as $n\to\infty$, $h_n\to0$, $h_n/\delta_n\to\infty$, and
\begin{equation}\label{eq:height-regularity}
 |h'(x)|+|h''(x)|\le Ch(x),\quad
 h_n\log(1/\delta_n)\to0,\quad h_{n+1}-h_n=o(\delta_n).
\end{equation}
The first estimate holds for sufficiently large real $x$; the last two are as $n\to\infty$.
Indeed, the mesh estimates give $\log m(x)=o(x)$ and hence $h(x)\sim xm(x)/\pi$ as $x\to+\infty$.
Differentiating~\eqref{eq:height} and using~\eqref{eq:mesh-reg} gives the derivative bounds.
The rapid decay of $m$ gives the second estimate, and the mean value theorem on $[x_n,x_{n+1}]$ gives the last one.
Moreover, $h(x)$ tends to zero faster than every negative power of $x$ as $x\to+\infty$, and $\log(1/\delta_n)=o(x_n^2)$ as $n\to\infty$.

From $x_n+\frac{\pi i}{2}$ remove a vertical segment directed down into the strip, of height $h_n$. Similarly, remove vertical slits extending upward from the lower wall, using its own sequence of vertices and heights. Omit finitely many initial slits so that every height is below $\pi/8$. Let $U$ be the resulting slit domain. Choose a conformal $\Phi:U\to S$ preserving the two ends and the upper and lower boundary arcs, and define
\begin{equation*}
 \tau_U=e^\Phi:U\to\HH\quad\text{and}\quad Q=\Rea\tau_U.
\end{equation*}
The positive harmonic function $Q$ vanishes on the sides of the slits and on the horizontal parts of the boundary.

\begin{figure}[H]
  \centering
  \includegraphics[width=0.9\linewidth]{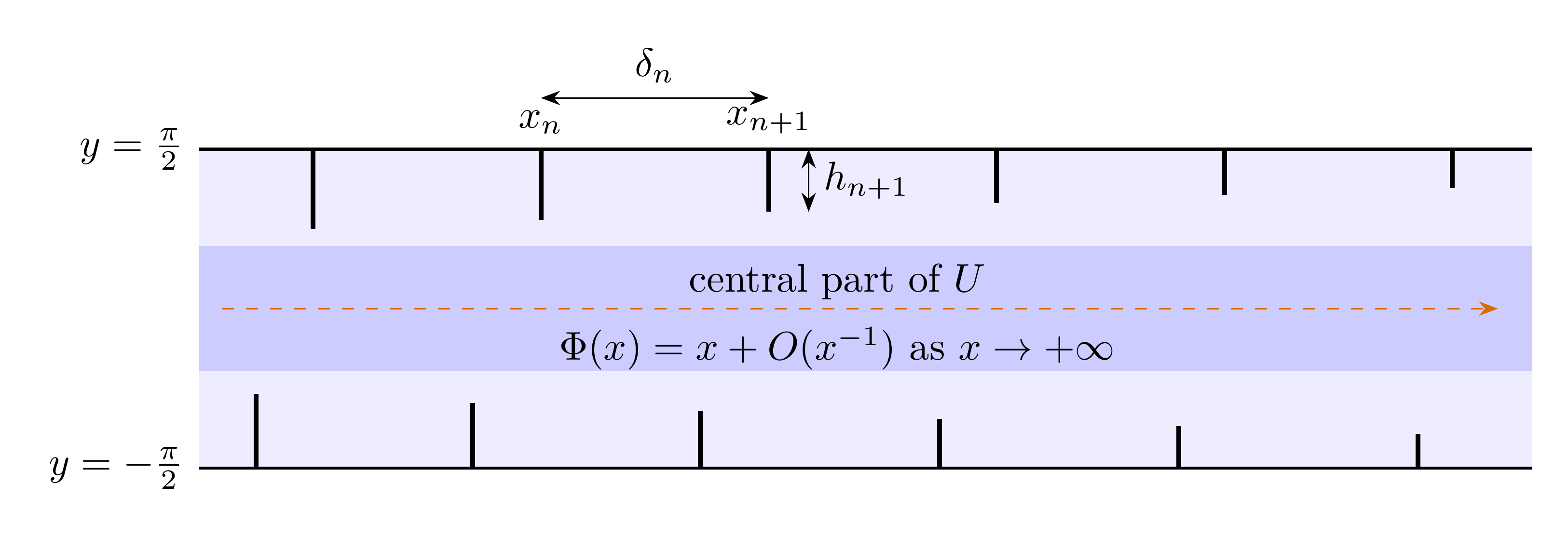}
  \caption[The slit domain before pullback]{%
    The slit domain $U\subset S$ before pullback to the physical spiral.
    Each wall has its own mesh, so the slit roots on the two walls
    need not have the same horizontal coordinates.}
  \label{fig:slit-strip}
\end{figure}

\begin{lemma}[Central normalization]\label{lem:central}
After a real translation of $\Phi$, one has
\begin{equation}\label{eq:central}
 \begin{aligned}
 &|\Ima\Phi(x+iy)-y|\le\frac C{9+x^2}
 &&\text{for }x+iy\in U,\\
 &\Phi(x)=x+O(x^{-1})
 &&\text{for }x\in\R,\ x\to+\infty.
 \end{aligned}
\end{equation}
\end{lemma}
\begin{proof}
Let \(v(x,y)=\Ima\Phi(x+iy)-y.\) Both terms are harmonic, so $v$ is harmonic in $U$. Since $U\subset S$ and $\Phi(U)=S$, one also has $|v|\le\pi$ throughout $U$.

We first identify the boundary values, understood from the corresponding boundary sides. On the upper and lower horizontal parts of the boundary, $\Ima\Phi=y=\pi/2$ and $\Ima\Phi=y=-\pi/2$, respectively, so $v=0$. At depth $t$ on an upper slit, with $0\le t\le h_n$, one has
\[
 y=\frac\pi2-t,\quad \Ima\Phi=\frac\pi2\quad\text{and}\quad v=t.
\]
Both sides of that slit belong to the upper boundary arc. On a lower slit the corresponding value is $v=-t$. Consequently, on each slit, \(|v|\le h_n.\) The same bounds hold at the roots and tips by boundary limits.

Choose
\[
 B(x,y)=K\frac{\cos(y/2)}{9+x^2},
\]
where the fixed constant $K>0$ will be chosen below. A direct calculation gives \(\Delta B\le-\frac18B.\) Thus $v-B$ and $-v-B$ are subharmonic in $U$.

Since $|y|\le\pi/2$, one has $\cos(y/2)\ge1/\sqrt2$. By~\eqref{eq:height}, together with~\eqref{eq:narrow-small} or~\eqref{eq:log-small}, for all sufficiently large $n$,
\[
 0\le(9+x_n^2)h_n
 \le\frac{9x_n+x_n^3}{\pi}\,\delta_n
 \longrightarrow0\text{ as }n\to\infty,
\]
where $\delta_n<1$. Hence
\[
 C_0:=\sup_{\text{slits on both walls}}(9+x_n^2)h_n<\infty.
\]
Taking $K\ge\sqrt2 C_0$ makes
\[
 B(x_n,y)\ge\frac{K}{\sqrt2(9+x_n^2)}
 \ge h_n\ge|v|
\]
on every slit. On the horizontal parts of the boundary, $v=0$ and $B>0$. Hence $v-B$ and $-v-B$ have nonpositive boundary values on all these boundary sides.

To handle the two unbounded ends, fix $\varepsilon>0$ and put
\[
 H_\varepsilon(x,y)=\varepsilon\cosh(x/2)\cos(y/2).
\]
This function is positive and harmonic in $S$, because its second derivatives in $x$ and $y$ cancel. On the vertical sides $x=\pm R$,
\[
 H_\varepsilon(x,y)\ge\frac{\varepsilon}{\sqrt2}\cosh(R/2)\ge\pi
\]
for all sufficiently large $R$. Since $|v|\le\pi$, both functions
\[
 v-B-H_\varepsilon\quad\text{and}\quad -v-B-H_\varepsilon
\]
are nonpositive on every boundary side of the bounded truncation
$U_R=U\cap\{-R<x<R\}$. They are subharmonic there, since subtracting the harmonic function $H_\varepsilon$ preserves subharmonicity. The maximum principle gives
\[
 |v|\le B+H_\varepsilon\quad\text{in }U_R.
\]
First let $R\to\infty$ with $\varepsilon$ fixed, and then let $\varepsilon\downarrow0$ at any fixed point of $U$. We obtain
\[
 |v(x,y)|\le B(x,y)\le\frac{K}{9+x^2},
\]
which is the first assertion of~\eqref{eq:central}.

We next pass from this bound to derivatives on the real axis. Every slit has height less than $\pi/8$, so the fixed central substrip $\{|y|<\pi/4\}$ is contained in $U$. In particular, the disk $D(x,r_*)$, with $r_*=\pi/8$, lies in $U$ for every real $x$. On this disk, the preceding bound is $O(x^{-2})$ as $x\to+\infty$. By the interior gradient estimate for harmonic functions, we obtain
\[
 |v_x(x,0)|+|v_y(x,0)|
 \le\frac{C}{r_*}\sup_{D(x,r_*)}|v|
 =O(x^{-2})\text{ as }x\to+\infty.
\]
Here the subscripts denote real partial derivatives. The Cauchy--Riemann equations give
\[
 \Phi'(x)-1=v_y(x,0)+i\,v_x(x,0)=O(x^{-2})\text{ as }x\to+\infty.
\]

Finally, $\Phi'(x)-1$ is integrable on a real tail. Hence there is a finite constant
\[
 c=\lim_{x\to+\infty}\bigl(\Phi(x)-x\bigr)
\]
and
\[
 \Phi(x)-x-c
 =-\int_x^\infty\bigl(\Phi'(t)-1\bigr)\,dt
 =O(x^{-1})\text{ as }x\to+\infty.
\]
This constant is real, because
\[
 \Ima\bigl(\Phi(x)-x\bigr)=v(x,0)=O(x^{-2})\longrightarrow0\text{ as }x\to+\infty.
\]
Replacing $\Phi$ by $\Phi-c$ is a real translation of the target strip. It preserves the two ends, the upper and lower boundary arcs, and the established imaginary-part bound. With this normalization, $\Phi(x)=x+O(x^{-1})$ as $x\to+\infty$, as required.
\end{proof}

\subsection{The entrance estimate}
Fix an upper gap. Let
\begin{equation}\label{eq:gap-points}
 B_n=\frac{x_n+x_{n+1}}2\quad\text{and}\quad
 C_n=B_n+i\left(\frac{\pi}{2}-h(B_n)-3\delta_n\right).
\end{equation}
$B_n$ is in the middle of the strip; $C_n$ is just below the slit tips, at the entrance to the narrow gap. The lower wall remains at a fixed distance from the comparison region.

\begin{lemma}[Uniform entrance distance]\label{lem:entrance-distance}
For either wall, as $n\to\infty$,
\begin{equation}\label{eq:entrance-distance}
 \rho_U(B_n,C_n)=\frac12\log(1/\delta_n)+O(1),
\end{equation}
with a constant independent of $n$.
\end{lemma}
\begin{proof}
Choose the window size and cutoff functions independently of $n$. Let $A>0$ satisfy $|m'|\le Am$, and choose a fixed $0<r_0<\pi/16$ so small that $e^{4Ar_0}<3/2$. Since $m(x_j)=\delta_j$, for mesh vertices in $[B_n-2r_0,B_n+2r_0]$ we have
\[
 \left|\log\frac{\delta_j}{\delta_k}\right|
 =\left|\int_{x_k}^{x_j}\frac{m'(x)}{m(x)}\,dx\right|
 \le A|x_j-x_k|\le4Ar_0.
\]
For large $n$, the preceding estimate applies to all mesh intervals meeting $[B_n-r_0,B_n+r_0]$. Their lengths are comparable to $\delta_n$, with gap ratio at most $e^{4Ar_0}<3/2$. The height estimates give
\[
 h(x)\asymp h_n\quad\text{and}\quad |h'(x)|+|h''(x)|\le Ch_n
 \text{ for }|x-B_n|\le2r_0
\]
for all sufficiently large $n$.

Take smooth cutoffs $\chi,\chi_0$, with $\chi=1$ on $[-r_0,r_0]$ and zero outside $[-2r_0,2r_0]$, and $\chi_0=0$ near and below $y=0$, $\chi_0=1$ near $y=\pi/2$. Extend $\chi(x-B_n)h(x)$ by zero outside the larger window. For $x,y\in\R$, define
\[
 F_n(x+iy)=x+i\bigl(y+\chi(x-B_n)\chi_0(y)h(x)\bigr).
\]
This is a global quasiconformal homeomorphism. Indeed, its real coordinate is unchanged, its vertical derivative is positive for large $n$, and its derivative is $I+O(h_n)$. Its dilatation satisfies $K(F_n)=1+O(h_n)$ as $n\to\infty$.
It fixes $B_n$, sends $C_n$ to $B_n+i(\pi/2-3\delta_n)$, and maps the nearby tip graph to $y=\pi/2$.

\emph{The upper bound.}
The domain
\[
 E_n=\{x+iy:|x-B_n|<r_0,\ -\pi/4<y<\pi/2-h(x)\}
\]
lies in $U$: it is below the upper slit tips and stays a fixed distance from the lower slits. After flattening and translation, it becomes the fixed rectangle
\[
 \mathcal R=F_n(E_n)-B_n
 =\{z:|\Rea z|<r_0,\ -\pi/4<\Ima z<\pi/2\}.
\]
The two points become $0$ and $z_n=i(\pi/2-3\delta_n)$, with
$\dist(z_n,\partial\mathcal R)=3\delta_n$ for large $n$.
The corresponding hyperbolic distance can be estimated directly from this Euclidean distance. Fix $r=r_0/2$ and put $a=i(\pi/2-r)$. Then
\[
 D(a,r)\subset\mathcal R\subset\mathcal H:=\{z:\Ima z<\pi/2\}.
\]
For $3\delta_n<r$, the explicit disk and half-plane formulas, together with domain monotonicity and the triangle inequality, give
\[
 \begin{aligned}
 \frac12\log\frac{\pi/2}{3\delta_n}
 &=\rho_{\mathcal H}(0,z_n)
 \le\rho_{\mathcal R}(0,z_n)\le\rho_{\mathcal R}(0,a)+\rho_{D(a,r)}(a,z_n)\\
 &=\rho_{\mathcal R}(0,a)
   +\frac12\log\frac{2r-3\delta_n}{3\delta_n}.
 \end{aligned}
\]
Here $\rho_{\mathcal R}(0,a)$ is fixed. Hence
\[
 \rho_{F_n(E_n)}\bigl(F_n(B_n),F_n(C_n)\bigr)
 =\rho_{\mathcal R}(0,z_n)
 =\frac12\log(1/\delta_n)+O(1).
\]

\emph{The lower bound.}
Remove only the upper-wall slit tips whose horizontal coordinates lie in $[B_n-r_0,B_n+r_0]$:
\[
 D_n=\C\setminus
 \bigl\{x_j+i(\pi/2-h_j):|x_j-B_n|\le r_0\bigr\}.
\]
Then $U\subset D_n$. Let $N_n$ be the number of selected tips, and let $a_n^-,a_n^+$ be their smallest and largest horizontal coordinates. The mesh estimates above imply
\[
 N_n\asymp r_0/\delta_n\quad\text{and}\quad
 a_n^\pm=B_n\pm r_0+O(\delta_n).
\]
Set
\[
 p_n=\frac{a_n^-+a_n^+}{2}=B_n+O(\delta_n)\quad\text{and}\quad
 s_n=\frac{a_n^+-a_n^-}{2}=r_0+O(\delta_n).
\]
The map $F_n$ preserves horizontal coordinates and sends the two extreme tips
$a_n^\pm+i(\pi/2-h(a_n^\pm))$ to $a_n^\pm+\pi i/2$.
The affine conformal map
\[
 G_n(Z)=-\frac{Z-p_n-\pi i/2}{s_n}
\]
sends these extreme punctures of $F_n(D_n)$ to $\mp1$, respectively, and the region below their line to the upper half-plane. Thus $G_n\circ F_n$ sends the original extreme tips to $\mp1$. In increasing order, write the normalized punctures as
\[
 -1=\zeta_{n,1}<\cdots<\zeta_{n,N_n}=1\quad\text{and}\quad
 \Omega_n:=G_n(F_n(D_n))
 =\C\setminus\{\zeta_{n,1},\ldots,\zeta_{n,N_n}\}.
\]
Their average gap is $2/(N_n-1)$, and their largest-to-smallest gap ratio is at most $3/2$. Therefore
\[
 \frac{4}{3(N_n-1)}
 \le\zeta_{n,j+1}-\zeta_{n,j}
 \le\frac{3}{N_n-1}
 \quad\text{for}\quad 1\le j<N_n.
\]
For $N_n\ge4$, these bounds imply the hypotheses of the punctured-plane estimate~\cite[Lemma~11.1]{AlbrechtBishop}:
\[
 \frac1{N_n}\le\zeta_{n,j+1}-\zeta_{n,j}\le\frac4{N_n}\quad\text{and}\quad
 \rho_{\Omega_n}(i,i/N_n)=\frac12\log N_n+O(1).
\]

The images of the actual points are
\[
 \widetilde B_n:=G_n(F_n(B_n))
 =\frac{p_n-B_n}{s_n}+i\frac{\pi}{2s_n},\ 
 \widetilde C_n:=G_n(F_n(C_n))
 =\frac{p_n-B_n}{s_n}+i\frac{3\delta_n}{s_n}.
\]
Their real parts are $O(1/N_n)$, while
\[
 \Ima\widetilde B_n\asymp1\quad\text{and}\quad
 \Ima\widetilde C_n\asymp1/N_n.
\]
It follows that $\widetilde B_n$ is at bounded hyperbolic distance from $i$ in the upper half-plane. For the second point, the upper-half-plane isometry $Z\mapsto N_nZ$ places $N_n\widetilde C_n$ in a fixed compact subset of that half-plane, so its distance from $i$ is also bounded. Thus, writing $\mathcal H_+=\{Z:\Ima Z>0\}$,
\[
 \rho_{\mathcal H_+}(\widetilde B_n,i)=O(1)\quad\text{and}\quad
 \rho_{\mathcal H_+}(\widetilde C_n,i/N_n)=O(1).
\]
Since $\mathcal H_+\subset\Omega_n$, these errors remain bounded in $\Omega_n$. The triangle inequality and the punctured-plane estimate yield
\[
 \begin{aligned}
 \rho_{F_n(D_n)}\bigl(F_n(B_n),F_n(C_n)\bigr)
 &=\rho_{\Omega_n}(\widetilde B_n,\widetilde C_n)=\rho_{\Omega_n}(i,i/N_n)+O(1)\\
 &=\frac12\log N_n+O(1)
  =\frac12\log(1/\delta_n)+O(1).
 \end{aligned}
\]
The first equality uses the conformality of $G_n$.

\emph{Returning to the original domains.}
For $\mathcal D=E_n$ or $D_n$, let
\[
 d_{\mathcal D}=\rho_{\mathcal D}(B_n,C_n)\quad\text{and}\quad
 d_{\mathcal D}'=\rho_{F_n(\mathcal D)}\bigl(F_n(B_n),F_n(C_n)\bigr).
\]
The disk estimate~\cite[Theorem~5.1]{EMM} gives
\[
 K(F_n)^{-1}d_{\mathcal D}-C\le d_{\mathcal D}'\le K(F_n)d_{\mathcal D}+C
 \quad\text{when}\quad\mathcal D=E_n\text{ or }D_n,
\]
where $C$ is independent of $n$, since $1\le K(F_n)\le2$ for large $n$.
For $E_n$, apply the estimate after conformal maps to the disk.
For $D_n$, lift $F_n$ to a $K(F_n)$-quasiconformal homeomorphism between disk universal covers.
Fix the starting lift, apply the disk estimate, and take the infimum over all endpoint lifts, using the covering formula~\cite[Remark~5.1]{Rugh}.

The preceding comparisons give $d_{\mathcal D}'=\frac12\log(1/\delta_n)+O(1)$.
Since $K(F_n)=1+O(h_n)$ and $h_n\log(1/\delta_n)\to0$ by~\eqref{eq:height-regularity}, we obtain, as $n\to\infty$,
\[
 d_{\mathcal D}=\frac12\log(1/\delta_n)
 +O\!\bigl(h_n\log(1/\delta_n)\bigr)+O(1)
 =\frac12\log(1/\delta_n)+O(1).
\]
Finally, $E_n\subset U\subset D_n$ gives
\[
 \rho_{D_n}(B_n,C_n)
 \le\rho_U(B_n,C_n)
 \le\rho_{E_n}(B_n,C_n),
\]
which proves~\eqref{eq:entrance-distance}. The lower-wall case follows by reflection.
\end{proof}

\begin{lemma}[Horizontal coordinate and entrance value]\label{lem:entrance-value}
With the normalization of Lemma~\ref{lem:central}, as $n\to\infty$,
\begin{align}
 |\Rea\Phi(C_n)-B_n|&\le CB_n^{-1}+CB_n^{-2}\log(1/\delta_n)=o(1),\label{eq:entrance-real}\\
 Q(C_n)&\asymp e^{x_n}\delta_n=M_n.\label{eq:entrance-value}
\end{align}
\end{lemma}
\begin{proof}
For large $n$, let
\[
 0\le y\le\Ima C_n=\frac\pi2-h(B_n)-3\delta_n\quad\text{and}\quad
 t=\frac\pi2-h(B_n)-y.
\]
Then $3\delta_n\le t\le\pi/2$. We first show that $D(B_n+iy,t/4)\subset U$, which gives the required boundary-distance estimate. By~\eqref{eq:height-regularity}, $|h'|\le1$ on $|x-B_n|\le\pi/8$ for all large $n$.
Hence each upper slit with $|x_j-B_n|<t/4$ satisfies
\[
 \left(\frac\pi2-h_j\right)-y
 =t+h(B_n)-h(x_j)
 \ge t-|x_j-B_n|>\frac{3t}{4}.
\]
The disk $D(B_n+iy,t/4)$ therefore misses every upper slit.
It lies in $-\pi/8<\Ima z<\pi/2$, whereas all lower slits lie at heights at most $-3\pi/8$.
Consequently,
\[
 \dist(B_n+iy,\partial U)
 \ge\frac14\left(\frac\pi2-h(B_n)-y\right),
\]
uniformly for the stated range of $y$.

Let $v(x,y)=\Ima\Phi(x+iy)-y$. By~\eqref{eq:central}, $|v|\le CB_n^{-2}$ on this disk, since its horizontal coordinates differ from $B_n$ by at most $\pi/8$.
The interior gradient estimate for harmonic functions gives
\[
 |v_x(B_n,y)|
 \le\frac C{t}\sup_{D(B_n+iy,t/4)}|v|
 \le\frac C{B_n^2\bigl(\pi/2-h(B_n)-y\bigr)}
 \quad\text{for}\quad 0\le y\le\Ima C_n.
\]
The Cauchy--Riemann equations give $\partial_y\Rea\Phi(x+iy)=-v_x(x,y)$. Hence
\[
 |\Rea\Phi(C_n)-\Rea\Phi(B_n)|\le\frac C{B_n^2}\int_0^{\Ima C_n}
       \frac{dy}{\pi/2-h(B_n)-y}=\frac C{B_n^2}
   \log\frac{\pi/2-h(B_n)}{3\delta_n}.
\]
As $n\to\infty$, this logarithm equals $\log(1/\delta_n)+O(1)$ because $h(B_n)\to0$.
Combining this bound with $\Rea\Phi(B_n)=B_n+O(B_n^{-1})$ proves~\eqref{eq:entrance-real}.
Its right-hand side is $o(1)$ since $B_n\sim x_n$ and $\log(1/\delta_n)=o(x_n^2)$ as $n\to\infty$.

Write $\Phi(B_n)=u_0+iv_0$ and $\Phi(C_n)=u_1+iv_1$.
By~\eqref{eq:central} and the estimate just proved, as $n\to\infty$ we have
$v_0\to0$, $u_1-u_0=o(1)$ and $v_1\to\pi/2$.
Conformal invariance and~\eqref{eq:strip-hyperbolic-distance} give
\[
 \cosh\bigl(2\rho_U(B_n,C_n)\bigr)
 =\frac{\cosh(u_1-u_0)-\sin v_1\sin v_0}
        {\cos v_1\cos v_0}
 \asymp\frac1{\cos v_1}.
\]
Lemma~\ref{lem:entrance-distance} makes the left-hand side comparable to $\delta_n^{-1}$, so $\cos v_1\asymp\delta_n$.
Finally, $u_1=B_n+o(1)=x_n+o(1)$ as $n\to\infty$, and therefore
$Q(C_n)=e^{u_1}\cos v_1\asymp e^{x_n}\delta_n=M_n$, proving~\eqref{eq:entrance-value}.
\end{proof}
Notice that the normalization has been proved before any conformal length of a trunk edge is used. It therefore does not assume the matching conclusion.

\subsection{Boundary derivatives from a rectangular gap}
Use downward depth $y$ from the upper wall and the horizontal coordinate $s=x-x_n$. In these coordinates, $Q_n(s,y):=Q\bigl(x_n+s+i(\pi/2-y)\bigr)$. Put $\ell=\min(h_n,h_{n+1})-\delta_n=h_n-\delta_n+o(\delta_n)$ as $n\to\infty$. In the rectangle $0<s<\delta_n$, $0<y<\ell$, the function $Q_n$ is zero on three sides. For its nonnegative crosscut data $g_0(s)=Q_n(s,\ell)$, define
\[
 b_k=\frac{2}{\delta_n}\int_0^{\delta_n} g_0(s)\sin(k\pi s/\delta_n)\,ds
 \text{ for }k\ge1.
\]
Solving this Dirichlet problem by separation of variables gives
\begin{equation}\label{eq:sine}
 Q_n(s,y)=\sum_{k\ge1}b_k\sin(k\pi s/\delta_n)
 \frac{\sinh(k\pi y/\delta_n)}{\sinh(k\pi\ell/\delta_n)}
 \quad\text{for}\quad 0<s<\delta_n,\ 0<y<\ell.
\end{equation}
Positivity and $|\sin(k\theta)|\le k\sin\theta$ for $0\le\theta\le\pi$ imply $|b_k|\le k b_1$. At $(\delta_n/2,\ell-2\delta_n)$, the first term dominates; a bounded Harnack chain reaches $C_n$. Lemma~\ref{lem:entrance-value} gives the first estimate below; the second then follows from~\eqref{eq:height} and $\ell=h_n-\delta_n+o(\delta_n)$:
\begin{equation*}
 b_1\asymp M_n\quad\text{and}\quad b_1e^{-\pi\ell/\delta_n}\asymp1.
\end{equation*}
Differentiating~\eqref{eq:sine} on a vertical side, the relative higher-mode error for $y\le\ell-2\delta_n$ is bounded, on a far tail, by a constant arbitrarily close to $\sum_{k\ge2}k^3e^{-2\pi(k-1)}<0.016$ as $n\to\infty$. Therefore
\begin{equation}\label{eq:side-der}
 \left|\tau_U'\!\left(x_j+i\left(\frac\pi2-y\right)\right)\right|
 \asymp\delta_n^{-1}\sinh(\pi y/\delta_n)
 \quad\text{for}\quad0<y\le h_n-4\delta_n,
\end{equation}
where $j=n$ or $n+1$ and the derivative is taken from inside this gap.
This is comparable to $y/\delta_n^2$ for $0<y\le\delta_n$, and to $\delta_n^{-1}e^{\pi y/\delta_n}$ for $\delta_n\le y\le h_n-4\delta_n$. Integrating the normal derivative on the trunk interval $I_n=[x_n,x_{n+1}]+\frac{\pi i}{2}$ gives
\begin{equation}\label{eq:trunk-size}
 |\tau_U(I_n)|=
 \sum_{k\text{ odd}}\frac{2b_k}{\sinh(k\pi\ell/\delta_n)}\asymp1.
\end{equation}
The first term also dominates this series.

\Needspace{8\baselineskip}
\begin{lemma}[The tip profile]\label{lem:tip}
Let $p=x_n+i(\pi/2-h_n)$ and $z=p+iu$ on the upper slit, where $u=h_n-y$ is the distance from its tip and $y$ is the downward depth from the upper wall. For $0<u\le c\delta_n$, the boundary derivative on either side satisfies
\begin{equation}\label{eq:tip-der}
 |\tau_U'(z)|\asymp\frac{M_n}{\sqrt{\delta_n u}}.
\end{equation}
For $c\delta_n\le u\le4\delta_n$, it satisfies $|\tau_U'(z)|\asymp M_n/\delta_n$. The same estimates hold for the lower slits after reflection.
\end{lemma}
\begin{proof}
We use a square-root coordinate to open the slit near its tip. Schwarz reflection extends the resulting map univalently to a disk. Koebe's distortion theorem and the chain rule then give the desired estimate for $\tau_U'$.

The mesh and height estimates allow us to choose a fixed small $r>0$ such that $D(p,r\delta_n)$ meets only this slit for all large $n$. The square-root coordinate $z=p+i\zeta^2$ opens the slit, with $\zeta=\pm\sqrt u$ corresponding to its two boundary sides. For $\Ima\zeta>0$ and $|\zeta|<\sqrt{r\delta_n}$, put
\[
 F(\zeta)=\tau_U(p+i\zeta^2).
\]
On the real diameter, $F$ has continuous, purely imaginary and injective boundary values. Schwarz reflection extends $F$ univalently to the full disk, since the open half-disks map into opposite half-planes. Figure~\ref{fig:tip-opening} shows how the square-root coordinate
separates the two sides of the slit near its tip.

\begin{figure}[htbp]
  \centering
  \includegraphics[width=0.9\linewidth]{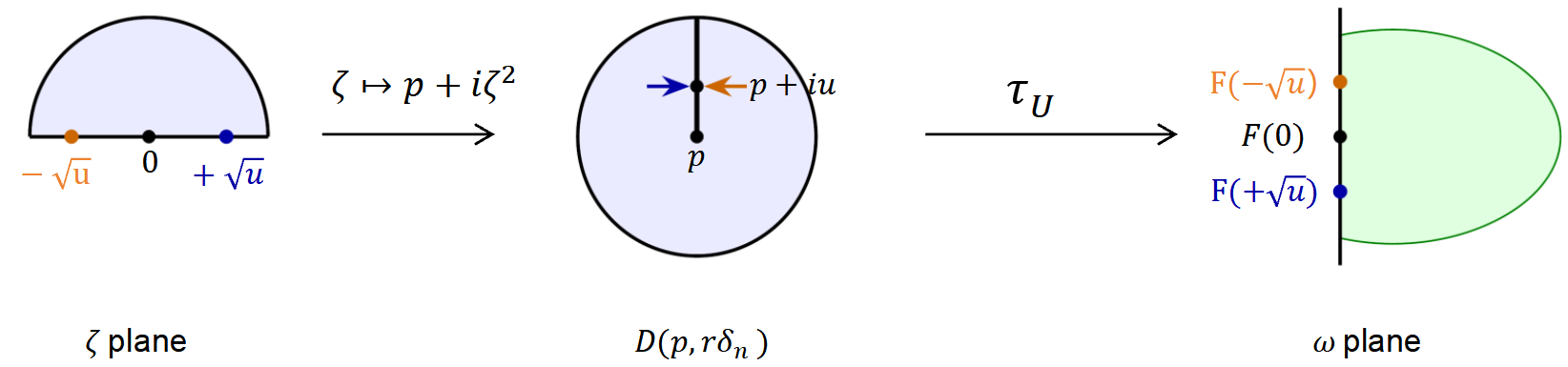}
  \caption[Opening a slit at its tip]{%
    The map $z=p+i\zeta^2$ takes the upper half-disk onto the slit
    disk centered at $p$. The origin corresponds to the tip $p$,
    and the points $\zeta=\pm\sqrt u$ correspond to the two sides
    of the slit at $p+iu$. The right panel shows the image of the
    upper half-disk under $\omega=F(\zeta)=\tau_U(p+i\zeta^2)$.
    The two marked sides have distinct images on the imaginary axis.}
  \label{fig:tip-opening}
\end{figure}

We next compare the value below the tip with the entrance value. The estimates $h(B_n)-h_n=o(\delta_n)$ and $h_{n+1}-h_n=o(\delta_n)$ hold as $n\to\infty$. The path from $p-ir\delta_n/4$ to $p-3i\delta_n$, then horizontally to $B_n+i(\pi/2-h_n-3\delta_n)$ and finally to $C_n$, has length $O(\delta_n)$ and stays at distance at least $c_r\delta_n$ from $\partial U$ for large $n$. Cover this path by overlapping disks of radii comparable to $\delta_n$, with their doubled disks contained in $U$. The number of disks and their overlap ratios can be chosen independently of $n$. Repeated application of Harnack's inequality, followed by Lemma~\ref{lem:entrance-value}, gives
\[
 Q(p-ir\delta_n/4)\asymp Q(C_n)\asymp M_n.
\]

Set $\zeta=\sqrt{r\delta_n}\,\xi$. The function
\[
 \frac{\Rea F(\sqrt{r\delta_n}\,\xi)}{M_n}
 \quad\text{for}\quad|\xi|<1,\ \Ima\xi>0
\]
is positive and harmonic in the unit upper half-disk. It vanishes on the real diameter.
The preceding Harnack comparison and~\eqref{eq:entrance-value} give
\[
 \frac{\Rea F(i\sqrt{r\delta_n}/2)}{M_n}
 =\frac{Q(p-ir\delta_n/4)}{M_n}
 \asymp\frac{Q(C_n)}{M_n}\asymp1.
\]
Thus its value at $i/2$ is comparable to $1$.
By boundary Harnack, its inward normal derivative at the origin is also comparable to $1$.
The Cauchy--Riemann equations and the scaling give \(|F'(0)|\asymp\frac{M_n}{\sqrt{\delta_n}}.\)
Koebe distortion for the reflected univalent map then yields
\[
 |F'(\zeta)|\asymp\frac{M_n}{\sqrt{\delta_n}}
 \quad\text{for}\quad|\zeta|\le\frac{\sqrt{r\delta_n}}2.
\]
Fix $0<c\le r/4$. For $0<u\le c\delta_n$, both $\zeta=\pm\sqrt u$ lie in this disk.
The chain rule on the corresponding boundary side gives
\[
 |\tau_U'(p+iu)|
 =\frac{|F'(\pm\sqrt u)|}{2\sqrt u}
 \asymp\frac{M_n}{\sqrt{\delta_n u}}
 \quad\text{for}\quad 0<u\le c\delta_n,
\]
proving~\eqref{eq:tip-der} on both sides.

For $c\delta_n\le u\le4\delta_n$, choose a fixed $\varepsilon>0$ sufficiently small compared with $c$. A disk of radius $2\varepsilon\delta_n$ centered at $z=p+iu$ meets only the straight slit and avoids its tip. On either side, the point $z\pm\varepsilon\delta_n$ can be joined to $C_n$ by a Harnack chain with a number of disks independent of $n$. Hence
\[
 Q(z\pm\varepsilon\delta_n)\asymp Q(C_n)\asymp M_n
 \quad\text{for}\quad z=p+iu,\ c\delta_n\le u\le4\delta_n,
\]
where the second comparison is Lemma~\ref{lem:entrance-value}.
Boundary Harnack in these side charts gives
\begin{equation}\label{eq:transition-derivative}
 |\tau_U'(z)|\asymp\frac{Q(z\pm\varepsilon\delta_n)}{\delta_n}
 \asymp\frac{M_n}{\delta_n},
\end{equation}
with constants uniform in $u$ and $n$.

For the ordinary part, the adjacent gap widths and entrance values satisfy
\[
 \frac1{\delta_{n-1}}-\frac1{\delta_n}=O(1),\quad
 \frac{\delta_{n-1}}{\delta_n}\to1 \quad\text{and}\quad
 \frac{M_{n-1}}{M_n}\to1
 \text{ as }n\to\infty.
\]
For $\delta_n\le y\le h_n-4\delta_n$, the corresponding exponential profiles differ by at most a factor $e^{Ch_n}$. For $0<y\le\delta_n$, the linear estimate $\sinh(\pi y/\delta_j)\asymp y/\delta_j$, with $j=n-1,n$, gives the same comparison. Together with the tip and transition estimates above, this proves comparability on the two sides for every $0<y<h_n$. The lower-wall argument follows by reflection.
\end{proof}

\subsection{One subdivision for the two sides}\label{subsec:subdivision}
We choose a common set of vertices on the two sides of each slit so that every edge has $\tau_U$-length bounded above and below by fixed positive constants on both sides. We then estimate the edge lengths to prove summability after pullback to the physical spiral.

\begin{figure}[H]
  \centering
  \includegraphics[width=0.7\linewidth]{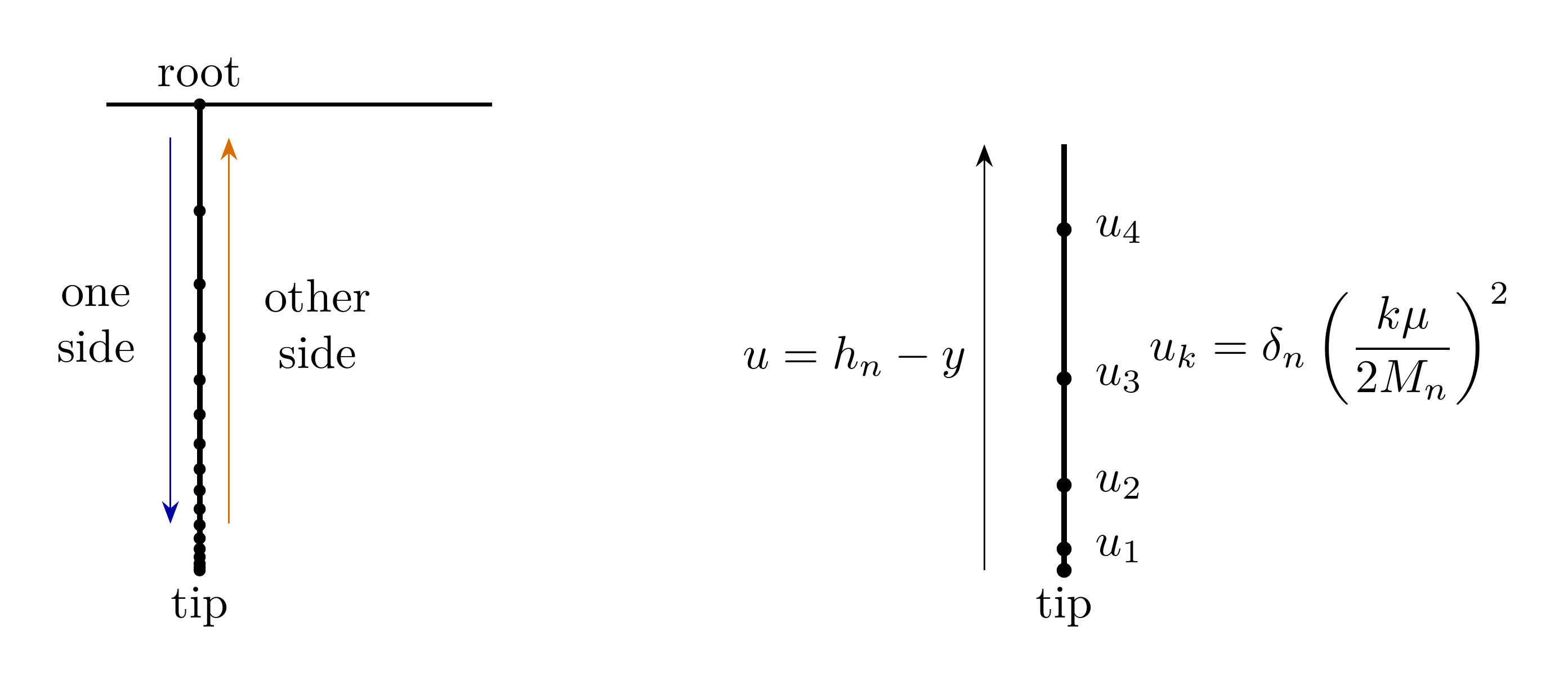}
  \caption[Subdivision of a slit and its tip region]{%
    The arrows trace the two sides of the slit through the same
    vertices in opposite directions. Each edge has density mass
    $\mu$, and its $\tau_U$-length on each side lies in $[c,C]$.
    The right panel magnifies the tip region, where the distance
    $u_k$ from the tip to the $k$th vertex is proportional to $k^2$.}
  \label{fig:shared-vertices}
\end{figure}

On the slit rooted at $x_n+\pi i/2$, define the positive density
\begin{equation}\label{eq:density}
 \varrho(y)=
 \begin{cases}
 \delta_n^{-1}e^{\pi y/\delta_n},&0\le y\le h_n-\delta_n,\\
 M_n/\sqrt{\delta_n(h_n-y)},&h_n-\delta_n<y<h_n.
 \end{cases}
\end{equation}
By~\eqref{eq:height}, its total mass is
\[
 J=\int_0^{h_n}\varrho(y)\,dy
 =M_n\left(2+\frac{e^{-\pi}}\pi\right)-\frac1\pi\asymp M_n.
\]
Set $N=\lceil J\rceil$ and divide the slit into $N$ consecutive edges, each of density mass $\mu=J/N\in[1/2,1]$ on a far tail.
The two boundary traversals and the tip spacing are illustrated in Figure~\ref{fig:shared-vertices}.

Since $1/2\le\mu\le1$, the first edge has depth
\begin{equation}\label{eq:first-edge-depth}
 y_1=\frac{\delta_n}{\pi}\log(1+\pi\mu)\asymp\delta_n \quad\text{and}\quad y_1<\delta_n.
\end{equation}
Thus the root estimate $|\tau_U'|\asymp y/\delta_n^2$ from~\eqref{eq:side-der} applies throughout $0<y\le y_1$ on either side. Its conformal length therefore satisfies
\[
 \int_0^{y_1}|\tau_U'|\,dy
 \asymp\int_0^{y_1}\frac{y}{\delta_n^2}\,dy
 =\frac{y_1^2}{2\delta_n^2}\asymp1.
\]
On every remaining edge, $y_1\le y<h_n$ and $y_1\asymp\delta_n$, so~\eqref{eq:side-der}, \eqref{eq:tip-der} and~\eqref{eq:transition-derivative} give $|\tau_U'|\asymp\varrho(y)$ on both sides. Its conformal length is therefore comparable to $\mu\asymp1$. Together with~\eqref{eq:trunk-size}, this places all slit and trunk side lengths in one fixed interval $[c,C]$.

For edges contained in $0\le y\le h_n-\delta_n$, consecutive lengths are comparable to $\delta_n/(j+1)$, where $j\ge0$ counts edges from the root. In the tip region $0\le u\le\delta_n$, vertices counted backward from the tip have distances
\begin{equation*}
 u_k=\delta_n\left(\frac{k\mu}{2M_n}\right)^2
 \quad\text{for}\quad  k=0,1,\ldots,\left\lfloor\frac{2M_n}{\mu}\right\rfloor.
\end{equation*}
The at most one edge crossing $y=h_n-\delta_n$ has length comparable to $\delta_n/M_n$, as do its neighboring edges.
This shows that adjacent lengths are comparable and, for every $s>1$,
\begin{equation*}
 \sum_{e\text{ in one strip slit}}(\diam e)^s
 \le C_s\delta_n^s\left(\sum_{j\ge1}j^{-s}+M_n^{1-s}\right)
 \le C_s\delta_n^s.
\end{equation*}
Pullback by $\psi$ multiplies lengths in that slit neighborhood by a factor comparable to the physical width. Since $\eta_n\delta_n\asymp\Delta_n$ in both models,
\begin{equation}\label{eq:physical-cost}
 A_s:=\sum_{e\text{ in the physical graph}}(\diam e)^s<\infty
 \quad\begin{cases}s>1,&\text{narrowing spiral},\\s>1+1/(q+1),&\text{logarithmic spiral}.
 \end{cases}
\end{equation}
Here and below the sum is over geometric edges; summing over boundary sides only multiplies it by at most two.

\subsection{Column intervals and a global bound}\label{subsec:columns}
We first locate a whole boundary column in the $\tau_U$ coordinate and estimate its length. We then bound $Q=\Rea\tau_U$ throughout $U$. These estimates will be used for boundary covering and for the growth bound, respectively.

A \emph{column} follows one side of a slit from its root to its tip, returns along the other side, and then follows the next trunk edge. The two slit sides are distinct boundary sides, although they lie on the same geometric segment. For either wall, let $\ell_n$ denote the Euclidean length of the interval traced by the $n$th column under $\tau_U$ on the imaginary axis. Thus $\ell_n$ includes the image lengths of both slit sides and the following trunk edge.

\begin{lemma}[Column intervals and a global bound]\label{lem:column-estimates}
With the normalization of Lemma~\ref{lem:central}, the following hold.
\begin{enumerate}[label=(\arabic*),leftmargin=*]
\item For the logarithmic model and all sufficiently large $n$, on either wall,
\begin{equation}\label{eq:column-scale}
 \ell_n\asymp M_n\quad\text{and}\quad|\Ima\zeta|\asymp e^{x_n},
\end{equation}
uniformly for $\zeta$ in the image interval of the $n$th column. The comparison constants are independent of $n$ and $\zeta$.
\item For either model,
\begin{equation}\label{eq:Qglobal}
 0<Q(w)\le e^x\cos y\le e^x
 \quad\text{for}\quad w=x+iy\in U.
\end{equation}
\end{enumerate}
\end{lemma}
\begin{proof}
\emph{(1) The column interval.}

For the length, Section~\ref{subsec:subdivision} divides the slit into $N\asymp M_n$ geometric edges, each with $\tau_U$-length comparable to $1$ on either side. The next trunk edge also has $\tau_U$-length comparable to $1$ by~\eqref{eq:trunk-size}. The boundary traversal is monotone on the imaginary axis, so these lengths add and give $\ell_n\asymp M_n$.

For the location, first consider an upper-wall column. Lemma~\ref{lem:entrance-value} and~\eqref{eq:central} give
\[
 \tau_U(C_n)=i e^{x_n+o(1)}+O(M_n)
 \text{ as }n\to\infty.
\]
Let $p=x_n+i(\pi/2-h_n)$. The interior gradient estimate and the bounds established in the proof of Lemma~\ref{lem:tip} give, by integration along the interior path and in the square-root coordinate,
\[
 |\tau_U(p)-\tau_U(C_n)|\le CM_n.
\]
Every point of the column interval is within distance $\ell_n=O(M_n)$ of $\tau_U(p)$. Thus, uniformly on that interval,
\[
 \zeta=i e^{x_n+o(1)}+O(M_n)\quad\text{and}\quad
 \frac{M_n}{e^{x_n}}=\delta_n\to0
 \text{ as }n\to\infty.
\]
This proves its location in~\eqref{eq:column-scale}. The lower-wall argument gives the corresponding interval on the negative imaginary axis.

\emph{(2) The global bound for $Q$.}
Fix any $w=x+iy\in U$. The real axis lies in $U$, so we may take real $T\to+\infty$ as a comparison point. Since $U\subset S$, conformal invariance and domain monotonicity give
\[
 \rho_S\bigl(\Phi(w),\Phi(T)\bigr)
 =\rho_U(w,T)\ge\rho_S(w,T).
\]
Lemma~\ref{lem:central} gives $\Phi(T)=T+o(1)$ as $T\to+\infty$. Hence the strip-distance formula~\eqref{eq:strip-hyperbolic-distance} yields
\[
 2e^{-T}\cosh\bigl(2\rho_U(w,T)\bigr)\to\frac1{Q(w)},\quad 2e^{-T}\cosh\bigl(2\rho_S(w,T)\bigr)\to\frac1{e^x\cos y},
\]
as $T\to+\infty$. Apply $\cosh(2\,\cdot)$ to the distance inequality, multiply by $2e^{-T}$, and let $T\to+\infty$. This gives
\[
 \frac1{Q(w)}\ge\frac1{e^x\cos y}.
\]
Since $Q(w)>0$ and $\cos y>0$ for $w=x+iy\in U\subset S$, taking reciprocals proves~\eqref{eq:Qglobal}.
\end{proof}

\section{Geometry of the graph and logarithmic area}\label{sec:graph-geometry}

Let \(T_o\) be the closed physical graph consisting of \(\partial\Omega\), including its cap point, and the images of the subdivided slits under \(\psi\). Its complementary components are the outside domain \(L\) and the slit corridor \(R_o=\psi(U)\). A slit root is its attachment point on the trunk; a tip is its free endpoint. The graph is connected, unbounded and locally finite, meaning that every compact subset of the plane meets finitely many edges.

For a graph \(T\) with edges \(e\) and a parameter \(r>0\), define its variable-scale neighborhood by
\begin{equation}
T(r)=\bigcup_e\{z:\dist(z,e)<r\diam(e)\}.
\label{eq:tube}
\end{equation}
The neighborhood radius depends on the edge size. This is the region where folding and the geometric correction will be supported.

\subsection{Comparison of nearby edges}\label{subsec:nearby-edges}

The geometry needed here includes comparable incident lengths, separation of nonincident edges relative to their sizes, uniformly bilipschitz vertex stars and bounded curvature of individual edges. The first property is contained in a stronger property
\begin{equation}
\dist(e,e')\le r\diam(e)
\quad\Longrightarrow\quad
\diam(e')\asymp_r\diam(e)
\label{eq:nearby-edges}
\end{equation}
for every fixed \(r\).

\begin{lemma}[Nearby edges on one strip wall]\label{lem:strip-nearby-edges}
The subdivision of Section~\ref{subsec:subdivision}, with \(|m'|\le Am\), satisfies \eqref{eq:nearby-edges} on either individual wall of the strip.
\end{lemma}

\begin{proof}
Translate or reflect the wall to \(y=0\), with slits pointing into \(y>0\). On one slit set \((\delta,h,M):=(\delta_n,h_n,M_n),\quad b_*=\delta/M^2\) and define the edge scale
\[
L_*(y)=
\begin{cases}
\delta e^{-\pi y/\delta},&0\le y\le h-\delta,\\
b_*+\sqrt{b_*(h-y)},&h-\delta<y\le h.
\end{cases}
\]
\(L_*(y)\) can be discontinuous, but the two expressions at the transition \(y=h-\delta\) are comparable to \(\delta/M\). Recall that every subdivision edge \(I\) has \(\varrho\)-mass \(\mu\in[1/2,1]\). 

\textbf{Step 1:} We claim that \(|I|\asymp L_*(y)\) for every \(y\in I\), with constants independent of the edge and the slit.

If \(I=[v,w]\) lies in the ordinary part, its mass equation gives
\[
\mu=\int_v^w\varrho(y)\,dy
=\int_v^w\frac1\delta e^{\pi y/\delta}\,dy
=\frac1\pi\left(e^{\pi w/\delta}-e^{\pi v/\delta}\right), 
\]
thus,
\[
\frac{\varrho(w)}{\varrho(v)}
=1+\pi\mu e^{-\pi v/\delta}\le1+\pi.
\]
Therefore, \(\varrho\) varies by a bounded factor on \(I\), and integration gives
\[
\frac{|I|}{1+\pi}\varrho(y)\le\mu\le(1+\pi)|I|\varrho(y).
\]
Since \(\frac12\le\mu\le1\), it follows that \(|I|\asymp1/\varrho(y)=L_*(y)\) throughout the edge.

For an edge \(I=[v,w]\) in the tip part, let \(u_1=h-w<u_2=h-v\) be the distances of its endpoints from the tip. By \eqref{eq:density}, 
\[
\varrho(y)=\frac{M}{\sqrt{\delta(h-y)}}=\frac1{\sqrt{b_*u}}\quad\text{and}\quad u=h-y,
\]
and the mass equation gives
\[
\mu=\int_{u_1}^{u_2}\frac{du}{\sqrt{b_*u}}=\frac2{\sqrt{b_*}}\left(\sqrt{u_2}-\sqrt{u_1}\right).
\]
Therefore, \(\sqrt{u_2}-\sqrt{u_1}=\frac{\mu\sqrt{b_*}}2\). By calculation,
\[
|I|=u_2-u_1
=\mu\sqrt{b_*u_1}+\frac{\mu^2b_*}{4}.
\]
Since \(\frac12\le\mu\le1\), we have \(|I|\asymp b_*+\sqrt{b_*u_1}\). Together with \(u_1\le u\le u_2\), we obtain
\[
b_*+\sqrt{b_*u_1}\le b_*+\sqrt{b_*u}=L_*(y)
\le b_*+\sqrt{b_*u_1}+\frac{b_*}{2}
\le\frac32\left(b_*+\sqrt{b_*u_1}\right).
\]
Therefore, \(\frac{1}{24}|L_*(y)|\le |I|\le |L_*(y)|\). Hence the comparison holds at every point, including the tip when it is an endpoint.

For the possible edge \(I=[v,w]\) crossing \(y_0=h-\delta\): when \(y\in[v,h-\delta]\), we have \(\varrho\asymp M/\delta\) and \(L_*\asymp\delta/M\) as before. When \(y\in[h-\delta,w]\), the inequality \(h-y\le \delta\) implies \(\varrho(y)=\frac M{\sqrt{\delta(h-y)}}\ge\frac M\delta\). Thus,
\[
\frac M\delta(w-y_0)\le\int_{y_0}^w\varrho(y)\,dy\le1.
\]
Since \(0\le v\le h-\delta\), we have \(M=e^{\pi h/\delta}>e^\pi>2\). Therefore, for any \(y\in[h-\delta,w]\), we have \(\frac\delta2\le h-y\le\delta\), and thus
\[
\varrho(y)\asymp\frac M\delta, \quad \text{and}\quad 
L_*(y)=\frac{\delta}{M^2}+\sqrt{\frac{\delta}{M^2}(h-y)}\asymp\frac\delta M.
\]
On the whole interval \(I\), \(\mu=\int_I\varrho(y)\,dy\asymp |I|\frac M\delta.\) Therefore, for all \(y\in I\), we have \(|I|\asymp\frac\delta M\asymp L_*(y)\). This proves the claim in this case as well.

\textbf{Step 2:} For any fixed \(c>0\), we next show that \(L_*(y')\asymp_c L_*(y)\) whenever \(|y-y'|\le cL_*(y)\), with comparison constants independent of the slit. 

If both points lie in the ordinary part, then \(L_*(y)\le\delta\), and
\[
\left|\log\frac{L_*(y')}{L_*(y)}\right|=\frac{\pi|y-y'|}{\delta}\le \frac{\pi cL_*(y)}{\delta}\le\pi c.
\]
Thus the ratio of the two scales is bounded above and below in terms of \(c\) alone.

If both points lie in the tip part, let \(s=\sqrt{(h-y)/b_*}\) and \(t=\sqrt{(h-y')/b_*}\).
Then \(L_*(y)=b_*(1+s)\), \(L_*(y')=b_*(1+t)\), and the proximity condition becomes
\[
|s^2-t^2|\le c(1+s).
\]
When \(s+t\ge1\), division by \(s+t\) gives
\[
|s-t|\le c\frac{1+s}{s+t}\le2c.
\]
When \(s+t<1\), it is easy that \(|s-t|<1\). Since \(1+s\) and \(1+t\) are both at least \(1\), they are comparable, proving the assertion on the tip part.

If the points lie on opposite sides of \(y_0=h-\delta\), first apply the comparison on the side containing \(y\), using \(|y-y_0|\le|y-y'|\le cL_*(y)\). If \(y\) in the ordinary part, \(L_*(y_0)\asymp_c L_*(y)\) as before. If \(y\) in the tip part,
\[
L_*(y)\asymp_c L_*(y_0^+)\asymp L_*(y_0^-)= L_*(y_0),
\]
where \( L_*(y_0^\pm)=\lim_{y\to y_0^\pm}L_*(y)\). Consequently,
\[
|y'-y_0|\le |y'-y|\le cL_*(y)\le C_c L_*(y_0).
\]
Apply the same argument on the other side gives
\[
L_*(y')\asymp_c L_*(y_0)\asymp_c L_*(y).
\]

Finally, let \(I,I'\) be subdivision edges of the same slit with \(\dist(I,I')\le r|I|\), and choose closest points with heights \(y\in I\) and \(y'\in I'\). The preceding estimate gives
\[
|y-y'|=\dist(I,I')
\le r|I|\le CrL_*(y),
\]
where \(C\) is a fixed constant. Applying the assertion with \(c=Cr\), we obtain
\[
|I'|\asymp L_*(y')\asymp_c L_*(y)\asymp |I|.
\]
Since these edges are straight segments, their lengths equal their diameters. This proves \eqref{eq:nearby-edges} for edges on the same slit.

\textbf{Step 3:}
For different slits, fix a small window of width \(R\) independent of \(n\) around a sufficiently far root \(x_n\). Since \(m(x_j)=\delta_j\), integrating the inequality \(|m'|\le Am\) gives
\[
\left|\log\frac{\delta_j}{\delta_n}\right|\le A|x_j-x_n|\le AR,
\]
Thus all base gaps \(\delta_j\) in this window are comparable to \(\delta_n\). In particular, distinct roots are separated by at least \(c_0\delta_n\). Meanwhile, \(|I|\le C\delta_n\) for every slit edge based in the window. Indeed, suppose \(I\) lies on the slit with root \(x_j\). Then
\[
|I|\le CL_*(y)\le C\max\left\{\delta_j e^{-\pi y/\delta_j},\frac{\delta_j}{M_j^2}
+\sqrt{\frac{\delta_j}{M_j^2}(h_j-y)}\right\}\le C\max\{\delta_j,2\delta_j \}\le C'\delta_n.
\]

Now let \(e\) and \(e'\) lie on different slits, rooted at \(x_n\) and \(x_j\), and suppose that \(\dist(e,e')\le r|e|\). It follows that \(c_0\delta_n\le |x_j-x_n|\le \dist(e,e')\le r|e|,\) thus, \(c_r\delta_n\le |e|\le C\delta_n\). This lower bound forces \(e\) to stay near its root. Because on the \(n\)-th slit the tip scale satisfies \(L_*(y)\le \frac{\delta_n}{M_n^2}+\frac{\delta_n}{M_n}=o(\delta_n),\) and \(|e|\asymp L_*(y)\) at every point in \(e\), the edge cannot meet the tip part once \(n\) is sufficiently large. On the ordinary part, the estimate \(c_r\delta_n\le |e|\le C\delta_n e^{-\pi y/\delta_n}\) implies \(y\le C_r\delta_n\) on \(e\).

Choose closest points of \(e\) and \(e'\), with depths \(y\) and \(y'\). Hence
\[
y'\le y+\dist(e,e')
\le y+r|e|
=O_r(\delta_n).
\]
Because \(\delta_j\asymp\delta_n\) and \(h_j/\delta_j\to\infty\), this point also lies in the ordinary part on a sufficiently far tail. It follows that
\[|e'|\asymp L_*(y')= \delta_j e^{-\pi y'/\delta_j}\asymp_r \delta_n\asymp_r |e|.
\]
 
For a slit--trunk pair, write \(I\) for the slit edge and \(J\) for the trunk edge. In the same window, \(|I|\le C\delta_n\) and \(|J|\asymp\delta_n\), so either \(\dist(I,J)\le r|I|\) or \(\dist(I,J)\le r|J|\) implies \(\dist(I,J)\le C_r\delta_n\). Since \(J\) lies on \(y=0\), the lower endpoint of \(I\) has depth at most \(C_r\delta_n\). On a sufficiently far tail this endpoint lies in the ordinary part, where \(|I|\asymp_r\delta_n\asymp|J|\) as before. 

Two trunk edges are compared directly using the estimate \(\delta_j\asymp\delta_n\).

\textbf{Step 4:}
Fix \(r>0\). Along the tail, both the base gaps and the maximal slit-edge lengths tend to zero. To locate the bases of two edges satisfying \(\dist(e,e')\le r|e|\), choose closest points on the two edges and project them onto the wall. Thus the distance between any two relevant roots or trunk endpoints is at most
\[
|e|+\dist(e,e')+|e'|
\le (r+1)|e|+|e'|.
\]
For fixed \(r\), this bound is uniformly smaller than the radius of the chosen window once both edges lie sufficiently far along the tail. All relevant objects therefore lie in one such window, and the preceding estimates apply with constants depending only on \(r\).

It remains to consider pairs consisting an edge from the finite initial part. Since this part consists of finitely many bounded edges and the tail edge lengths tend to zero, all edge lengths are bounded by some \(D<\infty\). Let \(K\) be the union of the initial edges. Regardless of which edge belongs to \(K\), the condition gives
\[
\dist(e,e')\le r|e|\le rD.
\]
Consequently, the other edge meets the closed \(rD\)-neighborhood of \(K\), which is compact. Local finiteness implies that only finitely many edges meet this neighborhood, so there are only finitely many remaining pairs to compare. Since their lengths are positive, both ratios \(|e'|/|e|\) and \(|e|/|e'|\) have a finite maximum over these pairs. This proves \eqref{eq:nearby-edges} for all the edges under consideration.
\end{proof}

\begin{lemma}[Transfer to the physical graph]\label{lem:physical-nearby-edges}
Both physical spiral graphs satisfy \eqref{eq:nearby-edges}.
\end{lemma}

\begin{proof}
Fix \(r>0\). We first prove that two physical edges sufficiently far satisfying \(\dist(e,e')\le r\diam(e)\) can be placed entirely into the same local coordinate neighborhood of the boundary. Thus, after pulling them back to the strip, we can apply Lemma~\ref{lem:strip-nearby-edges} to these two preimage edges.

At a root with strip coordinate \(x_n\) and physical scale \(\eta_n\), the boundary charts of Lemmas~\ref{lem:narrow-lip} and~\ref{lem:log-scale}, together with the estimates \(|I|\le C\delta_n\), give
\[
\diam(\text{whole slit})
\le C\eta_n h_n=o(\eta_n)\quad \text{and} \quad \diam(e)\le C\eta_n\delta_n=o(\eta_n).
\]
Thus every whole slit, as well as every trunk edge, is small compared with the local corridor width. By Section~\ref{sec:models}, the opposite wall is separated from the chosen wall by a positive multiple of \(\eta_n\). Different turns have a fixed positive separation in the narrowing model; and their separation is comparable to the physical radius in the logarithmic model, which is much larger than \(\eta_n\). Consequently, a disk whose radius is a sufficiently small fixed multiple of \(\eta_n\), centered at a far root, meets only one local wall and is covered by one of the boundary charts.

We also use the fact that \(\diam(\text{whole slit})\) tend to zero. In the narrowing model this follows from \(\eta_n\to0\) and \(h_n\to0\). In the logarithmic model, \eqref{eq:radius-width},\eqref{eq:outside-mesh}--\eqref{eq:log-small} and the definition~\eqref{eq:height} of \(h_n\) give
\[
\eta_n h_n\le C\Delta_n x_n=n^{1/(q+2)-1+o(1)}(\log n)^{q+o(1)}\longrightarrow0.
\]
The diameters of edges on physical trunk also tend to zero, since \(\eta_n\delta_n\asymp\Delta_n\to0\).

Now let \(e,e'\) be tail edges of the physical graph satisfying \(\dist(e,e')\le r\diam(e).\)
For a slit edge, define its base to be the root of the whole physical slit containing it. For a trunk edge, define its base to be the image under \(\psi\) of the left endpoint of its strip preimage. Let \(b_n,b_j\) denote the bases, \(x_n,x_j\) denote the real parts of the strip coordinates of the bases of \(e,e'\), respectively \(\bigl(d_n:=\diam(e\cup\{b_n\}),\ d_j:=\diam(e'\cup\{b_j\})\bigr)\).

We claim that the base of \(e'\) and the entire edge \(e'\) lie in the same \(o(\eta_n)\)-neighborhood. In the narrowing model, Choose closest points \(z\) and \(z'\) of \(e\) and \(e'\), respectively. Then,
\[|b_j-b_n|\le |b_n-z|+|z-z'|+|z'-b_j|\le d_n+rd_n+d_j=(r+1)d_n+d_j.\]
Since \(d_n\to0\), the two bases are at distance \(o(1)\). By definition, \(b_n=\Gamma_\pm(t_n), b_j=\Gamma_\pm(t_j).\)
Since \(|\Gamma_\pm(t)|=t+o(1),\) their parameters satisfy
\[|t_j-t_n|\le \bigl|t_j-|b_j|\bigr|
   +\bigl||b_j|-|b_n|\bigr|
   +\bigl||b_n|-t_n\bigr|\le e^{-t_j}+|b_j-b_n|+e^{-t_n}\to 0.
\]
Therefore \(\frac{\eta_j}{\eta_n}=e^{-(t_j-t_n)}\to 1.\) The whole slit or trunk edge associated with \(e'\) has diameter \(o(\eta_j)=o(\eta_n)\). Then, 
\[
\sup_{w\in e'\cup\{b_j\}}|w-b_n|\le  d_j+rd_n+d_n=o(\eta_n),
\]
implies the claim. In the logarithmic model, the physical diameters of all whole slits and trunk edges are uniform bounded by \(D\). By \eqref{eq:params} and~\eqref{eq:radius-width}, \(\eta_n\asymp\frac{r_n}{(\log r_n)^p}\to\infty.\) Thus, \(\frac{D}{\eta_n}\to 0.\)
For any \(w\in e'\cup\{b_j\}\),
\[
|w-b_n|\le |w-z'|+|z'-z|+|z-b_n|\le D+rD+D=(r+2)D=o(\eta_n).
\]
This implies the claim.

In both models, the local boundary geometry established in Lemmas~\ref{lem:narrow-lip} and~\ref{lem:log-scale} separates the two walls at scale \(\eta_n\). Moreover, distinct turns remain a fixed positive distance apart in the narrowing model, while their separation is comparable to the physical radius in the logarithmic model, as follows from the boundary asymptotics in the proof of Lemma~\ref{lem:outside-log}. Since both edges lie within \(o(\eta_n)\) of the same physical base, they must be associated with the same local wall.

On this chart, Lemmas~\ref{lem:narrow-lip} and~\ref{lem:log-scale} (\eqref{eq:narrow-wall-chart} and~\eqref{eq:log-wall-chart}) give a constant \(K\ge1\), independent of the root, such that
\[
K^{-1}\eta_n|u-v|\le |\psi(u)-\psi(v)|\le K\eta_n|u-v|.
\]
Write \(e=\psi(I)\) and \(e'=\psi(I')\), where \(I,I'\) are strip edges. It follows that
\[
\dist(I,I')\le K\eta_n^{-1}\dist(e,e')\le Kr\eta_n^{-1}\diam(e)\le K^2r\diam(I).
\]
Applying Lemma~\ref{lem:strip-nearby-edges}, and then using the chart estimates again, gives
\[
\diam(e')
\asymp \eta_n\diam(I')
\asymp_r \eta_n\diam(I)
\asymp \diam(e).
\]
This explains why the physical scale \(\eta_n\) cancels in the comparison.

Finally, all physical edge diameters are bounded by some \(D<\infty\), because they tend to zero along the tails and the initial part has only finitely many edges. If a nearby pair contains an initial edge, the other edge meets the closed \(rD\)-neighborhood of the finite initial part. This neighborhood is compact, so local finiteness leaves only finitely many such pairs. Enlarging the comparison constant to include their positive diameter ratios proves \eqref{eq:nearby-edges} for the entire physical graph.
\end{proof}

\textbf{Separation of nonincident edges.} We first work in strip coordinates on one wall. We show that there is a fixed constant \(c_{\mathrm{sep}}>0\) such that any two nonincident edges satisfy
\[
\dist(e,e')
\ge c_{\mathrm{sep}}
\max\{\diam(e),\diam(e')\}.
\]

If two nonincident edges lie on the trunk or the same slit, at least one complete edge lies between them. By \eqref{eq:consec-mesh} and Section~\ref{subsec:subdivision}, incident edge lengths are comparable, thus the edges next to \(e\) and \(e'\) have lengths comparable to \(|e|\) and \(|e'|\), respectively. \(\dist(e,e')\) is therefore bounded below by a fixed multiple of the larger edge length. 

For edges on different slits in a fixed tail window, the horizontal distance of their roots is at least \(c\delta_n\), while both edge lengths are at most \(C\delta_n\). The conclusion follows. For a slit edge and a nonincident trunk edge, either their horizontal projections are separated by at least one trunk edge, or the trunk edge has an endpoint at the slit root. In the latter case, the slit edge does not contain the root and is separated from it by at least the first subdivision edge on the slit, whose length is comparable to \(\delta_n\). Thus the required separation holds in both cases.

The argument in the first part of Lemma~\ref{lem:physical-nearby-edges} shows that any sufficiently far pair of edges satisfying \(\dist(e,e')\le r\diam(e)\) lies in a common boundary chart. Applying this conclusion with \(r=1\) in both orders, we see that two tail edges outside any common chart must satisfy \(\dist(e,e')\ge\max\{\diam(e),\diam(e')\}.\) Thus it remains to establish the separation estimate for edges lying in a common chart. By Lemmas~\ref{lem:narrow-lip} and~\ref{lem:log-scale}, the uniform bilipschitz boundary charts transfer the strip estimates to the physical graph, giving the required separation bound.

\medskip

\textbf{Uniformly bilipschitz vertex stars.} We next consider the geometry near each vertex. At a slit root, three edges meet: two trunk edges extending in opposite directions and the first slit edge pointing into the corridor. In strip coordinates, the slit meets the trunk at a right angle. Moreover, the two trunk edge lengths are comparable by \eqref{eq:consec-mesh}, and the first slit edge has length \(y_1=\frac{\delta_n}{\pi}\log(1+\pi\mu) \asymp\delta_n,\) as shown in~\eqref{eq:first-edge-depth}.

At the other subdivision vertex, the two incident edges lie on the same slit or trunk and have comparable lengths. Their union is locally a straight interval in strip coordinates, and hence a uniformly bilipschitz arc in the physical boundary chart. A slit tip has only one edge and requires no further argument. The finitely many initial vertices are handled by enlarging the constant.

By Lemmas~\ref{lem:narrow-lip} and~\ref{lem:log-scale}, all these vertex-star estimates hold uniformly on the physical graph.

\medskip

\textbf{Bounded curvature of individual edges.} 
In the narrowing model, the images of the strip slits under \(\psi\) are generally curved. By the local boundary charts of Lemma~\ref{lem:narrow-lip}, their physical curvature is \(O(\eta_n^{-1})\), which does not give a uniform bound since \(\eta_n\to0\). However, each subdivided tail edge \(e\) has length \(o(\eta_n)\). Indeed, the subdivision estimates of Section~\ref{subsec:subdivision} and the boundary charts of Lemma~\ref{lem:narrow-lip} give \(d_e\le C\eta_n\delta_n,\) where \(d_e\) is the arclength of \(e\). Hence the total turning of the tangent is at most \(O(d_e/\eta_n)=O(\delta_n)\to 0.\) Thus, sufficiently far along the tail, every slit edge is nearly straight.

Choose coordinates along the chord of such an edge, with its endpoints at \((0,0)\) and \((d'_e,0)\), where \(d'_e\) is the chord length. The edge can then be written as a graph \(y=p_e(x)\), \(0\le x\le d'_e\), satisfying
\[
p_e(0)=p_e(d'_e)=0,\quad
|p'_e(x)|\le\varepsilon_e,\quad\text{and}\quad
|p_e(x)|\le\varepsilon_e\min(x,d'_e-x),
\]
where \(\varepsilon_e\to0\) uniformly along the tail.

We now straighten these edges without changing their endpoints. The separation and vertex-star estimates established in Section~\ref{subsec:nearby-edges} allow us to choose thin diamond-shaped neighborhoods around the tail edges, whose interiors are pairwise disjoint and whose closures meet only at common endpoints. Using smooth cutoff functions, we construct local maps that straighten each edge onto its chord, fix its endpoints, and agree with the identity outside the corresponding neighborhood. The bounds on \(p_e\) and \(p'_e\) ensure that the derivatives of these local maps differ from the identity matrix by \(O(\varepsilon_e)\), uniformly over all edges.

Since \(\varepsilon_e\to0\), we may restrict the construction to a sufficiently far tail so that all these local maps are uniformly close to the identity. They combine with the identity outside their supports to give a global bilipschitz homeomorphism \(F\). We leave the finitely many initial edges unchanged. The map \(F\) fixes the trunk and \(L\), preserves \(\Omega\), and gives \(T=F(T_o)\) and \(R=F(R_o)\). All modified tail slit edges are now straight. The trunk has uniformly bounded curvature by its parametrization in \eqref{eq:narrow}, and the finitely many unchanged edges have bounded curvature. Hence the final graph has uniformly bounded curvature on every edge.

Since \(F\) is globally bilipschitz, the property \eqref{eq:nearby-edges}, the relative separation of nonincident edges, and the vertex-star bounds, are preserved with possibly different uniform constants. Together with the curvature bound obtained above, this proves all the required geometric properties of \(T\).

In the logarithmic model, no straightening is needed. The explicit parametrization \eqref{eq:log-model} and the derivative estimates in Lemma~\ref{lem:log-scale} already give uniformly bounded physical curvature of the tail edges. We set \(T=T_o\) and \(R=R_o\) for notational consistency.

\subsection{A uniform logarithmic-area estimate}

\begin{proposition}\label{prop:uniform-log-area}
For each fixed \(r>0\), and a sufficiently large fixed \(\Lambda\), both final graphs satisfy
\begin{equation}\label{eq:uniform-log-area}
\sup_{w\notin T(\Lambda r)}
\int_{T(r)}\frac{dA(z)}{|z-w|^2}<\infty.
\end{equation}
\end{proposition}

\begin{proof}
We first work with the graph \(T_o\) before straightening. Let \(l(s)\) be the length of the subdivision edge containing \(s\). Lemma~\ref{lem:strip-nearby-edges} gives \(l(z(y))\asymp\delta e^{-\pi y/\delta}\) if \(0\le y\le h-\delta\), and \(l(z(h-u))\asymp b_*+\sqrt{b_*u}\) if \(0\le u=h-y\le\delta\). Here a geometric column contains one slit and the following trunk edge; counting both slit sides changes only a constant. Since \(M\to\infty\),
\begin{equation}\label{eq:column-area-cost}
    \int_{\rm column}l(s)\,ds
 \le C\left(\delta^2+b_*\delta+\sqrt{b_*}\,\delta^{3/2}\right)
 \le C\delta^2.
\end{equation}
The first term also includes the trunk edge, whose length is comparable to \(\delta\).

\textit{Reduction to a line integral.}
For a straight edge \(e\) of length \(\ell\), its \(r\ell\)-neighborhood \(N_r(e)\) has area at most \(2r\ell^2+\pi r^2\ell^2\). If \(z\in N_r(e)\) and \(w\notin N_{\Lambda r}(e)\) with \(\Lambda\ge2\), then \(|z-w|\ge \dist(w,e)/2\), whereas
\[
\int_e\frac{\ell\,ds(\zeta)}{|\zeta-w|^2}
\ge\frac{\ell^2}{(\dist(w,e)+\ell)^2}.
\]
Consequently, writing \(d=\dist(w,e)\) and using \(d\ge2r\ell\), we obtain
\begin{align*}
\int_{N_r(e)}\frac{dA(\zeta)}{|\zeta-w|^2}
&\le 4(2r+\pi r^2)\frac{\ell^2}{d^2}
   =4(2r+\pi r^2)
   \frac{(d+\ell)^2}{d^2}
   \frac{\ell^2}{(d+\ell)^2}\\
&\le4(2r+\pi r^2)
   \left(1+\frac1{2r}\right)^2 \int_e\frac{\ell\,ds(\zeta)}{|\zeta-w|^2},
\end{align*}
Summing over the edges of a column, with \(l(\zeta)\) denoting the length of the edge containing \(\zeta\), bounds its contribution by
\begin{equation}
 C_r\int_{\mathrm{column}}
 \frac{l(\zeta)}{|\zeta-w|^2}\,|d\zeta|.
 \label{eq:column-line-integral}
\end{equation}
Uniformly bilipschitz charts give the same comparison for the corresponding physical arcs, with adjusted constants and a fixed enlargement of \(\Lambda\).

\textit{One slit and its trunk.}
Fix a column and a point \(w\) outside the enlarged neighborhoods \(N_{\Lambda r}(e)\) of all its subdivision edges, with \(\Lambda\ge2\). Let \(S_{\mathrm{ord}}\) be the closed ordinary segment of the slit, parametrized by height \(0\le y\le h-\delta\). Choose \(z_0=z(y_0)\in S_{\mathrm{ord}}\) nearest to \(w\). Then
\[
 |z(y)-w|^2
 \ge \dist(w,S_{\mathrm{ord}})^2+(y-y_0)^2\quad \text{and} \quad
 \dist(w,S_{\mathrm{ord}})
 \ge \Lambda r\,l(z_0).
\]
For \(|y-y_0|\le\delta\), by the proof of Lemma~\ref{lem:strip-nearby-edges}, \(l(z(y))\le Cl(z_0)\). On the remaining ordinary segment, \(|z(y)-w|^2\ge\delta^2\). Therefore
\[
 \begin{aligned}
 \int_{S_{\mathrm{ord}}}
 \frac{l(\zeta)}{|\zeta-w|^2}\,|d\zeta|
 &\le
 Cl(z_0)\int_{\R}
 \frac{dt}{\dist(w,S_{\mathrm{ord}})^2+t^2}
 +\frac{C}{\delta^2}
 \int_{S_{\mathrm{ord}}}l(\zeta)\,|d\zeta|\\
 &\le C_r\left(
 \frac{l(z_0)}{\dist(w,S_{\mathrm{ord}})}+1
 \right)
 \le C_r.
 \end{aligned}
\]

For the closed tip segment \(S_{\mathrm{tip}}\), choose its point \(z_1\) nearest to \(w\). Let \(y_1\) be the height of \(z_1\), and set \(u_1=h-y_1\). For a point \(\zeta\in S_{\mathrm{tip}}\) with \(u=h-y\), then
\[
|\zeta-w|^2\ge |z_1-w|^2+(u-u_1)^2\quad \text{and} \quad
l(z(h-u))\le C\bigl(l(z_1)+\sqrt{b_*|u-u_1|}\bigr),
\]
with \(l(z_1)\ge cb_*\) and \(|z_1-w|\ge c_r l(z_1)\). We obtain
\[
 \begin{aligned}
 \int_{S_{\mathrm{tip}}}
 \frac{l(\zeta)}{|\zeta-w|^2}\,|d\zeta|
 &\le C\int_{\R}
 \frac{l(z_1)+\sqrt{b_*|u-u_1|}}
 {|z_1-w|^2+(u-u_1)^2}\,du\\
 &\le C_r\left(
 \frac{l(z_1)}{|z_1-w|}
 +\sqrt{\frac{b_*}{|z_1-w|}}
 \right)
 \le C_r.
 \end{aligned}
\]
The square-root term follows by substituting \(u-u_1=|z_1-w|v\), since \( \int_{\R}\frac{|v|^{1/2}}{1+v^2}\,dv<\infty.\) By the proof of Lemma~\ref{lem:strip-nearby-edges}, any subdivision edge crossing \(y=h-\delta\) has length comparable to \(\delta/M\) and to both transition scales. Hence the integral may be split at \(h-\delta\).

Finally, let \(z_2\) be the point of the following trunk edge \(e_{\mathrm{tr}}\) nearest to \(w\).
Since \(l\asymp\delta\) on this edge and \(|z_2-w|\ge c_r\delta\), the same estimate gives
\[
 \int_{e_{\mathrm{tr}}}
 \frac{l(\zeta)}{|\zeta-w|^2}\,|d\zeta|
 \le C\delta\int_{\R}
 \frac{dt}{|z_2-w|^2+t^2}
 \le C_r\frac{\delta}{|z_2-w|}
 \le C_r.
\]
Together with \eqref{eq:column-line-integral}, these estimates bound the contribution of a column by \(C_r\).

\textit{Several columns and physical charts.}
We estimate the contribution of the tubes meeting a fixed small strip window containing \(w\). Use a slightly larger window to include the columns producing these tubes. This is possible on the far tail, since the slit heights, trunk lengths, and tube radii tend to zero. By \eqref{eq:mesh-reg}, all root gaps in the larger window are comparable to a common \(\delta\).

First consider columns whose base coordinates satisfy \(|x_n-\Rea w|\le A\delta\), where \(A\) is a sufficiently large fixed constant. Since adjacent roots are separated by at least \(c\delta\), the number of these columns is uniformly bounded. The preceding single-column estimate therefore bounds their total contribution by \(C_r\). For the remaining columns, order these columns by increasing \(|x_n-\Rea w|\). Since adjacent roots are separated by at least \(c\delta\) by \eqref{eq:mesh-reg},
\[
 \dist(w,\text{the \(j\)th remaining column})\ge cj\delta,
 \qquad j\ge1.
\]
By \eqref{eq:column-area-cost} and \eqref{eq:column-line-integral}, their total contribution is at most
\[
 C_r\sum_{j\ge1}\frac{\delta^2}{(j\delta)^2}
 =C_r\sum_{j\ge1}\frac1{j^2}
 \le C_r.
\]
Thus the tubes meeting the smaller window have uniformly bounded total contribution.

By Lemmas~\ref{lem:narrow-lip} and~\ref{lem:log-scale}, the normalized wall charts are uniformly bilipschitz. Hence the preceding strip estimate transfers to physical neighborhoods of radius \(c\eta\), with uniform constants after a fixed enlargement of \(\Lambda\). The wall separation established in the proof of Lemma~\ref{lem:physical-nearby-edges} bounds the number of wall charts needed.

\textit{Estimates in a fixed physical neighborhood.}
Choose a fixed \(r_*>0\) sufficiently small that the enlarged disk \(D(w,10r_*)\) meets at most the two walls of one turn on the far tail and lies within the reflected outside boundary windows. For an observation point \(w\) satisfying \(\dist(w,\text{tail graph})<r_*\), choose a nearest point of the tail graph to \(w\). If it lies on a wall, denote it by \(p\); if it lies on a slit, let \(p\) be the root of that slit. Write \(p=\Gamma_\pm(t)\) and set \(\eta_*=\eta(t)\). Throughout the enlarged disk the wall parameters differ by \(O(r_*/t)\), since \(|\Gamma_\pm'(t)|\sim t\) as in the proof of Lemma~\ref{lem:narrow-scale}. So their widths are comparable to \(\eta_*\). Write \(\Delta^+\) and \(\Delta^-\) for the local physical trunk-edge scales on the two walls, respectively.

By \eqref{eq:narrow-der}, \(c\le|\sigma_L'|\le C\). Along a wall arc of length \(O(r_*)\), integration gives
\[
 |\Rea\sigma_L(p)-\Rea\sigma_L(q)|\le  Cr_*.
\]
Since \(|\tau_L'|=e^{\Rea\sigma_L}|\sigma_L'|\), the values of \(|\tau_L'|\) on this enlarged arc have bounded ratio. By the boundary partition chosen in Section~\ref{subsec:narrow-model}, each trunk edge has \(\tau_L\)-length exactly \(2\pi\), so its physical length is comparable to a common scale \(\Delta^\pm\) on that wall.

Fix a sufficiently small \(c>0\) as allowed by the local charts, and let \(c\eta_*\le r_{\mathrm d}\le r_*\). A column on either wall has physical tube area at most \(C_r\eta_*^2\delta_n^2\le C_r(\Delta^\pm)^2\), by the estimate above and \eqref{eq:narrow-mesh}. By the localization argument in the proof of Lemma~\ref{lem:physical-nearby-edges} and the preceding width comparison, any column whose tube meets \(D(w,r_{\mathrm d})\) has a root \(p'\) satisfying
\[
 |p'-w|\le r_{\mathrm d}+C_r\eta_*(h_n+\delta_n)\le2r_{\mathrm d}
\]
on the far tail, since \(h_n+\delta_n\to0\) and \(r_{\mathrm d}\ge c\eta_*\). By the boundary parametrization following \eqref{eq:narrow}, the wall arc in \(D(w,2r_{\mathrm d})\) has length \(O(r_{\mathrm d})\). The number of its roots is therefore at most \(C(1+r_{\mathrm d}/\Delta^\pm)\). Consequently,
\begin{equation}\label{eq:local-tube-area}
 \area\bigl(T_o(r)\cap D(w,r_{\mathrm d})\bigr)\le C_r\sum_{\pm}\left(1+\frac{r_{\mathrm d}}{\Delta^\pm}\right)(\Delta^\pm)^2 \le C_r r_{\mathrm d}(\Delta^++\Delta^-).
\end{equation}
For the last inequality use \(\Delta^\pm=o(\eta_*)\) from \eqref{eq:narrow-mesh}--\eqref{eq:narrow-small} and \(r_{\mathrm d}\ge c\eta_*\). 

Let \(R_k=2^kc\eta_*\) and \(A_k=\{z:R_k\le|z-w|<2R_k\}.\) Choose \(K\) so that \(R_K<r_*\le2R_K\). Applying the disk bound~\eqref{eq:local-tube-area} to these annuli, with the fixed window enlarged to cover the last one, gives
\[
 \begin{aligned}
 \int_{T_o(r)\cap \left(D(w,r_*)\setminus D(w,c\eta_*)\right)} \frac{dA(z)}{|z-w|^2}
 &\le \sum_{k=0}^{K} \int_{T_o(r)\cap A_k}\frac{dA(z)}{|z-w|^2}\\
 &\le C_r\sum_{k=0}^{K} \frac{R_k(\Delta^++\Delta^-)}{R_k^2}
 \le C_r\frac{\Delta^++\Delta^-}{\eta_*}.
 \end{aligned}
\]
This is bounded and tends to zero on the tail.

If \(\dist(w,T_o)\le2c\eta_*\), the local wall-chart estimates above give 
\[
 \int_{T_o(r)\cap D(w,c\eta_*)} \frac{dA(z)}{|z-w|^2}\le C_r.
\]
Otherwise, this disk misses \(T_o(r)\) on the far tail, since the tube radii are \(o(\eta_*)\). Combining this bound with the annular estimates yields
\[
 \int_{T_o(r)\cap D(w,r_*)} \frac{dA(z)}{|z-w|^2}\le C_r.
\]

In the logarithmic model \(\eta(x)\to\infty\), so these local charts eventually cover any fixed physical neighborhood meeting the tail graph. No intermediate range is needed there. 

\textit{The distant part and straightening.}
By the tube definition \eqref{eq:tube} and \eqref{eq:physical-cost} with \(s=2\), the graph \(T_o\) in either model satisfies
\[
 \area T_o(r) \le \pi(1+r)^2\sum_e(\diam e)^2<\infty.
\]
In particular, for any fixed \(a>0\),
\[
 \int_{T_o(r)\setminus D(w,a)}\frac{dA(z)}{|z-w|^2} \le a^{-2}\area T_o(r),
\]
uniformly in \(w\). This handles the region beyond the fixed local radius. Centers located a fixed distance from the tail graph are handled in the same way, since the tail tube radii tend to zero. 

The remaining initial part lies in a fixed bounded region. Only finitely many edge tubes meet this region, by local finiteness and the uniform bound on edge diameters. For each of these edges, let \(N_r(e)\) denote its \(r\diam(e)\)-neighborhood. Since \(w\notin T_o(\Lambda r)\), \(|z-w|\ge(\Lambda-1)r\diam(e)\) for \(z\in N_r(e)\). Using \(\area N_r(e)\le\pi(1+r)^2(\diam e)^2\), we obtain
\[
 \int_{N_r(e)}\frac{dA(z)}{|z-w|^2} \le \frac{\pi(1+r)^2(\diam e)^2}{(\Lambda-1)^2r^2(\diam e)^2} \le C_r.
\]
Summing over the finitely many initial edges gives a bound independent of \(w\). Together with the preceding local and distant estimates, this proves \eqref{eq:uniform-log-area} for \(T_o\).

In the narrowing model, \(T=F(T_o)\), where \(F\) is globally bilipschitz by Section~\ref{subsec:nearby-edges}. By \eqref{eq:tube}, tube parameters change only by fixed factors under \(F\). Thus, for sufficiently large \(\Lambda\), the change of variables \(z=F(\zeta)\), with \(w_o=F^{-1}(w)\), gives
\[
 \int_{T(r)}\frac{dA(z)}{|z-w|^2} \le C_F\int_{T_o(C_Fr)} \frac{dA(\zeta)}{|\zeta-w_o|^2} \le C_r.
\]
In the logarithmic case \(T=T_o\). This proves \eqref{eq:uniform-log-area} for both final graphs.
\end{proof}

\subsection{Harmonic measure and conformal side lengths}

We first prove the following harmonic-measure estimate for both \(T_o\) and \(T\). The proof follows the Beurling argument in \cite[Lemma~3.2]{AlbrechtBishop}, using the edge-length and separation estimates established in Section~\ref{sec:graph-geometry}, which hold for both graphs.

\begin{lemma}[A uniform boundary side near the graph]\label{lem:uniform-harmonic-side}
For either model and each fixed \(r>0\), there is \(\varepsilon_r>0\) such that, whenever \(z\in T(r)\) lies in a complementary component \(D\) of \(T\), some boundary side \(\mathfrak{s}\) of an edge satisfies
\begin{equation}\label{eq:uniform-side-harmonic-measure}
 \omega_D(z,\mathfrak{s})\ge\varepsilon_r.
\end{equation}
The same statement holds with \(T\) replaced by \(T_o\). The neighborhoods are defined as in \eqref{eq:tube}. The constant can be chosen for both graphs and depends only on \(r\) and the fixed geometric constants, not on \(z\) or the edge indices.
\end{lemma}

\begin{proof}
The preceding geometric arguments give the following facts for \(T\). Its vertex degrees are at most \(v_*=3\). There is \(c_{\mathrm{sep}}>0\) such that nonincident edges satisfy
\begin{equation}\label{eq:edge-separation}
\dist(e,e')\ge c_{\mathrm{sep}}\max\{\diam(e),\diam(e')\}.
\end{equation}
For every fixed \(a>0\), choose \(C_a\ge1\) so that \eqref{eq:nearby-edges} means
\[
\dist(e,e')\le a\diam(e)
\quad\Longrightarrow\quad
C_a^{-1}\diam(e)\le\diam(e')\le C_a\diam(e).
\]
These facts also hold for \(T_o\) before applying \(F\).

Choose an edge \(e\) with \(\dist(z,e)<r\diam(e)\), and let \(p\) be a nearest point to \(z\) in \(T\). The graph is closed, so \(p\) exists and \(p\in\partial D\). Choose an edge \(e_0\) containing \(p\). The nearby-edge comparison gives
\begin{equation*}
C_r^{-1}\diam(e)\le \diam(e_0)\le C_r\diam(e),
\end{equation*}
with \(|z-p|\le rC_r\diam(e_0)\). 

Let \(C_B\) be the universal constant in the Beurling escape estimate used below. Fix \(A>\max\{1,2rC_r,4C_B^2rC_r\}.\) By \eqref{eq:nearby-edges}, every edge meeting \(D(p,A\diam(e_0))\) has diameter between \(C_A^{-1}\diam(e_0)\) and \(C_A\diam(e_0)\). The bounded vertex degree and the separation~\eqref{eq:edge-separation} of nonincident edges proved after Lemma~\ref{lem:physical-nearby-edges} therefore give a uniform bound:
\begin{equation*}
 \#\{e':e'\cap D(p,A\diam(e_0))\ne\varnothing\}\le N_r.
\end{equation*}
Here \(N_r\) depends only on \(r\) and the fixed geometric constants, not on \(p\) or \(\diam(e_0)\).

Because \(T\) is connected and unbounded, a graph path from \(p\) to a point outside \(D(p,A\diam(e_0))\) contains a compact connected set joining \(p\) to its boundary. By the Beurling projection theorem, as used in \cite[proof of Lemma~3.2]{AlbrechtBishop}, planar Brownian motion starting at \(z\) has probability at most
\[
 C_B\sqrt{\frac{|z-p|}{A\diam(e_0)}}
 \le C_B\sqrt{\frac{rC_r}{A}}
\]
of leaving \(D(p,A\diam(e_0))\) before hitting this connected set. Here Brownian hitting probabilities are the harmonic measures of the corresponding boundary sets \cite[Appendix~F]{GarnettMarshall}. Since this connected set is contained in \(T\), any Brownian path that leaves the disk without hitting \(T\) also avoids this set. Thus the probability of leaving the disk before hitting \(T\) is at most \(C_B\sqrt{rC_r/A}<1/2\), where the last inequality follows from our choice \(A>4C_B^2rC_r\). Hence
\begin{equation}\label{eq:local-boundary-measure}
\omega_D\bigl(z,\partial D\cap D(p,A\diam(e_0))\bigr)\ge\tfrac12.
\end{equation}

The boundary part in \eqref{eq:local-boundary-measure} belongs to at most \(N_r\) geometric edges. Each edge has at most two boundary sides. Vertices have zero harmonic measure. Thus one full side has harmonic measure at least \(1/(4N_r)\). Set \(\varepsilon_r=1/(4N_r)\). This proves \eqref{eq:uniform-side-harmonic-measure}. All choices depend only on the metric constants and \(r\), so no absolute curvature bound has entered the proof. The same proof applies to \(T_o\), since it satisfies the same metric conditions.
\end{proof}

For a conformal map \(\tau:D\to\HH\) whose boundary-side intervals have length at most \(C_*\), Lemma~\ref{lem:uniform-harmonic-side} and \eqref{eq:poisson-half-plane} immediately give a coordinate-strip bound. If \(\tau(z)=x+iy\) and \(I=\tau(\mathfrak{s})\), then
\begin{equation}\label{eq:coordinate-strip-bound}
\varepsilon_r\le\omega_{\HH}(x+iy,I)
\le\frac{|I|}{\pi x}\le\frac{C_*}{\pi x},
\end{equation}
i.e., \(0<x\le\frac{C_*}{\pi\varepsilon_r}\).

For the narrowing model define the conformal coordinate \(\tau_o=\tau_U\circ\sigma_R:R_o\to\HH.\) Section~\ref{subsec:subdivision}, together with the finite initial part, provides constants \(c_o,C_o>0\) with
\begin{equation}\label{eq:original-side-lengths}
c_o\le|\tau_o(\mathfrak{s})|\le C_o
\end{equation}
Although \(F\) controls physical distances, it does not directly control the interval lengths under the conformal maps before and after straightening. But the coordinate-strip bound shows that \(\tau_o(R_o\cap T_o(r))\) lies in a vertical strip of uniformly bounded width. We use this control below to compare the conformal coordinates before and after straightening.

\begin{proposition}[Boundary lengths after geometric straightening]\label{prop:straightened-side-lengths}
Let \(R=F(R_o)\) in the narrowing model. For a conformal \(\tau_R:R\to\HH\) taking infinity to infinity, all boundary-side lengths have a uniform positive lower bound and finite upper bound, after a fixed normalization of \(\tau_R\).
\end{proposition}

\begin{proof}
By the construction of \(F\) in Section~\ref{subsec:nearby-edges}, the tapered neighborhoods used to define \(F\) lie in \(T_o(r_F)\) for some fixed \(r_F>0\). Thus the essential support of the Beltrami coefficient \(\mu_F\) of \(F\) lies in this neighborhood. By \eqref{eq:coordinate-strip-bound}, there is a fixed \(B>0\) such that
\begin{equation}\label{eq:original-support-strip}
\tau_o\bigl(T_o(r_F)\cap R_o\bigr)
\subset\{0<\Rea\zeta<B\}.
\end{equation}
Now set 
\[
G=\tau_R\circ F\circ\tau_o^{-1}:\HH\to\HH.
\]
Since \(\tau_R\) and \(\tau_o^{-1}\) are conformal, \eqref{eq:original-support-strip} implies \(\operatorname{supp}\mu_G\subset\{0<\Rea\zeta<B\}.\) Extend \(G\) to the other half-plane by reflection \(G(-\overline{\zeta})=-\overline{G(\zeta)}.\) This is a global quasiconformal homeomorphism, and its Beltrami coefficient is supported in \(|\Rea\zeta|<B\).

By calculation for \(r_\infty>2B\),
\[
\int_{\{\zeta:|\Rea\zeta|<B\}\setminus D(0,r_\infty)}
\frac{dA(\zeta)}{|\zeta|^2}\le
\int_{-B}^{B}\int_{|v|>r_\infty/2}
\frac{dv\,du}{v^2}=\frac{8B}{r_\infty}<\infty.
\]
By the Teichm\"uller--Wittich--Belinskii theorem \cite[Chapter~V, \S6]{LehtoVirtanen}, there exists \(c\in\C\setminus\{0\}\) such that \(G(\zeta)/\zeta\to c\ne0\). The reflection identity gives \(c=\overline c\). Since \(G(\HH)=\HH\) and \(c\ne0\), we have \(c>0\). Rescale \(\tau_R\) by \(c^{-1}\), so that \(G(\zeta)/\zeta\to1\). For \(|\Rea w|\ge2B\), integrating first in \(\Ima\zeta\) gives
\begin{equation*}
\int_{|\Rea\zeta|<B}
\frac{dA(\zeta)}{|\zeta-w|^2}
=\pi\log\frac{|\Rea w|+B}{|\Rea w|-B}
\le\pi\log3.
\end{equation*}
By Lemma~\ref{lem:normalized-distortion}, \(G\) is bilipschitz on \(\{w:|\Rea w|\ge2B\}\).

Fix \(D>2B\). For \(p\in i\R\), reflection and the bilipschitz estimate give
\begin{equation*}
2\Rea G(p+D)
=|G(p+D)-G(p-D)|\asymp D.
\end{equation*}
Since \(G(p)\in i\R\), \(|G(p+D)-G(p)|\ge \Rea G(p+D) \ge cD.\) Quasisymmetry, comparing \(p\) and \(p-D\) from \(p+D\), gives \(|G(p+D)-G(p)|\le C'|G(p+D)-G(p-D)|\le CD.\) Hence \(|G(p+D)-G(p)|\asymp D\) uniformly for \(p\in i\R\).

Let \(q\in i\R\) satisfy \(c_o\le|p-q|\le C_o\). Applying quasisymmetry with base point \(p\) in both orders gives
\[
 \frac{1}{\vartheta_K(D/c_o)}
 \le
 \frac{|G(p)-G(q)|}{|G(p)-G(p+D)|}
 \le
 \vartheta_K(C_o/D).
\]
Since \(D,c_o,C_o\) are fixed and \(|G(p)-G(p+D)|\asymp D\), we obtain \(c'\le|G(p)-G(q)|\le C'\), with constants independent of \(p,q\).
Every side interval has endpoints of this kind by \eqref{eq:original-side-lengths}. The boundary identity \(G\circ\tau_o=\tau_R\circ F\) gives \(G\bigl(\tau_o(\mathfrak{s})\bigr)=\tau_R\bigl(F(\mathfrak{s})\bigr).\) Thus \(G\) maps each original side interval to the corresponding side interval after straightening. This proves the proposition.
\end{proof}

In the logarithmic model put \(\tau_R=\lambda\tau_U\circ\psi^{-1}\) with \(\lambda>0\) fixed; the side bounds follow from Section~\ref{subsec:subdivision}. In the narrowing model use Proposition~\ref{prop:straightened-side-lengths}. A fixed positive rescaling in each case gives
\begin{equation}\label{eq:final-side-lengths}
2\pi\le|\tau_R(\text{boundary side})|\le C_0.
\end{equation}

\section{The corrected folding model and the large-value tract}\label{sec:folded-tract}

The preceding geometry and side lengths provide the graph input for folding. We first compare the variable-scale neighborhoods of the graph with vertical strips in the conformal coordinate.

For every fixed \(r_1>0\), apply Lemma~\ref{lem:uniform-harmonic-side} to the final graph and \eqref{eq:coordinate-strip-bound} with \(\tau=\tau_R\) and \(C_*=C_0\). This gives
\begin{equation}\label{eq:tube-to-strip}
\tau_R\bigl(T(r_1)\cap R\bigr)
\subset\{0<\Rea\zeta<B_{r_1}\}.
\end{equation}
The resulting \(B_{r_1}\) is independent of the edge indices. 

Conversely, for every fixed \(D>0\) there is \(r_D<\infty\) such that
\begin{equation}\label{eq:strip-to-tube}
\tau_R^{-1}\bigl(\{0<\Rea\zeta\le D\}\bigr)
\subset T(r_D).
\end{equation}
To see this, for \(\zeta=x+iy\) with \(0<x\le D\), let \(J=i[y-(D+C_0),y+(D+C_0)]\). By \eqref{eq:poisson-half-plane},
\[
 \omega_{\HH}(\zeta,J)
 =\frac2\pi\arctan\frac{D+C_0}{x}
 \ge\frac12.
\] 
By \eqref{eq:final-side-lengths}, each partition interval has length at least \(2\pi\). Thus at most \(2+\lceil(D+C_0)/\pi\rceil\) partition intervals meet \(J\). Set \(N_D=2+\lceil(D+C_0)/\pi\rceil\). Since \(\omega_{\HH}(\zeta,J)\ge1/2\) and at most \(N_D\) partition intervals meet \(J\), at least one full partition interval \(I\) satisfies \(\omega_{\HH}(\zeta,I)\ge\frac1{2N_D}.\) Let \(z=\tau_R^{-1}(\zeta)\), and let \(\mathfrak{s}\) be the boundary side corresponding to \(I\), with physical edge \(e\). Conformal invariance and the Beurling projection estimate (see \cite[proof of Lemma~2.1]{BishopFolding}) give
\[
 \frac1{2N_D}
 \le\omega_R(z,\mathfrak{s})
 \le C\sqrt{\frac{\diam(e)}{\dist(z,e)}}.
\]
Hence \( \dist(z,e)\le4C^2N_D^2\diam(e).\) Choosing \(r_D>4C^2N_D^2\), independently of \(\zeta\), gives \(z\in T(r_D)\) and proves \eqref{eq:strip-to-tube}.

\subsection{An explicit bounded interpolation}

Let \(U_\pm=\{z:\pm\Ima z>0,\ |z|>1\}\) and set \(\mathcal M(z)=(z-1)/(z+1)\). \(\mathcal M\) maps \(U_+\) to the first quadrant. On the first quadrant define
\[
A_*(re^{i\theta})=re^{i\beta(\theta)}
\quad \text{with} \quad
\beta(\theta)=
\begin{cases}
\theta, & 0\le\theta\le\pi/4,\\
3\theta-\pi/2, & \pi/4\le\theta\le\pi/2.
\end{cases}
\]
It is easy to calculate that \(A_*\) is a \(3\)-quasiconformal map from the first quadrant onto the upper half-plane.

Let \(H_+=\mathcal M^{-1}\circ A_*\circ \mathcal M,\) and \(H_-(z)=\overline{H_+(\overline z)}\). They fix the real boundary outside \([-1,1]\), and
\begin{equation}\label{eq:interpolation-pairing}
    H_+(e^{it})=H_-(e^{-it})\in[-1,1] \quad \text{for} \quad 0\le t\le\pi.
\end{equation}
The endpoints \(\pm1\) are fixed.

The region where \(H_+\) can differ from the identity satisfies
\[
\arg \mathcal M(z)>\pi/4
\quad\Longleftrightarrow\quad
|z|^2-1<2\Ima z.
\]
It lies in the disk centered at \(i\) of radius \(\sqrt2\). Both \(H_\pm\) are the identity for \(|z|\ge1+\sqrt2\) on their domains, and injectivity gives
\begin{equation}\label{eq:interpolation-map-bound}
|H_\pm(z)|\le\max\{|z|,1+\sqrt2\}.
\end{equation}

We now define the interpolation in the final right-half-plane coordinate \(\zeta\) of the folding construction. This coordinate is obtained from \(\tau_R\) by the boundary-partition adjustments and the inverse folding map; these steps are described in the proof of Proposition~\ref{prop:folding} below. Thus \(\zeta\) is the coordinate used to define \(\nu\), rather than the original conformal coordinate \(\tau_R(z)\).

At this stage, subdivide the boundary \(i\R\) at the points \(k\pi i\), \(k\in\Z\), into intervals
\[
 I_k=\{it:k\pi\le t\le(k+1)\pi\}.
\]
Each interval has length \(\pi\), half the imaginary period of the exponential map. Classify \(I_k\) according to the region on the opposite side of its corresponding physical boundary arc. It is type~1 if that region is an \(R\)-component, and type~2 if it is the \(L\)-component.

Assign to each half-strip \(\{\zeta=x+iy:x>0,\ k\pi<y<(k+1)\pi\}\) the type of its boundary interval \(I_k\). Exponentiation maps this half-strip onto \(U_+\) when \(k\) is even and onto \(U_-\) when \(k\) is odd. On \(k\pi<\Ima\zeta<(k+1)\pi\), define
\begin{equation}\label{eq:folding-interpolation}
\nu(\zeta)=
\begin{cases}
e^\zeta, & \text{type 2},\\
H_+(e^\zeta), & \text{type 1 and \(k\) even},\\
H_-(e^\zeta), & \text{type 1 and \(k\) odd}.
\end{cases}
\end{equation}
On \(\{\zeta=x+k\pi i,\ x>0\}\), \(e^\zeta\) is real with modulus greater than one; \(H_\pm\) are the identity in a neighborhood of each such point. The map \(\nu(\zeta)\) therefore is a \(3\)-quasiregular map, with
\begin{equation}\label{eq:interpolation-bounds}
\nu(\zeta)=e^\zeta\ \text{for}\ \Rea\zeta\ge1
\quad \text{and} \quad
\sup_{0<\Rea\zeta\le D}|\nu(\zeta)|
\le\max\{e^D,1+\sqrt2\}.
\end{equation}
The boundary identifications will follow from the simultaneous parameter adjustments in the proof of Proposition~\ref{prop:folding}: type~2 sides will match the \(L\)-side exponential, while type~1 sides will match by \eqref{eq:interpolation-pairing}.

We apply quasiconformal folding in the \(L/R\) form of \cite[Theorem~3.4]{AlbrechtBishop}, using the corrected boundary identifications of \cite[Section~4]{BishopCorrections} and the interpolation \eqref{eq:folding-interpolation}. Solving the associated Beltrami equation then gives an entire function.

\begin{proposition}[Bounded folding]\label{prop:folding}
For either graph, this construction yields a quasiregular map \(g\), a quasiconformal homeomorphism \(\varphi\), and an entire function \(f=g\circ\varphi^{-1}\) such that, for a fixed \(r>0\),
\begin{equation}\label{eq:folded-model}
S(f)\subset\{-1,0,1\},
\quad
g=e^{\tau_R}\text{ on }R\setminus T(r)
\quad\text{and}\quad
\operatorname{supp}\mu_g\subset T(r).
\end{equation}
Furthermore, for every fixed \(r_1>0\),
\begin{equation*}
\sup_{T(r_1)\cup L}|g|<\infty.
\end{equation*}
\end{proposition}

\begin{proof}
We choose the boundary adjustments and folding maps simultaneously on the \(L\)- and \(R\)-sides, as in \cite[Theorem~3.4]{AlbrechtBishop}. The endpoint parities are chosen consistently at every shared vertex.
Write \(it_j\) for the ordered endpoints of the boundary partition induced by \(\tau_R\). First choose a quasiconformal map \(\iota:\HH\to\HH\) that sends these endpoints to \(\pi i m_j\), where the integers \(m_j\) are increasing, and is linear between consecutive boundary endpoints. The choices in \cite[Lemmas~3.1--3.2]{BishopFolding}, with the normalization used here, give uniformly bounded endpoint displacements.

Next choose a quasiconformal map \(\lambda_{\mathrm{adj}}:\HH\to\HH\) fixing these new endpoints. It reparametrizes each boundary interval so that physical arclength is proportional to coordinate length. More precisely, if \(z(t)\), \(0\le t\le\ell\), parametrizes a physical boundary side by arclength, and \(a,b\) are its image endpoints after applying \(\iota\circ\tau_R\), then
\[
 \lambda_{\mathrm{adj}}\bigl(\iota(\tau_R(z(t)))\bigr)
 =a+\frac{t}{\ell}(b-a).
\]
The required extension from the boundary to \(\HH\) is supplied by \cite[Section~4]{BishopFolding}.

Finally, for each resulting boundary interval \(J\), take the square in \(\HH\) having \(J\) as its side on \(i\R\). Inside the union of these squares, insert finite trees attached to \(i\R\), and let \(W\subset\HH\) be the region remaining after these trees removal. The trees and compatible triangulations of \(\HH\) and \(W\) are chosen as in \cite[Sections~5--6]{BishopFolding}. Mapping corresponding triangles affinely gives a quasiconformal homeomorphism
\[
 \psi_F:\HH\to W
\]
with uniformly bounded dilatation. Figure~\ref{fig:bishop-simple-fold} illustrates the folding step when the added tree consists of a single edge. For a point on an added tree edge, approaching from its two sides gives two boundary limits under \(\psi_F^{-1}\). The boundary pairing is chosen so that \(\nu\) takes the same value at these two limits, by \eqref{eq:interpolation-pairing} and \eqref{eq:folding-interpolation}. Hence \(\nu\circ\psi_F^{-1}\) extends continuously across the edge.

\begin{figure}[htbp]
  \centering
  \includegraphics[width=0.5\linewidth]{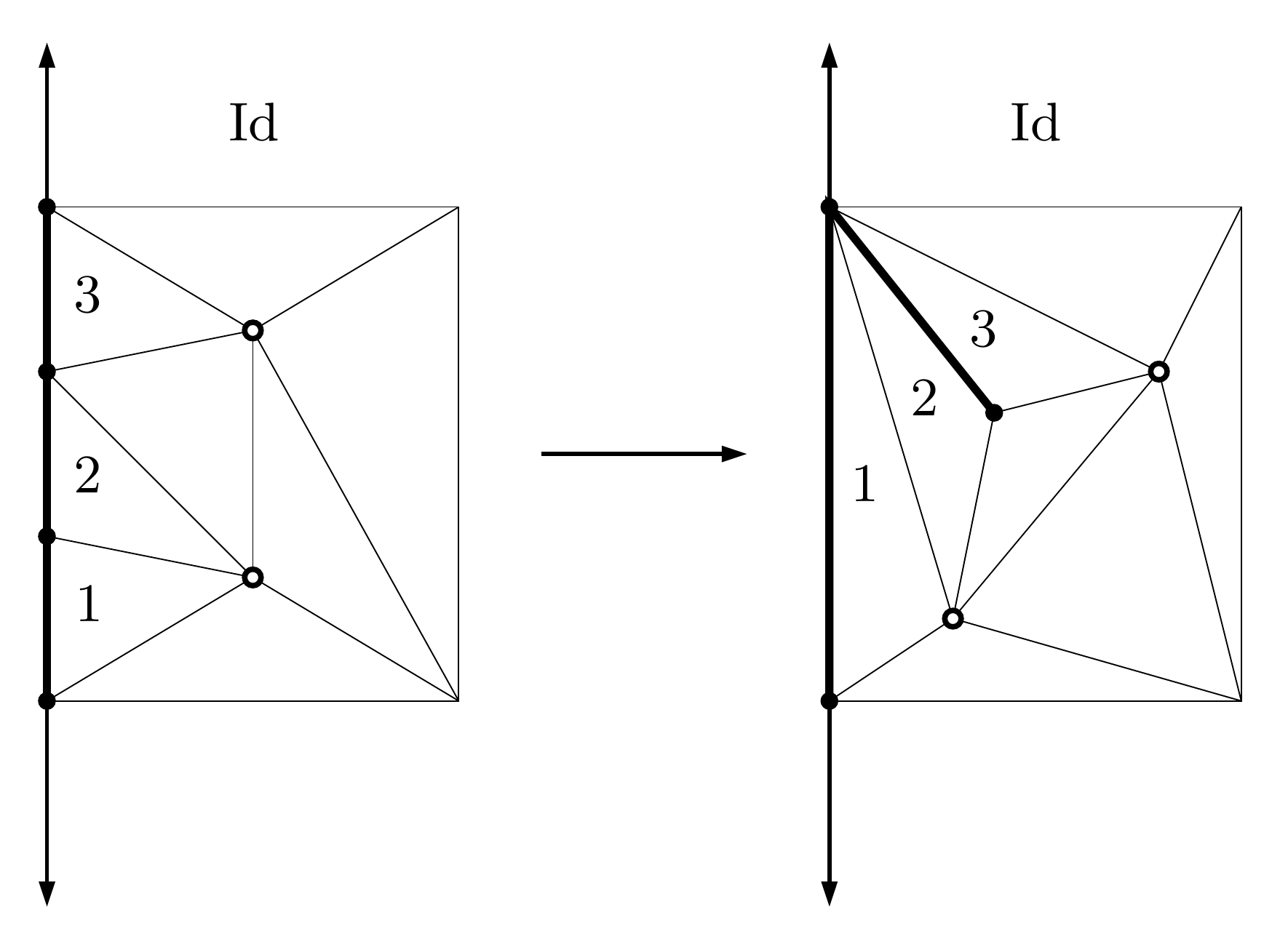}
  \caption[A simple fold along one boundary interval]{%
    A simple fold of a boundary interval $J$ of length $3\pi$.
    The map is affine on each triangle. It maps subinterval~$1$
    onto all of $J$ and subintervals~$2$ and~$3$ onto opposite
    sides of the slit shown in bold.
    Reproduced from Bishop \cite[p.~7, Figure~3]{BishopFolding}.}
  \label{fig:bishop-simple-fold}
\end{figure}

We therefore set
\[
\eta_R=\psi_F^{-1}\circ\lambda_{\mathrm{adj}}\circ\iota\circ\tau_R\quad \text{and} \quad g=\nu\circ\eta_R,
\]
on \( R'=\{z\in R:\lambda_{\mathrm{adj}}(\iota(\tau_R(z)))\in W\}.\) Thus \(\psi_F^{-1}\) is applied on \(W\) and returns the coordinate to \(\HH\), where \(\nu\) is defined. 
Let \(\eta_L:L\to\HHm\) be the adjusted coordinate supplied by the same construction, and set \(g=e^{\eta_L}\) on \(L\). The map \(\eta_L\) is uniformly quasiconformal and agrees with \(\tau_L\) outside a fixed graph neighborhood.

On each common wall edge, the two adjusted coordinates use the same oriented normalized arclength parameter on intervals of length \(2\pi\). Each such edge therefore corresponds to two consecutive type~2 intervals of length \(\pi\), and
\[
 \eta_R(z)-\eta_L(z)\in2\pi i\Z
 \qquad(z\text{ on a common wall edge}).
\]
Hence \(\nu(\eta_R(z))=e^{\eta_L(z)}\) there.
The type~1 sides match by \eqref{eq:interpolation-pairing}, with the paired parameters chosen by the folding construction.
The compatible endpoint parities give the same value \(1\) or \(-1\) at every shared vertex.
Thus the definitions of \(g\) agree across all boundary sides and extend continuously across the graph.

The map \(\iota\) equals the identity whenever \(\Rea\zeta\ge D_\iota\), for some fixed \(D_\iota>0\). The maps \(\lambda_{\mathrm{adj}}\) and \(\psi_F\) equal the identity outside the unions of their corresponding squares. The upper bound in \eqref{eq:final-side-lengths}, together with the uniformly bounded displacement of partition endpoints under \(\iota\), puts all such squares into a strip of fixed width \(D_0\). Each of \(\iota\), \(\lambda_{\mathrm{adj}}\), and \(\psi_F^{-1}\) is injective and fixes \(\{\Rea\zeta\ge D_0\}\). Hence, for \(D'>D_0\), each maps the strip \(\{0<\Rea\zeta<D'\}\) into itself wherever defined.

For a point of \(T(r_1)\cap R\), \eqref{eq:tube-to-strip} bounds its initial real coordinate. Choose \(D'>\max\{D_0,B_{r_1},1\}\). The coordinate changes are injective and equal the identity outside the strip \(\{0<\Rea\zeta<D'\}\). Hence all intermediate coordinates remain in this strip wherever defined, and \eqref{eq:interpolation-bounds} (see also~\eqref{eq:interpolation-map-bound}) gives \(|g|\le\max\{e^{D'},1+\sqrt{2}\}\). Continuity extends the bound across the inserted trees and the original graph. On \(L\), the adjusted model remains exponential in a left half-plane, so \(|g|\le1\). 

Again by \eqref{eq:interpolation-bounds}, \(g=e^{\tau_R}\) whenever \(\Rea\tau_R>D_0\). The remaining region satisfies \( \{z\in R:0<\Rea\tau_R(z)\le D_0\}\subset T(r_{D_0})\) by \eqref{eq:strip-to-tube}. Thus \(g=e^{\tau_R}\) outside this graph neighborhood in \(R\), and its Beltrami coefficient vanishes there. This proves \(g=e^{\tau_R}\) on \(R\setminus T(r)\) and \(\operatorname{supp}\mu_g\subset T(r)\). 

Apply the folding theorem \cite[Theorem~3.4]{AlbrechtBishop}, and solve the Beltrami equation yield \(f=g\circ\varphi^{-1}\) is entire. At each folding vertex, the standard boundary coordinate is \(k\pi i\) for some integer \(k\). Since \(e^{k\pi i}=(-1)^k\) and \(H_\pm\) fix \(\pm1\), the value of \(g\) at such a vertex is \(\pm1\). Thus branching at these vertices can only produce the critical values \(\pm1\). On the \(L\)-side, the exponential model tends to \(0\) as its real coordinate tends to \(-\infty\), accounting for the possible asymptotic value \(0\). Inside each half-strip, \(\nu\) is either the exponential map or its composition with the homeomorphism \(H_+\) or \(H_-\). These maps are locally injective, and the formulas agree with the exponential near the horizontal seams. Hence the interpolation creates no additional interior branch points. After solving the Beltrami equation, the resulting entire function therefore satisfies \(S(f)\subset\{-1,0,1\},\) as asserted in \eqref{eq:folded-model}.
\end{proof}

\subsection{All preimages of sufficiently large values}

By Proposition~\ref{prop:folding}, \(\varphi\) is conformal outside \(T(r)\). The finite area of \(T(r)\) permits the normalization \(\varphi(z)/z\to1\), as recalled in Section~\ref{subsec:qc-lower-bounds}. Proposition~\ref{prop:uniform-log-area} and Lemma~\ref{lem:normalized-distortion} then give a uniform bilipschitz bound for \(\varphi\) on \(\C\setminus T(\Lambda r)\).

Let \(B_1\) bound \(\Rea\tau_R\) on \(T(\Lambda r)\cap R\), and let \(C_1\) bound \(|g|\) on \(T(\Lambda r)\cup L\). Define
\begin{equation*}
M_0=\max\{1,B_1,\log C_1\} \quad\text{and}\quad \Psi:=\tau_R^{-1},\quad V=\Psi\bigl(\{\Rea\zeta>M_0\}\bigr).
\end{equation*}
The set \(V\) is the region selected in the coordinates of the model map \(g\). Under the change of coordinates defining \(f=g\circ\varphi^{-1}\), it corresponds to \(\widehat V:=\varphi(V)\). Since \(V\subset\C\setminus T(\Lambda r)\), the preceding estimates show that \(\varphi\) is conformal and bilipschitz on \(V\).

\begin{lemma}[Exhaustive large-value identity]\label{lem:large-values}
One has
\begin{equation}\label{eq:large-values}
f(\varphi(z))=e^{\tau_R(z)}\text{ for }z\in V\quad \text{and}\quad f^{-1}(A_{e^{M_0}})= \widehat V.
\end{equation}
Consequently \(\widehat V\) is a logarithmic tract and, on a simply connected subset of \(A_{e^{M_0}}\), every inverse branch is
\begin{equation}\label{eq:inversebranches}
h_k(w)=\varphi\bigl(\Psi(\log w+2\pi i k)\bigr)\text{ for }k\in\Z.
\end{equation}
\end{lemma}

\begin{proof}
Since \(V\cap T(\Lambda r)=\varnothing\), proposition~\ref{prop:folding} therefore gives \(g=e^{\tau_R}\) on \(V\). If \(z\notin V\), either \(z\in T(\Lambda r)\cup L\), where \(|g(z)|\le C_1\le e^{M_0}\), or \(z\) lies in the unchanged region of \(R\) with \(\Rea\tau_R(z)\le M_0\), where \(|g(z)|\le e^{M_0}\). Thus \(g^{-1}(A_{e^{M_0}})=V\) and \(f^{-1}(A_{e^{M_0}})=\widehat V\). 

The map \(\tau_R\circ\varphi^{-1}:\widehat V\to\{\Rea\zeta>M_0\}\) is conformal, and \(f=\exp\circ\tau_R\circ\varphi^{-1}\) there. Hence \(f\) maps \(\widehat V\) onto \(A_{e^{M_0}}\) as a universal covering, so \(\widehat V\) is a logarithmic tract. Applying the inverse coordinate map to \(\log w+2\pi ik\) gives all branches in \eqref{eq:inversebranches}. Infinitely many distinct preimages of every large value also show that \(f\) is transcendental.
\end{proof}

\section{The Hausdorff dimension of the escaping set}\label{sec:dimensions}
We complete the proof of Theorem~\ref{thm:construction} by treating $d=1$, $1<d<2$ and $d=2$ separately.

\subsection{The case \texorpdfstring{$d=1$}{d = 1}}\label{sec:dimone}
In this subsection $f$ is the fixed function from the narrowing model. For every $s>1$,~\eqref{eq:physical-cost} gives $A_s<\infty$.

\paragraph{Step 1. Sampling the inverse derivative.}
We first prove that the inverse-derivative sum $P_s(w)$ tends to zero as $|w|\to\infty$. Recall that $\Psi=\tau_R^{-1}:\HH\to R$, where $\HH=\{\Rea\zeta>0\}$. The final graph $T$ consists of the corridor boundary and the subdivided slits after straightening in Section~\ref{subsec:nearby-edges}. Its edges still satisfy~\eqref{eq:physical-cost}, since the straightening is bilipschitz.

For each edge $e$ and each side $\mathfrak{s}$ of $e$ facing $R$, let $I_{e,\mathfrak{s}}=\tau_R(\mathfrak{s})\subset i\R$. Here $s>1$ is the fixed exponent, whereas $\mathfrak{s}$ denotes a boundary side. Boundary values are understood through the side parametrization. Thus $\Psi(I_{e,\mathfrak{s}})=e$ as physical point sets. A slit edge has two sides facing $R$; their coordinate intervals have disjoint interiors, but both map to the same physical edge. These side intervals partition $i\R$, and~\eqref{eq:final-side-lengths} gives $2\pi\le |I_{e,\mathfrak{s}}|\le C_0$.

\emph{Sampling on $\Rea\zeta=1$.}
For every $y\in\R$, the Poisson kernel gives
\[
 \omega_{\HH}\bigl(1+iy,i[y-2,y+2]\bigr)
 =\frac1\pi\int_{-2}^{2}\frac{dt}{1+t^2}
 =\frac{2}{\pi}\arctan2.
\]
The window has length $4<2\pi$, so it meets at most two side intervals in positive length. At least one of these full intervals therefore satisfies
\[
 \omega_{\HH}(1+iy,I_{e,\mathfrak{s}})\ge\frac{\arctan2}{\pi}.
\]
The full interval need not be contained in the window. Set $z=\Psi(1+iy)\in R$. Conformal invariance and the Beurling estimate~\eqref{eq:beurling-side} give
\[
 \frac{\arctan2}{\pi}\le\omega_R(z,\mathfrak{s})
 \le\frac4\pi\sqrt{\frac{\diam e}{\dist(z,e)}}.
\]
Hence $\dist(z,\partial R)\le\dist(z,e)\le C\diam e$. Koebe's estimate~\eqref{eq:koebe-distance}, with $X=1$, now gives $|\Psi'(1+iy)|\le C\diam e$. Since the selected edge depends on $y$, we record the uniform bound
\[
 |\Psi'(1+iy)|^s\le C_s
 \sum_{(e,\mathfrak{s}):\,I_{e,\mathfrak{s}}\cap i[y-2,y+2]\ne\varnothing}
 (\diam e)^s\quad\text{for}\quad y\in\R.
\]

Fix $\theta\in\R$ and sample at $y=\theta+2\pi k$, $k\in\Z$. If $I_{e,\mathfrak{s}}=i[a,b]$, it meets the sampling window only when $\theta+2\pi k\in[a-2,b+2]$. This interval has length at most $C_0+4$. Thus each side interval occurs at most $1+\lfloor(C_0+4)/(2\pi)\rfloor$ times, independently of $\theta$. Summing the preceding bound and using the fact that each geometric edge has at most two sides facing $R$, we obtain
\begin{equation}\label{eq:one-sampling}
 \sup_{\theta\in\R}\sum_{k\in\Z}
 |\Psi'(1+i(\theta+2\pi k))|^s\le C_sA_s.
\end{equation}
The lower bound on side lengths controls the number of intervals meeting one window; the upper bound controls the number of samples associated with one interval.

\emph{Moving to $\Rea\zeta=x$.}
For every $x\ge1$ and $y\in\R$, we claim that
\begin{equation}\label{eq:halfplane-distortion}
 |\Psi'(x+iy)|\le x|\Psi'(1+iy)|.
\end{equation}
To apply disk distortion, set
\[
 A_y(\zeta)=\frac{1+\zeta}{1-\zeta}+iy,
 \quad F_y=\Psi\circ A_y\quad\text{and}\quad r=\frac{x-1}{x+1}\in[0,1).
\]
The map $A_y$ sends $\mathbb D$ conformally onto $\HH$, with $A_y(0)=1+iy$ and $A_y(r)=x+iy$. The chain rule gives
\[
 |F_y'(0)|=2|\Psi'(1+iy)|\quad\text{and}\quad |F_y'(r)|=\frac{2}{(1-r)^2}|\Psi'(x+iy)|.
\]
Since $F_y$ is univalent, Koebe's derivative distortion theorem yields
\[
 \frac{2}{(1-r)^2}|\Psi'(x+iy)|
 \le\frac{1+r}{(1-r)^3}\,2|\Psi'(1+iy)|.
\]
After cancellation, the factor on the right is $(1+r)/(1-r)=x$. This proves~\eqref{eq:halfplane-distortion}.

\emph{The inverse branches of $f$.}
Let $w$ be a target value with $|w|>e^{M_0}$. Choose a local branch of $\log$ near $w$ and put $x=\log|w|>M_0\ge1$ and $\theta=\Ima\log w$. By Lemma~\ref{lem:large-values}, all inverse branches are
\[
 h_k(w)=\varphi\bigl(\Psi(\log w+2\pi ik)\bigr)\quad\text{and}\quad \log w+2\pi ik=x+i(\theta+2\pi k),\quad k\in\Z.
\]
These logarithmic points lie in $\{\Rea\zeta>M_0\}$, so their images under $\Psi$ belong to $V$. The map $\varphi$ is conformal and Lipschitz there. If $L_0$ is its Lipschitz constant on $V$, then $|\varphi'|\le L_0$ on $V$, and the chain rule gives
\[
 h_k'(w)=
 \frac{\varphi'\bigl(\Psi(\log w+2\pi ik)\bigr)
       \Psi'(\log w+2\pi ik)}{w}.
\]
Applying~\eqref{eq:halfplane-distortion} and then~\eqref{eq:one-sampling}, we obtain
\begin{equation}\label{eq:onepressure}
\begin{aligned}
 P_s(w)=\sum_{k\in\Z}|h_k'(w)|^s
 &\le\left(\frac{L_0x}{|w|}\right)^s
      \sum_{k\in\Z}|\Psi'(1+i(\theta+2\pi k))|^s\\
 &\le C_sL_0^sA_s\left(\frac{\log|w|}{|w|}\right)^s
 \longrightarrow0\text{ as }|w|\to\infty.
\end{aligned}
\end{equation}
All constants are independent of $\theta$, so the convergence is uniform in $\arg w$. It remains to supply the finite-cost initial cover required by Theorem~\ref{prop:covering-principle}.

\paragraph{Step 2. The initial cover and the dimension conclusion.}
Fix $s>1$ and a desired mesh $\varepsilon>0$, where the mesh is the supremum of the diameters of the covering sets. We first construct a cover of $\Omega$ with finite $s$-cost, then transfer it to $\widehat V$. Enlarge the Lipschitz constant $L_0$ above, if necessary, so that $L_0\ge1$.

Recall that $\gamma(t)=te^{it}$ and $\eta(t)=e^{-t}$. The center curve has arclength element $|\gamma'(t)|\,dt=\sqrt{1+t^2}\,dt\asymp t\,dt$. Along the far corridor, place disks with diameters comparable to the local width $\eta(t)$ and with centers spaced by a comparable arclength. Each disk covers a center-curve arc of length comparable to $\eta(t)$. Along each arc, the parameter changes by $O(\eta(t)/t)$. Hence $\eta$ stays comparable to its value at the center. We need $O(t/\eta(t))$ disks per unit increase in $t$. Each disk contributes $O(\eta(t)^s)$ to the $s$-cost. The contribution over a parameter increment $dt$ is therefore bounded by a constant times $(t\,dt/\eta(t))\eta(t)^s=t\eta(t)^{s-1}\,dt$.

Choose $t_0>T_0$ so large that all these disk diameters are at most $\varepsilon/L_0$ for $t\ge t_0$. If $D_j$ are the disks covering this tail, then
\begin{equation*}
\begin{aligned}
 \sum_j(\diam D_j)^s
 &\le C_s\int_{t_0}^{\infty}\frac{t}{\eta(t)}\eta(t)^s\,dt
 =C_s\int_{t_0}^{\infty}t\eta(t)^{s-1}\,dt\\
 &=C_s\int_{t_0}^{\infty}te^{-(s-1)t}\,dt<\infty\quad\text{for}\quad s>1.
\end{aligned}
\end{equation*}
The remaining part of $\Omega$, including the cap, is bounded. Cover it by finitely many disks of diameter at most $\varepsilon/L_0$. Together these disks form a countable cover $\{D_j\}$ of $\Omega$ with mesh at most $\varepsilon/L_0$ and finite $s$-cost.

The geometric straightening maps $\Omega$ onto itself, so $V\subset R\subset\Omega$. Thus the nonempty sets $E_j=\varphi(D_j\cap V)$ cover $\widehat V=\varphi(V)$. Since $\varphi$ is $L_0$-Lipschitz on $V$, we have
\[
 \diam E_j\le L_0\diam D_j\le\varepsilon\quad\text{and}\quad
 \sum_j(\diam E_j)^s\le L_0^s\sum_j(\diam D_j)^s<\infty.
\]
For every $R_\ast>e^{M_0}$ and $z\in X_{R_\ast}$, the definition~\eqref{eq:exterior-orbits} gives $|f(z)|>R_\ast>e^{M_0}$. Hence~\eqref{eq:large-values} implies $z\in\widehat V$. The same sets $E_j$ therefore cover $X_{R_\ast}$, and choosing $\varepsilon\le1/16$ supplies the initial cover required by Theorem~\ref{prop:covering-principle}.

By~\eqref{eq:onepressure}, we may choose $R_\ast$ sufficiently large that the successive inverse-image covers in Theorem~\ref{prop:covering-principle} have both mesh and total $s$-cost tending to zero. The theorem then gives $\dimH I(f)\le s$. The graph and the function $f$ are fixed for all $s>1$; the cutoff $R_\ast$ may depend on $s$. Letting $s\downarrow1$ gives $\dimH I(f)\le1$. Combining this with the universal lower bound~\cite[Theorem~9.1]{BergweilerRempe}, we obtain
\begin{equation*}
 \dimH I(f)=1.
\end{equation*}

\subsection{The case \texorpdfstring{$1<d<2$}{1 < d < 2}}\label{sec:intermediate}\label{sec:blocks}
Fix $1<d<2$, put $q=(d-1)^{-1}>1$, and use the logarithmic model with parameters~\eqref{eq:params}. We first establish the boundary-block estimates, then prove the upper and lower dimension bounds for the resulting fixed function. For $s>1+1/(q+1)$ and $\varepsilon>0$, the goal is
\begin{equation}\label{eq:targetsampling}
 \sup_{\theta\in\R}\sum_{k\in\Z}
 |\Psi'(X+i(\theta+2\pi k))|^s\le C_{s,\varepsilon}X^{1-s+\varepsilon}
 \quad\text{for}\quad X\ge2.
\end{equation}
This sharper estimate uses whole columns rather than individual small edges. The scaling of $\tau_R$ has already been fixed and is absorbed into comparison constants.

\paragraph{Inverse derivative along one slit.}
\begin{lemma}[Root and tip control]\label{lem:root-tip}
On either side of a tail slit, let $t\ge0$ be conformal distance along $i\R$ from that side's root value toward the tip. For $t>0$ on this side,
\begin{equation}\label{eq:root-tip}
 \left|\frac{d\Psi}{dt}\right|\le
 \frac{C\Delta_n}{\sqrt{t(1+t)}}.
\end{equation}
\end{lemma}
\begin{proof}
Fix one side of an upper slit and put $z(y)=x_n+i(\pi/2-y)$ for $0\le y\le h_n$. Here $y$ is the depth from the wall. All boundary values and derivatives are taken from this side. Since $\tau_R=\lambda\tau_U\circ\psi^{-1}$ with fixed $\lambda>0$, monotonicity of the boundary image gives the exact expression
\[
 t(y)=\lambda|\tau_U(z(y))-\tau_U(z(0))|
 =\lambda\int_0^y|\tau_U'(z(v))|\,dv.
\]
Thus $dt/dy=\lambda|\tau_U'(z(y))|$ for $0<y<h_n$. Differentiating $\Psi(\lambda\tau_U(z(y)))=\psi(z(y))$ along this side and using~\eqref{eq:log-width}, we obtain
\[
 \left|\frac{d\Psi}{dt}\right|
 =\frac{|\psi'(z(y))|}{\lambda|\tau_U'(z(y))|}
 \asymp\frac{\eta(x_n)}{|\tau_U'(z(y))|}.
\]

For $0<y\le h_n-4\delta_n$,~\eqref{eq:side-der} and the integral expression for $t$ give
\[
 t(y)\asymp\cosh\frac{\pi y}{\delta_n}-1\quad\text{and}\quad
 |\tau_U'(z(y))|\asymp\frac1{\delta_n}\sinh\frac{\pi y}{\delta_n}.
\]
Using $\sinh^2 v=(\cosh v-1)(\cosh v+1)$ and $\eta(x_n)\delta_n\asymp\Delta_n$ from~\eqref{eq:log-small}, we get
\[
 \left|\frac{d\Psi}{dt}\right|
 \asymp\frac{\eta(x_n)\delta_n}{\sinh(\pi y/\delta_n)}
 \asymp\frac{\Delta_n}{\sqrt{t(1+t)}}.
\]

For $0<u=h_n-y\le4\delta_n$,~\eqref{eq:tip-der} and~\eqref{eq:transition-derivative} together give $|\tau_U'(z(y))|\ge cM_n/\delta_n$. The chain rule therefore yields $|d\Psi/dt|\le C\Delta_n/M_n$. By Lemma~\ref{lem:column-estimates}\,(1), the distance $t$ is at most the full column length: $t\le\lambda\ell_n\le CM_n$. Since $M_n\to\infty$ as $n\to\infty$, we have $\sqrt{t(1+t)}\le CM_n$, and hence
\[
 \left|\frac{d\Psi}{dt}\right|
 \le C\frac{\Delta_n}{M_n}
 \le\frac{C'\Delta_n}{\sqrt{t(1+t)}}.
\]
The estimate at the tip follows by continuity. This proves~\eqref{eq:root-tip}. The same estimates hold on the other side, with $t$ measured anew from that side's root. Lower slits are treated by reflection.
\end{proof}

\paragraph{The cost of boundary blocks.}
For $X\ge2$ partition $i\R$ into $J_j=i[jX,(j+1)X]$, $j\in\Z$, and define
\begin{equation*}
 B_s(X)=\sum_{j\in\Z}(\diam\Psi(J_j))^s.
\end{equation*}
Boundary values are interpreted through the side parametrization.
Let $K_n$ be the interval of the $n$th whole column. From~\eqref{eq:column-scale}, with the final scaling absorbed,
\begin{equation}\label{eq:Kn}
 |K_n|\asymp M_n\quad\text{and}\quad|\Ima\zeta|\asymp e^{x_n}\quad\text{for}\quad\zeta\in K_n.
\end{equation}

\begin{lemma}[Uniform block cost]\label{lem:block-cost}
For every $s>1+1/(q+1)$ and $\varepsilon>0$,
\begin{equation}\label{eq:block-cost}
 B_s(X)\le C_{s,\varepsilon}X^\varepsilon\quad\text{for}\quad X\ge2.
\end{equation}
\end{lemma}
\begin{proof}
Call a column \emph{long} if $|K_n|\ge X$ and \emph{short} otherwise.

\emph{Long columns.}
By Lemma~\ref{lem:root-tip} and direct integration, the part of either slit side corresponding to $0\le t\le2X$ has physical length at most
\[
 C\Delta_n\int_0^{2X}\frac{dt}{\sqrt{t(1+t)}}
 \le C\Delta_n\log(2+X).
\]
This part meets only a bounded number of grid blocks. For the remaining pieces in the ranges $kX\le t\le(k+1)X$, $k\ge2$, the same lemma gives physical length at most $C\Delta_n/k$. Each range also meets only a bounded number of blocks. The following trunk edge has physical length $O(\Delta_n)$ and conformal length $O(1)$. Splitting at the tip and the range endpoints gives only a bounded number of pieces in each block. Applying these bounds on both sides and adding the trunk contribution, we obtain
\begin{equation*}
\begin{aligned}
 \sum_j\bigl(\diam\Psi(J_j\cap K_n)\bigr)^s
 &\le C_s\Delta_n^s
 \left((\log(2+X))^s+\sum_{k\ge2}k^{-s}+1\right)\\
 &\le C_s\Delta_n^s(\log(2+X))^s.
\end{aligned}
\end{equation*}
Since $\sum_n\Delta_n^s<\infty$ by~\eqref{eq:log-mesh-decay}, the total contribution from long columns is at most $C_s(\log(2+X))^s$.

\emph{Short columns.}
We bound the number of grid blocks meeting short columns and the physical size of all these columns. By~\eqref{eq:Kn}, a short column satisfies $e^{x_n}\delta_n=M_n\le CX$. The estimate $\log\delta_n\ge-Cx_n^{1/q}$ from~\eqref{eq:log-mesh-decay} therefore gives $x_n-Cx_n^{1/q}\le\log(CX)$. Since $q>1$, this first implies $x_n=O(\log(2+X))$. Substituting this bound back into the same inequality yields
\begin{equation*}
 x_n\le\log X+C(\log(2+X))^{1/q}\quad\text{for}\quad |K_n|<X\text{ and } X\ge2.
\end{equation*}
By~\eqref{eq:Kn}, all short columns lie in the interval
\[
 |\Ima\zeta|\le CX\exp\!\left(C(\log(2+X))^{1/q}\right).
\]
They therefore meet at most $C\exp(C(\log(2+X))^{1/q})$ grid blocks, since each block has length $X$. By~\eqref{eq:log-radius}, their physical images lie in a disk of radius $C\exp(C(\log(2+X))^{1/q})$. This also covers the following trunk edges, where $x\in[x_n,x_{n+1}]$ and $x_{n+1}-x_n=\delta_n\to0$ as $n\to\infty$.

The initial boundary portion has bounded conformal coordinates and a bounded physical image. Thus the same bounds apply to this portion. We multiply the number of blocks by the $s$th power of the diameter bound. The total $s$-cost of the short columns and the initial portion is therefore at most
\begin{equation*}
 C_s\exp\!\left(C_s(\log(2+X))^{1/q}\right)
 \le C_{s,\varepsilon}X^\varepsilon.
\end{equation*}
The last inequality follows from $q>1$.

Combining the above estimates, we obtain
\[
 B_s(X)\le C_s(\log(2+X))^s
       +C_s\exp\!\left(C_s(\log(2+X))^{1/q}\right)\le C_{s,\varepsilon}X^\varepsilon\text{ for } X\ge2.
\]
This proves~\eqref{eq:block-cost}.
\end{proof}

\paragraph{From blocks to derivatives.}
\begin{lemma}[Uniform lattice sampling]\label{lem:block-sampling}
Estimate~\eqref{eq:targetsampling} holds.
\end{lemma}
\begin{proof}
The argument is similar to the boundary sampling argument in Section~\ref{sec:dimone}.
At $X+iy$, the interval $i[y-X,y+X]$ has harmonic measure $1/2$. It meets at most four grid blocks; their physical images join at endpoints. The Beurling estimate and Koebe's theorem~\eqref{eq:koebe-distance} give
\[
 |\Psi'(X+iy)|^s\le C_sX^{-s}
 \sum_{J_j\cap i[y-X,y+X]\ne\varnothing}(\diam\Psi(J_j))^s.
\]
A fixed block is seen only for $y\in[jX-X,(j+1)X+X]$, an interval of length $3X$. A lattice of spacing $2\pi$ with any phase $\theta$ has at most $CX$ samples there. Summing at $y=\theta+2\pi k$ gives
\[
 \sum_k|\Psi'(X+i(\theta+2\pi k))|^s
 \le C_sX^{1-s}B_s(X)\le C_{s,\varepsilon}X^{1-s+\varepsilon}.
\]
The counting constant is independent of the phase.
\end{proof}
The power $X^{-s}$ converts a boundary diameter into a derivative. The additional factor $X$ counts the lattice samples seeing one block. Together they explain $X^{1-s}$.

\paragraph{Part 1. The upper bound.}
We prove that $\dimH I(f)\le d$ by constructing a cover with finite cost and then applying Theorem~\ref{prop:covering-principle}. Fix $s>d=1+1/q$. The boundary estimates above apply because $s>1+1/(q+1)$. Choose
\begin{equation}\label{eq:epsilon-choice}
 0<\varepsilon<q(s-1)-1.
\end{equation}

\emph{The inverse-derivative sum.}
Let $|w|>e^{\max\{M_0,2\}}$ and choose a local branch of $\log$ near $w$. By Lemma~\ref{lem:large-values}, all inverse branches are $h_k(w)=\varphi(\Psi(\log w+2\pi ik))$, $k\in\Z$. The chain rule gives
\[
 h_k'(w)=\varphi'\!\left(\Psi(\log w+2\pi ik)\right)
          \Psi'(\log w+2\pi ik)\frac1w.
\]
Since $\Rea(\log w+2\pi ik)>M_0$, the argument of $\varphi'$ lies in $V$. The Lipschitz bound for $\varphi$ on $V$ therefore gives a uniform bound for this factor. Apply Lemma~\ref{lem:block-sampling} with $X=\log|w|$ and $\theta=\Ima\log w$. Recalling that $P_s(w)=\sum_k|h_k'(w)|^s$, we obtain
\begin{equation}\label{eq:logpressure}
 P_s(w)\le C_{s,\varepsilon}|w|^{-s}(\log|w|)^{1-s+\varepsilon}
 \longrightarrow0\text{ as }|w|\to\infty.
\end{equation}
The factor $|w|^{-s}$ comes from $(\log w)'=1/w$. To apply Theorem~\ref{prop:covering-principle}, we still need an initial cover of small mesh and finite $s$-cost.

\emph{A cover obtained by one inverse image.}\label{subsec:firstpullback}
On the center curve,~\eqref{eq:log-model} gives $\arg z=a^{-1}\log|z|$ for a continuous choice of argument. Thus, with $r=|z|$, this curve is parametrized by $z(r)=r\exp(i\log r/a)$ and satisfies $|z'(r)|=\sqrt{1+a^{-2}}$. For every sufficiently large integer $m$, its length in $e^m\le|z|\le e^{m+1}$ is $O(e^m)$. The corridor width there is comparable to $e^m/m^{q-1}$. We can therefore cover this part of the corridor by $O(m^{q-1})$ disks whose radii are comparable to the width.

Intersect these disks with $V$, apply $\varphi$, and cover each nonempty image by a disk centered in that image. The Lipschitz bound on $V$ and $\varphi(z)/z\to1$ as $|z|\to\infty$ give $O(m^{q-1})$ target disks covering the corresponding part of $\widehat V$. Their centers and radii satisfy $|w|\asymp e^m$ and $R_m\asymp e^m/m^{q-1}$. Hence $R_m/|w|\asymp m^{1-q}\to0$ as $m\to\infty$, uniformly over these disks. For large $m$, every doubled disk lies in $\{|w|>e^{M_0}\}$, where all inverse branches extend. Koebe distortion gives $\diam h_k(D)\le C R_m|h_k'(w)|$ for $D=D(w,R_m)$. Together with~\eqref{eq:logpressure}, this yields
\begin{equation}\label{eq:firstdisk}
 \sum_k(\diam h_k(D))^s
 \le C_sR_m^sP_s(w)
 \le C_{s,\varepsilon}
 \left(\frac{e^m}{m^{q-1}}\right)^s e^{-ms}m^{1-s+\varepsilon}
 =C_{s,\varepsilon}m^{1-qs+\varepsilon}.
\end{equation}
Multiplying by the $O(m^{q-1})$ target disks in each annulus, we obtain the finite total cost
\begin{equation}\label{eq:initialseries}
 \sum_m C_{s,\varepsilon}m^{q-1}m^{1-qs+\varepsilon}
 =\sum_m C_{s,\varepsilon}m^{q(1-s)+\varepsilon}<\infty.
\end{equation}
Indeed,~\eqref{eq:epsilon-choice} gives $q(1-s)+\varepsilon<-1$. It also gives $1-qs+\varepsilon<-q<0$, so~\eqref{eq:firstdisk} shows that every individual inverse-image diameter tends uniformly to zero as $m\to\infty$.

Fix $R_\ast>e^{\max\{M_0,2\}}$ and a small mesh. The preceding estimates handle a sufficiently far tail. Cover the remaining bounded part of $\widehat V\cap\{|w|>R_\ast\}$ by finitely many disks $D=D(w,r)$, with $D(w,2r)\subset\{|w|>e^{M_0}\}$. By Koebe distortion and~\eqref{eq:logpressure},
\[
 \sum_k\bigl(\diam h_k(D)\bigr)^s
 \le C_s r^s P_s(w)\le C_s' r^s.
\]
Choosing the radii sufficiently small gives the prescribed mesh for all inverse images, with finite total $s$-cost.

\emph{Applying the covering principle.}
Recall that $X_{R_\ast}=\{z:|f^n(z)|>R_\ast\text{ for every }n\ge0\}$. If $R_\ast>e^{M_0}$, Lemma~\ref{lem:large-values} gives $X_{R_\ast}\subset\widehat V\cap\{|w|>R_\ast\}$. Since $f(X_{R_\ast})\subset X_{R_\ast}$, we have
\[
 X_{R_\ast}\subset
 f^{-1}\!\left(\widehat V\cap\{|w|>R_\ast\}\right).
\]
The preceding cover therefore covers $X_{R_\ast}$. Together with~\eqref{eq:logpressure}, this verifies both hypotheses of Theorem~\ref{prop:covering-principle} and gives $\dimH I(f)\le s$. The function $f$ is fixed while $s>d$ varies. Letting $s\downarrow d$ yields
\begin{equation}\label{eq:upperd}
 \dimH I(f)\le d.
\end{equation}
The threshold comes from~\eqref{eq:initialseries}: one can choose $\varepsilon>0$ to make the series converge precisely when $q(s-1)>1$, or $s>d$.

\paragraph{Part 2. The lower bound.}
We estimate the growth of $f$ and apply the lower bound~\eqref{eq:tract-lower-bound}. For a continuous map $u$ on the plane, write $M(r,u)=\max_{|z|=r}|u(z)|$. We use this notation for both $f$ and the quasiregular model $g$.

In the unchanged part of $R$,~\eqref{eq:Qglobal} gives $\log|g(z)|=\Rea\tau_R(z)\le Ce^x$, where $x=\Rea\psi^{-1}(z)$. By the logarithmic model, $x\asymp(\log|z|)^q$ as $|z|\to\infty$ in the corridor. On the remaining region, $g$ is bounded by Proposition~\ref{prop:folding}. Consequently,
\[
 \log\log M(r,g)\le C(\log r)^q
 \quad\text{for}\quad r\text{ sufficiently large}.
\]
Since $\varphi(z)/z\to1$ as $|z|\to\infty$, its inverse also satisfies $\varphi^{-1}(z)/z\to1$. In particular, $|\varphi^{-1}(z)|\le2|z|$ for large $|z|$. The identity $f=g\circ\varphi^{-1}$ from~\eqref{eq:folding-factorization} therefore gives
\begin{equation*}
 \log\log M(r,f)\le C'(\log r)^q
 \quad\text{for}\quad r\text{ sufficiently large}.
\end{equation*}
For every $\varepsilon>0$, we have $C'(\log r)^q\le(\log r)^{q+\varepsilon}$ for sufficiently large $r$. Exponentiating twice verifies~\eqref{eq:tract-growth}.

By Lemma~\ref{lem:large-values}, $\widehat V$ is a logarithmic tract. Its boundary is a regular level curve of $|f|$, since $e^{M_0}$ is larger than the moduli of all singular values. The entire function $f$ is continuous on its closure and satisfies $|f|=e^{M_0}$ on the boundary. Its complement is unbounded because it contains $\varphi(L)$. Thus the tract hypotheses of~\cite[Theorem~1.2]{BergweilerKarpinskaStallard} hold. Applying~\eqref{eq:tract-lower-bound} with $D=\widehat V$, we obtain
\begin{equation*}
 \dimH I(f)\ge\dimH I(f,\widehat V)\ge1+\frac1q=d.
\end{equation*}
Together with~\eqref{eq:upperd}, this proves $\dimH I(f)=d$ for the intermediate-dimensional construction.

\subsection{The case \texorpdfstring{$d=2$}{d = 2}}\label{sec:scope}
For $d=2$, take $f(z)=e^z$. Then $S(f)=\{0\}$, so $f\in\calS$, and $\dimH I(f)=2$ by McMullen's result~\cite[proof of Theorem~1.2]{McMullen}.
Together with the constructions in Sections~\ref{sec:dimone} and~\ref{sec:intermediate}, this proves Theorem~\ref{thm:construction} for every $d\in[1,2]$.

\section*{Acknowledgments}

The main proof strategies were proposed by the authors. We used ChatGPT (OpenAI) to simply some calculations and polish partial text and figures. The authors independently verified all results and take full responsibility for the mathematical content.

\phantomsection

\end{document}